\documentclass[10pt,twoside,a4paper,reqno,final]{amsart}
\usepackage[latin1]{inputenc}
\usepackage[T1]{fontenc}

\usepackage[centertags]{amsmath}
\usepackage{amsthm}
\usepackage{amssymb}
\usepackage{amsfonts,a4}
\usepackage{amsxtra}
\usepackage{amscd}
\usepackage{mathabx} 
\usepackage{esint}

\usepackage[english]{babel}

\usepackage{algorithm}
\usepackage{algorithmic}

\usepackage{color}
\usepackage[mathscr]{eucal}

\usepackage{bbm}
\usepackage{enumitem}

\usepackage{graphicx}

\usepackage{xspace} %

\DeclareMathAlphabet{\mathpzc}{OT1}{pzc}{m}{it}

\newtheorem{thrm}{Theorem}
\newtheorem{crllr}[thrm]{Corollary}
\newtheorem{lmm}[thrm]{Lemma}
\newtheorem{prpstn}[thrm]{Proposition}

\newtheorem{rmrk}[thrm]{Remark}
\theoremstyle{definition}
\newtheorem{ass}[thrm]{Assumption}
\newtheorem{xmpl}[thrm]{Example}

\numberwithin{equation}{section}

\newcommand{\Ax}{\mathcal{A}}

\newcommand{\Ao}{\ensuremath{\mathcal{A}}}
\newcommand{\ao}{\ensuremath{\mathfrak{a}}}
\newcommand{\Aon}[1][n]{\ensuremath{\Ao_{n}}}

\newcommand{\abs}[1]{\ensuremath{\left|#1\right|}}
\newcommand{\asdefined}{\mathrel{=:}}

\newcommand{\bsS}{\ensuremath{\boldsymbol{S}}}
\newcommand{\bsT}{\ensuremath{\boldsymbol{T}}}
\newcommand{\bSn}[1][n]{\ensuremath{\boldsymbol{S}^{#1}}}
\newcommand{\bSnG}[1][n]{\ensuremath{\boldsymbol{S}^{#1}_{\grid}}}
\newcommand{\bSk}[1][k]{\ensuremath{\boldsymbol{S}_{{#1}}}}
\newcommand{\bsD}{\ensuremath{\boldsymbol{D}}}

\newcommand{\bbI}{\ensuremath{\boldsymbol{1}}}

\newcommand{\Bogn}[1][n]{\mathfrak{B}^{#1}}
\newcommand{\BogG}[1][\grid]{\mathfrak{B}^{#1}}
\newcommand{\Bogk}[1][k]{\mathfrak{B}^{#1}}
\newcommand{\bweak}{\ensuremath{\overset{b}{\rightharpoonup}}}

\newcommand{\Cleq}{\ensuremath{\lesssim}}

\renewcommand{\d}{\ensuremath{{\rm d}}}

\newcommand{\definedas}{\mathrel{:=}}

\newcommand{\dual}[2]{\ensuremath{\left\langle #1,\,#2\right\rangle}}

\newcommand{\dx}{\,\mathrm{d}x}
\newcommand{\ds}{\,\mathrm{d}s}

\DeclareMathOperator{\divo}{div}

\DeclareMathAlphabet{\lf}{OT1}{pzc}{m}{it}

\newcommand{\eg}{e.\,g.\xspace}
\newcommand{\elm}{\ensuremath{E}\xspace}

\newcommand{\est}{\mathcal{E}}

\newcommand{\frakB}{\operatorname{B}}

\newcommand{\grids}{\ensuremath{\mathbb{G}}\xspace}
\newcommand{\grid}{\mathcal{G}}
\newcommand{\gridn}[1][n]{\grid_{#1}}
\newcommand{\gridk}[1][k]{\grid_{#1}}

\newcommand{\helm}[1][\ell]{\ensuremath{h_{\elm}}}

\newcommand{\hG}[1][\grid]{\ensuremath{h_{#1}}}

\newcommand{\PG}{\ensuremath{\mathfrak{I}}}
\newcommand{\PGdiv}[1][\grid]{\ensuremath{\PG^{#1}_{\divo}}}
\newcommand{\PGdivk}[1][k]{\ensuremath{{\PG}^{\grid_{#1}}_{\divo}}}
\newcommand{\PGQ}[1][\grid]{\ensuremath{{\PG}^{#1}_\Q}}

\newcommand{\id}{\operatorname{id}}

\newcommand{\jump}[1]{\left[\negthinspace\left[{#1}\right]\negthinspace\right]}

\newcommand{\mat}[1]{\boldsymbol{#1}}
\newcommand{\mcD}{\mathcal{D}}

\newcommand{\marked}{\ensuremath{\mathcal{M}}}

\newcommand{\MARK}{\texttt{MARK}\xspace}

\newcommand{\Mk}[1][k]{\marked_{#1}}

\newcommand{\nablas}{\operatorname{D}\hspace{-1pt}}

\newcommand{\N}{\ensuremath{\mathbb{N}}}

\newcommand{\norm}[2][\Omega]{\ensuremath{\left\|#2\right\|_{#1}}}
\newcommand{\normbig}[2][\Omega]{\ensuremath{\big\|#2\big\|_{#1}}}

\DeclareMathOperator{\osc}{osc}

\renewcommand{\P}{\ensuremath{\mathbb{P}}}

\renewcommand{\paragraph}[1]{\noindent\raisebox{0pt}[10pt][0pt]{\textbf{#1.}}}

\newcommand{\PGS}[1][\grid]{\ensuremath{\Pi_{#1}}}
\newcommand{\pzb}{\mathpzc{b}}
\newcommand{\PnG}[1][n]{\ensuremath{P^{#1}_\grid}} 

\newcommand{\Q}{\ensuremath{\mathbb{Q}}}
\newcommand{\QG}[1][\grid]{\ensuremath{\mathbb{Q}(#1)}}
\newcommand{\QoG}[1][\grid]{\ensuremath{\mathbb{Q}_0(#1)}}

\newcommand{\REFINE}{\texttt{REFINE}}
\newcommand{\R}{\ensuremath{\mathbb{R}}}
\newcommand{\Rd}{\ensuremath{\mathbb{R}^d}}
\newcommand{\Rdd}{\ensuremath{\mathbb{R}^{d\times d}}}
\newcommand{\Rdds}{{\ensuremath{\Rdd_{\text{sym}}}}}
\newcommand{\Res}{\ensuremath{\mathcal{R}}}

\newcommand{\sides}{\mathcal{S}}

\DeclareMathOperator{\supp}{supp}
\newcommand{\SOLVE}[1]{\texttt{SOLVE}\ensuremath{{#1}}}

\newcommand{\tldt}{\ensuremath{\tilde t}}
\newcommand{\tr}{\ensuremath{\tilde r}}

\newcommand{\Tri}{\mathcal{B}}
\newcommand{\Trilin}[4][]{\ensuremath{\Tri_{#1}[#2,\,#3,\,#4]}}

\newcommand{\Un}[1][n]{\ensuremath{\vec{U}^{#1}}}
\newcommand{\UnG}[1][n]{\ensuremath{\vec{U}_{\grid}^{#1}}}
\newcommand{\Uk}[1][k]{\ensuremath{\vec{U}_{#1}}}

\renewcommand{\vec}[1]{\ensuremath{\boldsymbol{#1}}}
\newcommand{\vecu}{\ensuremath{\vec{u}}}

\newcommand{\vecE}{\ensuremath{\vec{E}}}
\newcommand{\vecv}{\vec{v}}
\newcommand{\vecV}{\vec{V}}

\newcommand{\vecf}{\vec{f}}
\newcommand{\vecw}{\vec{w}}

\newcommand{\V}{\ensuremath{\mathbb{V}}}

\newcommand{\VG}[1][\grid]{\V(#1)}
\newcommand{\VGk}[1][\grid_k]{\V(#1)}
\newcommand{\VoG}[1][\grid]{\V_0(#1)}

\providecommand{\vec}[1]{\boldsymbol{#1}}

\newcommand{\weak}{\ensuremath{\rightharpoonup}}

\algsetup{indent=2em}

\begin{document}


\title[AFEM for implicit power-law fluids]{Adaptive Finite element approximation
  of\\ steady flows of incompressible fluids with\\ implicit power-law-like rheology}


\author[C.~Kreuzer]{Christian Kreuzer}
\address{Christian Kreuzer,
 Fakult\"at f\"ur Mathematik,
 Ruhr-Universit\"at Bochum,
 Universit\"atsstrasse 150, D-44801 Bochum, Germany
 }%
\urladdr{http://www.ruhr-uni-bochum.de/ffm/Lehrstuehle/Kreuzer/index.html}
\email{christan.kreuzer@rub.de}

\author[E.~S\"uli]{Endre S\"uli}
\address{Endre S\"uli,
  Mathematical Institute, University of Oxford,
  24-29 St Giles', Oxford OX1 3LB, UK
}%
\urladdr{http://people.maths.ox.ac.uk/suli/index.html}
\email{Endre.Suli@maths.ox.ac.uk}



\begin{abstract}
  We develop the a posteriori error analysis of finite element
  approximations to implicit power-law-like models for viscous
  incompressible fluids in $d$ space dimensions, $d \in \{2,3\}$. The Cauchy stress and the symmetric part of
  the velocity gradient in the class of models under consideration are
  related by a, possibly multi-valued, maximal monotone $r$-graph,
  with $\frac{2d}{d+1}<r<\infty$. We establish upper and lower bounds on the finite
  element residual, as well as the local stability of the error
  bound. We then consider an adaptive finite element approximation of
  the problem, and, under suitable assumptions, we show the weak
  convergence of the adaptive algorithm to a weak solution of the
  boundary-value problem. The argument is based on a variety of weak
  compactness techniques, including Chacon's biting lemma and a finite
  element counterpart of the Acerbi--Fusco Lipschitz truncation of
  Sobolev functions, introduced by L. Diening, C. Kreuzer and
  E. S\"uli [Finite element approximation of steady flows of
  incompressible fluids with implicit power-law-like rheology. SIAM
  J. Numer. Anal., 51(2), 984--1015].
\end{abstract}

\maketitle

\section{Introduction\label{s:introduction}}

Typical physical models of fluid flow rely on the assumption that the Cauchy stress is an explicit function of
the symmetric part of the velocity gradient of the fluid. This constitutive hypothesis then leads to the
Navier--Stokes system and its nonlinear generalizations, such as
fluids with shear-rate-dependent viscosity including
  power-law fluids with constant or variable power-law index.
It is known however that the framework of classical continuum
mechanics, built upon the notions of current and reference
configuration and an
explicit constitutive equation for the Cauchy stress, is too narrow
for the accurate description of inelastic behavior of
solid-like materials or viscoelastic properties of materials. Our
starting point in this paper is therefore a
generalization of the classical framework of continuum mechanics,
referred to as implicit constitutive theory, which
was proposed recently in a series of papers by Rajagopal and
collaborators; see, for example,
\cite{Rajagopal:03,Rajagopal:06,RajagopalSrinivasa:08}.
The underlying principle
of implicit constitutive theory in the context of viscous flows is the following: instead of demanding that
the Cauchy stress is an explicit function of the symmetric part of the velocity gradient, one may allow an implicit
relationship between these quantities. This then leads to a general theory, which
admits fluid flow models with implicit and possibly discontinuous power-law-like rheology; see, \cite{Malek:07, Malek:08}.
Very recently a rigorous mathematical existence theory was developed for these models by Bul\'{\i}\v{c}ek, Gwiazda, M\'alek,
and \'Swierczewska-Gwiazda in \cite{BulGwiMalSwi:09}, for $r>\frac{2d}{d+2}$;
for the range $\frac{2d}{d+2}<r\le\frac{3d}{d+2}$ the Acerbi--Fusco Lipschitz truncation \cite{AcerbiFusco:88} was used in
order to prove the existence of a weak solution.
In \cite{DieningKreuzerSueli:2013}, using a variety of weak compactness
techniques, we showed that a subsequence of the sequence of
finite element solutions converges weakly to a weak solution of the problem as the finite element discretization parameter
$h$ tends to $0$. A key new technical tool in the analysis presented in \cite{DieningKreuzerSueli:2013} was
a finite element counterpart of the Acerbi--Fusco Lipschitz truncation
of Sobolev functions. However, in the case of velocity approximations that are not exactly
divergence-free the convergence theory developed there was restricted to the range $\frac{2d}{d+1}<r<\infty$.

The focus of the present paper is on the adaptive finite
  element approximation of implicitly constituted
power-law-like models for viscous incompressible fluids.
As in \cite{DieningKreuzerSueli:2013}, the implicit constitutive
relation between the stress and the symmetric part of the velocity
gradient is approximated by an explicit (smooth) constitutive law. The
resulting steady non-Newtonian flow problem is then discretized by a
mixed finite element method. Guided by an a posteriori error analysis, we propose
a numerical method with competing adaptive strategies for the mesh refinement
and the approximation of the implicit constitutive law, and we present a
rigorous convergence proof generalizing the ideas in \cite{MoSiVe:08}
and \cite{Siebert:11}. More precisely, we show that a subsequence
of the adaptively generated sequence of discrete approximations
converges, in the weak topology of
the ambient function space, to a weak solution of the model when
  $\frac{2d}{d+1}<r<\infty$. In contrast with \cite{DieningKreuzerSueli:2013}, stimulated by ideas
  from \cite{BulGwiMalSwi:12} we shall be able to avoid resorting to the theory of Young measures.
We
emphasize that even in the case when the weak solution is unique
we have only weak convergence of a subsequence; in this case, however,
such a subsequence can be identified with the aid of the a
posteriori bounds derived herein; cf. Remark~\ref{R:filter-conv}.

The paper is structured as follows. In
  Section~\ref{sec:preliminaries} we shall formulate the
  problem under consideration and will introduce some known mathematical results.
  In Section~\ref{s:fem_approx} we define the finite element
  approximation of the problem and present related technical
  properties and tools, such as the discrete Lipschitz truncation from
  \cite{DieningKreuzerSueli:2013}.  Section~\ref{sec:aposteriori} is
  concerned with the a posteriori error analysis for both the
   error in the approximation of the graph and the finite element approximation. The adaptive
   algorithm together with our main result are stated in
   Section~\ref{sec:afem+conv}; for the sake of clarity of the presentation
   certain technical parts of the proof are deferred to Section~\ref{sec:aux}.
   We conclude the paper by discussing concrete
   graph approximations for certain problems of practical relevance.
While the emphasis in this paper is on the mathematical analysis of adaptive finite element algorithms for implicitly constituted fluid flow
models, the ideas developed herein may be of more general
interest in the convergence analysis of adaptive finite element
methods for other nonlinear problems in continuum mechanics with
possibly nonunique weak solutions.

\section{Implicitly constituted power-law-like fluids}
\label{sec:preliminaries}
In this section we introduce the variational model of steady flow,
in a bounded open Lipschitz domain $\Omega \subset \mathbb{R}^d$,
$d \in \{2,3\}$, with polyhedral boundary $\partial \Omega$,
of an incompressible fluid with an implicit constitutive law given
by a maximal monotone $x$-dependent $r$-graph.
We then recall the existence result from \cite{BulGwiMalSwi:09}
together with some known results and mathematical
tools from the literature.

\subsection{The variational formulation\label{s:problem}}
Before stating the weak formulation of the problem we
need to introduce basic notations and recall some well-known
properties of Lebesgue and Sobolev function spaces.

For a measurable subset $\omega\subset\R^d$, we denote the classical
spaces of  Lebesgue and vector-valued Sobolev functions
by $(L^s(\omega) :=L^s(\omega;\R),\norm[s;\omega]{\cdot})$ and
$(W^{1,s}(\omega)^d :=W^{1,s}(\omega;\Rd),\norm[1,s;\omega]{\cdot})$, $s\in[1,\infty]$, respectively.
Henceforth $\omega$ will be assumed to have Lipschitz continuous boundary.
We denote the space of functions in $W^{1,s}(\omega)^d$ with zero
trace by  $W^{1,s}_0(\omega)^d$  and let
$W^{1,s}_{0,\divo}(\omega)^d\definedas \{\vecv\in
W^{1,s}_0(\omega)^d\colon \divo\vecv = 0\}$. Moreover, we denote
the space of functions in $L^s(\omega)$ with zero
integral mean by  $L^s_0(\omega)$. For $s,s'\in(1,\infty)$ with  $\frac1s+\frac1{s'}=1$
we have that $L^{s'}(\Omega)$ and $L^{s'}_0(\Omega)$ are the dual
spaces of $L^{s}(\Omega)$ and $L^{s}_0(\Omega)$, respectively. We
  have, for such $s$ and $s'$,
that $W^{1,s}_0(\Omega)^d$ is the closure of $\mathcal{D}(\Omega)^d\definedas C^\infty_0(\omega)^d$
and its dual  is denoted by $W^{-1,s'}(\Omega)^d$.
For $\omega=\Omega$ we omit the domain in our notation for norms; \eg, we write
$\norm[s]{\cdot}$ instead of $\norm[s,\Omega]{\cdot}$.

For $r\in(1,\infty)$, we define $r'\in(1,\infty)$ by $\tfrac1r+\tfrac1{r'}=1$, and set
%
\begin{align}\label{eq:tr}
  \tr \definedas
  \begin{cases}
    \tfrac{1}{2} \tfrac{dr}{d-r}\quad&\text{if}~r\leq \frac{3d}{d+2}
    \\
    r'\quad &\text{otherwise}.
  \end{cases}
\end{align}

With such $r$, $r'$ and $\tilde{r}$, we shall consider the following boundary-value problem.

\smallskip

 \paragraph{Problem}
 For $\vecf\in L^{r'}(\Omega)^d$ find
  $(\vec{u},p,\mat{S})\in W_0^{1,r}(\Omega) ^d\times L_0^{
    \tr}(\Omega)\times L^{r'}(\Omega;\Rdds)$ such that
  \begin{align}\label{eq:implicit}
    \begin{alignedat}{2}
      \divo(\vec{u}\otimes\vec{u} +  p\bbI -\vec{S})&=\vecf &\quad
      &\text{in}~ \mcD'(\Omega)^d,
      \\
      \divo\vec{u} &=0& \quad&
      \text{in}~ \mcD'(\Omega),
     \\
       (\nablas \vec{u}(x), \mat{S}(x)) &\in\Ax(x)&\quad&\text{for almost
         every}~x\in\Omega.
    \end{alignedat}
  \end{align}
  Here, $\nablas \vec{u}\definedas \frac12(\nabla \vec{u}+(\nabla
  \vec{u})^{\rm T})\in \Rdds:=\{\vec{\delta}\in\R^{d\times d}:
  \vec{\delta}=\vec{\delta}^{\rm T}\}$
  signifies the symmetric part of the gradient of $\vec{u}$.
  As is indicated in our
  choice of the solution space for the velocity $\vec{u}$ in the statement of the above
  boundary-value problem, we shall suppose a homogenous Dirichlet boundary condition for $\vecu$.
  The integrability of the pressure $p$ is inherited from the
    convective term and therefore the definition~\eqref{eq:tr} of $\tr$ is a
    consequence of the embedding $W^{1,r}_0(\Omega)\hookrightarrow
    L^{2\tr}(\Omega)$.
  The implicit constitutive law, which relates the shear rate to the shear stress,
  is given by an inhomogeneous maximal monotone $r$-graph
  $\Ax:x\mapsto\Ax(x) \subset\Rdds\times\Rdds$. In particular, we assume that for almost
  every $x\in\Omega$ the following properties hold:
  \begin{enumerate}[leftmargin=1cm,itemsep=1ex,label=(A\arabic{*})]
\item\label{A1} $(\vec{0},\vec{0})\in\Ax(x)$;
\item\label{A2} For all $(\vec{\delta}_1,\vec{\sigma}_1),(\vec{\delta}_2,\vec{\sigma}_2)\in\Ax(x)$,
  \begin{align*}
    (\vec{\sigma}_1-\vec{\sigma}_2):(\vec{\delta}_1-\vec{\delta}_2)\geq 0 \qquad\text{($\Ax(x)$ is a
      monotone graph)};
  \end{align*}
\item\label{A3}
  If $(\vec{\delta},\vec{\sigma})\in\Rdd_{\text{sym}}\times\Rdd_{\text{sym}}$ and
  \begin{align*}
    (\bar{\vec{\sigma}}-\vec{\sigma}):(\bar{\vec{\delta}}-
    \vec{\delta})\geq 0\quad\text{for all}~
      (\bar{\vec{\delta}},\bar{\vec{\sigma}})\in \Ax(x),
  \end{align*}
  then $(\vec{\delta},\vec{\sigma})\in\Ax(x)$ (i.e., $\Ax(x)$ is a maximal monotone graph);
  \item\label{A4}There exists a nonnegative function $m\in L^1(\Omega)$ and a constant
    $c>0$, such
    that for all $(\vec{\delta},\vec{\sigma})\in\Ax(x)$ we have
    \begin{align*}
      \vec{\sigma}:\vec{\delta}\ge - m(x)+c(|\vec{\delta}|^r+|\vec{\sigma}|^{r'})\qquad\text{(i.e., $\Ax(x)$
        is an $r$-graph);}
    \end{align*}
    \item\label{A5}
      The set-valued mapping $\Ax:\Omega\to\Rdds\times\Rdds$ is
      measurable, i.e.,  for any closed sets
      $\mathcal{C}_1,\mathcal{C}_2\subset\Rdds$, we have that
      \begin{align*}
        \big\{x\in\Omega:
        \Ax(x)\cap(\mathcal{C}_1\times\mathcal{C}_2)\neq\emptyset\big\}
      \end{align*}
      is a Lebesgue measurable subset of $\Omega$.
\end{enumerate}

The following existence result was originally
proved by Bul\'{\i}\v{c}ek, Gwiazda, M\'alek,
and \'Swierczewska-Gwiazda in
\cite{BulGwiMalSwi:09} assuming additionally that if
$\vec{\delta}_1\neq\vec{\delta}_2$ and $\vec{\sigma}_1\neq
\vec{\sigma}_2$, then the inequality
in~\ref{A2} is strict. In fact, based on a generalization of the fundamental
theorem on Young measures, this condition was required in order to prove
that the implicit constitutive law is satisfied. For the unsteady case,
they presented a new technique in
\cite{BulGwiMalSwi:12} avoiding the additional condition. This
technique can also be applied to steady problems~\eqref{eq:implicit};
compare with~\cite{BulGwiMalRajSwi:12}.

\begin{prpstn}
  For $r>\frac{2d}{d+2}$ there exists a (not necessarily unique)
  weak solution to problem~\eqref{eq:implicit}.
\end{prpstn}

 \begin{rmrk}\label{rem:Aprop}
  Two remarks concerning the definition of $x$-dependent maximal monotone graph
  $\Ax$ are now in order:
  \begin{itemize}[leftmargin=0.5cm]
  \item Let $\mathcal{D}\subset\Rdds\times\Rdds$ be a closed set; then,
    for a.e. $x\in\Omega$ we have that the set $\Ax(x)\cap\mathcal{D}$
    is closed. To see this, we assume w.l.o.g. that
    $\Ao(x)\cap\mathcal{D}\neq\emptyset$ and let
    $\{(\vec{\delta}_k,\vec{\sigma}_k)\}_{k\in\N}\subset
    \Ao(x)\cap\mathcal{D}$, such that $\vec{\delta}_k\to
    \vec{\delta}\in\Rdds$ and $\vec{\sigma}_k\to\vec{\sigma}\in\Rdds$
    as $n\to\infty$.  Let
    $(\bar{\vec{\delta}},\bar{\vec{\sigma}})\in\Ao(x)$ be arbitrary. We then
    have that
    \begin{align*}
      0\le(\bar{\vec{\sigma}}-\vec{\sigma}_k):(\bar{\vec{\delta}}-
      \vec{\delta}_k)\to(\bar{\vec{\sigma}}-\vec{\sigma}):(\bar{\vec{\delta}}-
      \vec{\delta})
    \end{align*}
    as $k\to\infty$. This proves that
    $(\vec{\delta},\vec{\sigma})\in\Ao(x)\cap\mathcal{D}$ thanks to condition~\ref{A3} and the
    closedness of $\mathcal{D}$.
    Taking $\mathcal{D}:=\{\vec{\sigma}\}\times\Rdds$, we then deduce that the set
    \begin{align*}
      \{\vec{\sigma}\in\Rdd_{\text{sym}}:(\vec{\delta},\vec{\sigma})\in\Ax(x)\}
    \end{align*}
    is closed. This is condition $(A5)(i)$ of
    \cite{BulGwiMalSwi:09}.
    \item According to \cite[Theorem 8.1.4]{AubinFrankowska:2009} property \ref{A5}
      is equivalent to the fact that the graph of the set-valued map
      $\Ax(x)$ belongs to the product $\sigma$-algebra
      $\mathfrak{L}(\Omega)\otimes\mathfrak{B}(\Rdd_{\text{sym}})
      \otimes\mathfrak{B}(\Rdd_{\text{sym}})$. Here
      $\mathfrak{L}(\Omega)$ denotes the Lebesgue measurable subsets of
      $\Omega$ and
      $\mathfrak{B}(\Rdd_{\text{sym}})$  the Borel subsets of $\Rdd_{\text{sym}}$.
      With the same argument it follows that~\ref{A5} is equivalent
      to the fact that, for any closed
      $\mathcal{C}\subset\Rdd_{\text{sym}}$, the sets
      \begin{align*}
        \hspace{2cm}&\big\{(x,\vec{\sigma})\in\Omega\times\Rdd_{\text{sym}}: \text{there exists }
        \vec{\delta}\in\mathcal{C},\text{ such that }
        (\vec{\delta},\vec{\sigma})\in\Ax(x)\big\},
        \\
        \hspace{2cm}&\big\{(x,\vec{\delta})\in\Omega\times\Rdd_{\text{sym}}: \text{there exists }
        \vec{\sigma}\in\mathcal{C},\text{ such that }
        (\vec{\delta},\vec{\sigma})\in\Ax(x)\big\}
      \end{align*}
      are measurable relative to
      $\mathfrak{L}(\Omega)\otimes\mathfrak{B}(\Rdd_{\text{sym}})$. These
      equivalent  conditions imply that there exist measurable
      functions (so-called selections)
      $\bsS^\star,\bsD^\star:\Omega\times\Rdds\to \Rdds$ such that
      $\big(\vec{\delta},\bsS^\star(x,\vec{\delta})\big),\big(\bsD^\star(x,\vec{\sigma}),\vec{\sigma}\big)\in\Ax(x)$
      for a.e. $x\in\Omega$ and all
      $\vec{\delta},\vec{\sigma}\in\Rdds$; compare  also with \cite[Remark 1.1]{BulGwiMalSwi:12}.
  \end{itemize}

\end{rmrk}

\subsection{Analytical framework\label{s:inf-sup}}
We shall briefly recall some results that are crucial for the existence
theory for problem~\eqref{eq:implicit}.

  \paragraph{Inf-sup condition} The inf-sup condition has a central role in the analysis of the
  Stokes problem. It states that, for
  $s,s'\in(1,\infty)$ with  $\frac1s+\frac1{s'}=1$,
  there exists an $\alpha_s>0$ such that
  \begin{align}\label{eq:inf-sup}
    \sup_{0\neq \vecv\in W^{1,s}_0(\Omega)^d}\frac{\int_\Omega q\divo \vecv\dx}{\norm[1,s]{\vecv}}\geq
    \alpha_s \norm[s']{q}\qquad\text{for all}~q\in L^{s'}_0(\Omega).
  \end{align}
   This is the consequence of the existence of the {\em
        Bogovski\u\i} operator $\mathfrak{B}:L^s_0(\Omega)\to
      W^{1,s}_0(\Omega)^d$, with
    \begin{align*}
      \divo \mathfrak{B}h=h
        \qquad\text{and}\qquad\alpha_s\norm[{1,s}]{\mathfrak{B}h}\le
        \norm[s]{h}
    \end{align*}
    for all $s\in(1,\infty)$; compare e.g. with \cite{DieSchuRu:10,Bogovskii:79}.
    It follows from
    \cite[\S II, Proposition 1.2]{BrezziFortin:91} that condition~\eqref{eq:inf-sup}
    is equivalent to the isomorphism
    \begin{align}\label{eq:isom}
      L_0^{s'}(\Omega)\cong \left\{\vecv'\in W^{-1,s'}(\Omega)^d\colon
      \langle \vecv',\,\vecw\rangle=0
      ~\text{for all $\vecw\in W^{1,s}_{0,\divo}(\Omega)^d$}\right\}.
  \end{align}
  \paragraph{Korn's inequality} According to \eqref{eq:implicit} the
  maximal monotone graph defined in \ref{A1}--\ref{A5} provides
  control only of the symmetric part of the velocity gradient. Korn's inequality states
  that this already suffices to control the norm of a Sobolev
  function; \hbox{i.\,e.},\xspace for $s\in(1,\infty)$, there exists a $\gamma_s>0$ such that
  \begin{align}
    \label{eq:korn}
    \gamma_s \norm[1,s]{\vecv}\le\norm[s]{\nablas \vecv}\qquad
    \text{for all}~\vecv\in W^{1,s}_0(\Omega)^d;
  \end{align}
  compare e.\,g. with \cite{DieSchuRu:10}.

  We conclude this subsection with Chacon's biting lemma and a
  corollary of it that is relevant for our purposes; compare e.g. with
  \cite{BrooksChacon:80} and \cite[Lemma 7.3]{GwiaMalSwier:07}.
   \begin{lmm}[Chacon's biting lemma]\label{l:biting}
   Let $\Omega$ be a bounded domain in $\R^d$ and let
   $\{v_n\}_{n\in\N}$ be a bounded sequence in $L^1(\Omega)$. Then,
   there exists a nonincreasing sequence of measurable subsets
   $E_j\subset\Omega$ with $\abs{E_j}\to0$ as $j\to\infty$, such that $\{v_n\}_{n\in\N}$ is
   precompact in the weak topology of $L^1(\Omega\setminus E_j)$, for
   each $j\in\N$.

   In other words, there exists a $v\in L^1(\Omega)$,
   such that for a subsequence (not relabelled) of $\{v_n\}_{n \in \N}$,
   $v_n\weak v$ weakly in $L^1(\Omega\setminus E_j)$ as $n\to\infty$ for all
   $j\in\N$. We denote this by writing
   \begin{align*}
     v_n\bweak v\quad\text{in}~L^1(\Omega)
   \end{align*}
   and call $v$ the biting limit of $\{v_n\}_{n\in\N}$.
 \end{lmm}

 \begin{lmm}\label{l:bit=>L1}
   Let $\{v_n\}_{n\in\N}\subset L^1(\Omega)$ be a sequence of nonnegative
   functions such that $ v_n\bweak v$ for some $v\in L^1(\Omega)$.
   Then,
   \begin{align*}
   \lim_{n\to\infty}\int_\Omega (v_n-v)\dx =0\quad\text{implies that}\quad
   v_n\weak v \quad\text{weakly in}~L^1(\Omega)~\text{as}~n\to\infty.
 \end{align*}

 \end{lmm}

\section{Finite Element Approximation\label{s:fem_approx}}
This section is concerned with approximating problem~\eqref{eq:implicit}
by the finite element method. To this end, we first approximate~\eqref{eq:implicit}
 by an explicitly constituted problem.
We then introduce a general finite element
framework for inf-sup stable Stokes elements.
This, together with some representative examples of velocity-pressure pairs of
finite element spaces, is the subject of Section
\ref{ss:fem_spaces}. 
The finite element approximation of~\eqref{eq:implicit} is stated
in Section \ref{ss:Galerkin}.

\subsection{Approximation of maximal monotone
  $r$-graphs\label{ss:approxAx}}
In general an $x$-depen\-dent maximal monotone
$r$-graph $\Ax$ satisfying~\ref{A1}--\ref{A5} cannot be represented
in an explicit fashion. However, it can be approximated by a regular single-valued
monotone tensor field based on a regularized {\em
  measurable selection} $\bsS^\star$ with the following properties; compare with
\cite{BulGwiMalSwi:09,BulGwiMalSwi:12} and Remark~\ref{rem:Aprop}.
\begin{lmm}[\!\!{\cite[Lemma 2.2]{BulGwiMalSwi:12}}]\label{lem:Sstar}
  Let $\bsS^\star:\Omega\times\Rdds\to\Rdds$ be a measurable
  selection of the $x$-dependent maximal monotone $r$-graph $\Ax$ with
  the properties \ref{A1}--\ref{A5}. Then, for
  $\vec{\delta},\vec{\sigma}\in\Rdds$, the following two statements
  are equivalent for almost all $x\in\Omega$:
  \begin{itemize}
    \item
    $
    \big(\vec{\sigma}-\bsS^\star(x,\bsD)\big):(\vec{\delta}-\bsD)\ge
    0\qquad \text{for all}~\bsD\in\Rdds$;
    \item
      $(\vec{\delta},\vec{\sigma})\in\Ax(x).$
  \end{itemize}

\end{lmm}

In \cite{GwiaMalSwier:07,GwiazdaZatorska:07,BulGwiMalSwi:09} the
selection $\bsS^\star$ is used to approximate the
maximal monotone graph $\Ax$ by a single-valued monotone mapping
$\bSn:\Omega\times\Rdds\to\Rdds$ based
on a mollification technique. In order to allow for different
practical implementations of such an approximation, we shall formulate
its required properties and demonstrate in Section~\ref{sec:graph-approximation}
how such graph-approximations can
be constructed
for some typical problems of practical
interest within the class of problems under consideration.

\begin{ass}\label{ass:Sn}
  For $n\in\N$, there exists a  mapping
  $\bSn:\Omega\times\Rdds\to\Rdds$, such that
  \begin{itemize}
  \item $\bSn(\cdot,\vec{\delta}):\Omega\to\Rdds$ is measurable for
    all $\vec{\delta}\in\Rdds$;
  \item $\bSn(x,\cdot):\Rdds\to\Rdds$ is continuous for almost every $x\in\Omega$;
  \item $\bSn$ is strictly monotone; i.e., for almost every $x\in\Omega$
    we have
    \begin{align*}
      \left(\bSn(x,\vec{\delta}_1)-\bSn(x,\vec{\delta}_2)\right):\left(\vec{\delta}_1-\vec{\delta}_2\right)
      >0 \quad\text{for
        all}\quad\vec{\delta}_1\neq\vec{\delta}_2\in\Rdds;
    \end{align*}
  \item
    There exist constants $\tilde c_1,\tilde c_2>0$ and nonnegative functions
  $\tilde m\in L^{1}(\Omega)$, $\tilde k\in L^{r'}(\Omega)$ such that,
  uniformly in $n\in\N$, we have
  \begin{alignat*}{2}
    \abs{\bSn(x,\vec{\delta})}&\le \tilde c_1 \abs{\vec{\delta}}^{r-1}+\tilde
    k(x)
    \quad\text{and}\quad
        \bSn(x,\vec{\delta}):\vec{\delta}\ge \tilde c_2 \abs{\vec{\delta}}^{r}-\tilde
    m(x)
  \end{alignat*}
  for all $ \vec{\delta}\in\Rdds$ and almost every $x\in\Omega$.
  \end{itemize}
\end{ass}
We emphasize that in contrast with
  \cite{BulGwiMalSwi:09,BulGwiMalSwi:12} we assume that $\bSn$ is
  strictly monotone. Of course we have to assume additionally that the
  graph of $\bSn$ converges to $\Ax$ is some sense. This will
  be specified in Assumption~\ref{ass:estAconv} with the aid of the a
  posteriori
  graph approximation indicator formulated in Section~\ref{sec:graph-error}.

Having at hand a mapping $\bSn:\Omega\times\Rdds\to\Rdds$ as in
  Assumption~\ref{ass:Sn}, we
aim to approximate the solution of \eqref{eq:implicit} by solving the
following explicitly constituted nonlinear boundary-value problem:
For $\vecf\in L^{r'}(\Omega)^d$ find
$(\vec{u},p,\mat{S})\in W_0^{1,r}(\Omega) ^d\times L_0^{
  \tr}(\Omega)\times L^{r'}(\Omega,\Rdds)
$ such that
\begin{align}\label{eq:explicit}
  \begin{alignedat}{2}
      \divo(\vec{u}\otimes\vec{u} +  p\bbI -\bsS)&=\vecf &\quad
      &\text{in}~ \mcD'(\Omega)^d,
      \\
      \divo\vec{u} &=0& \quad&
      \text{in}~ \mcD'(\Omega),
      \\
      \bsS(x)&=\bSn(x,\nablas\vecu(x))&\quad&\text{for almost every}~x\in\Omega.
    \end{alignedat}
  \end{align}

\subsection{Domain partition and refinement framework\label{sec:grid}}
In this section we provide the framework for adaptive grid
refinement. For the sake of simplicity, we restrict our presentation
to conforming simplicial meshes and refinement by bisection. To be more
precise, let $\grid_0$
be a regular conforming partition of $\Omega$ into closed simplexes, the
so-called macro mesh. Each simplex in the partition is referred to as
an element.
We assume that there exists a refinement routine
$\REFINE$ with the following properties.
\begin{itemize}[leftmargin=0.5cm]
\item The refinement routine has two input arguments: a regular
  conforming partition $\grid$
  and a subset $\marked\subset\grid$ of marked elements. The output is
  a refined regular conforming triangulation of $\Omega$, where all
  elements in $\marked$ have been bisected at least once.
  The input
  grid can be $\grid_0$ or the output of a previous application of
  \REFINE.
\item {\em Shape-regularity:} We call $\grid'$ a refinement of $\grid$
  (briefly $\grid'\ge\grid$), when it can be produced
  from $\grid$ by a finite number of applications of \REFINE. The set
  \begin{align*}
    \grids\definedas \{\grid\colon \grid~\text{is a refinement of}~\grid_0\}
  \end{align*}
  is shape-regular, \hbox{i.\,e.},\xspace for any element $\elm\in\grid$ with $\grid\in\grids$, the
  ratio of its diameter to the diameter of the largest inscribed ball
  is bounded uniformly with respect to all partitions $\grid$ with $\grid\in\grids$.
\end{itemize}

For the proof of existence of such a procedure, we refer to
\cite{Baensch:91}, \cite{Kossaczky:94}, \cite{Stevenson:08} or the monograph
\cite{SchmidtSiebert:05} and the references therein.

  For every element $\elm\in\grid$, $\grid\in\grids$,
  there exists an invertible affine mapping
  \begin{align*}
    \boldsymbol{F}_\elm:\elm\to \hat \elm,
  \end{align*}
  where $\hat\elm$ is the standard reference $d$-simplex.
The neighbourhood of an element $\elm\in\grid$, with $\grid \in\grids$, is denoted by
\begin{align*}
  \mathcal{N}^\grid(\elm):=\{\elm'\in\grid\colon \elm'\cap\elm\neq\emptyset\}.
\end{align*}
Let $\omega\subset\Omega$ and define
$\mathcal{U}^\grid(\omega)\definedas \bigcup\left\{\elm\in\grid\mid \elm\cap
  \omega\neq\emptyset\right\}$.
For subsets $\mathcal{M}\subset\grid$, let
\begin{align*}
  \Omega({\mathcal{M}})\definedas \bigcup\{\elm\mid
  \elm\in\mathcal{M}\} \subset\Omega\quad\text{and}\quad
  \mathcal{U}^\grid(\mathcal{M})\definedas \mathcal{U}^\grid(\Omega(\mathcal{M})) \subset\Omega,
\end{align*}
i.e., we have $\Omega=\Omega(\grid)$. Thanks to the shape-regularity of
$\grids$, we have that $\# \mathcal{N}^\grid(\elm)\le C$ and
$\abs{\mathcal{U}^\grid(\elm)}=\abs{\Omega(\mathcal{N}^\grid(\elm))}\le
C\abs{\elm}$ with a constant $C>0$ independent of $\grid\in\grids$.
For
$\grid\in\grids$, we define the mesh-size function
\begin{align*}
  \Omega\ni x\mapsto h_\grid(x)\definedas \abs{\mathcal{U}^\grid({\{x\})}}^{1/d}.
\end{align*}
For $x\in\operatorname{interior}(\elm)$, this coincides with the usual
definition $h_\grid(x)=\abs{\elm}^{1/d}=:\!h_\elm$. The mesh-size
function is monotonically decreasing under refinement.

We call the $(d-1)$-dimensional sub-simplexes of any simplex $\elm
\in\grid$, whose interiors
lie inside $\Omega$, the sides of
$\grid$ and denote the set of all of them by $\sides(\grid)$. For $S\in\sides(\grid)$,
we define $h_S\definedas
|S|^{1/(d-1)}$ and observe for $x\in S$ that $c
h_S\le\hG(x)\le Ch_S$, with constants $C,c>0$
depending solely on the shape-regularity of~$\grid$.

\subsection{Finite element spaces\label{ss:fem_spaces}}
Denote by $\P_m$ the space of polynomials of degree at most
 $m\in\N$. For a given grid $\grid\in\grids$ and certain subspaces
 $\Q\subseteq L^\infty(\Omega)$ and $\V\subseteq W^{1,\infty}_0(\Omega)^d$ the
 finite element spaces are given by
\begin{subequations}\label{df:VQ}
  \begin{align}
    \VG&\definedas \left\{\vecV\in\V~\colon~
      \vecV|_{\elm}\circ\boldsymbol{F}_\elm^{-1}\in\hat{\mathbb{P}}_\V,
      ~\elm \in\grid~\text{and}~
      \vecV|_{\partial\Omega}=0\right\},\label{df:Vn}
    \\
    \QG&\definedas \left\{Q\in \Q~\colon
      ~Q|_{\elm}\circ\boldsymbol{F}_\elm^{-1}\in\hat{\mathbb{P}}_\Q,~\elm\in\grid
    \right\},
    \label{df:Qn}
  \end{align}
\end{subequations}
where $\hat{\mathbb{P}}_\V\subset W^{1,\infty}(\hat\elm)^d$ and
$\hat{\mathbb{P}}_\Q\subset L^\infty(\hat\elm)$
are finite-dimensional subspaces such that
$\mathbb{P}_1^d\subseteq\hat{\mathbb{P}}_\V\subseteq\mathbb{P}_\ell^d$
and
$\mathbb{P}_{0}\subseteq\hat{\mathbb{P}}_\Q\subseteq\mathbb{P}_\jmath$ for some
$\ell\ge \jmath\in\N$.
For convenience, we introduce the space of piecewise polynomials of
degree at most $m\in\N$ over $\grid$ by
\begin{align*}
  \P_m(\grid):=\{R:\bar\Omega\to\R~\colon R|_{\elm}\in\P_m,~\elm\in\grid\}.
\end{align*}
Note
that $\QG\subset L^\infty(\Omega)\cap\P_{\jmath}(\grid)$ and
since $\VG\subset C_0(\bar\Omega)^d\cap\P_\ell(\grid)^d$ it follows that $\VG\subset
W^{1,\infty}_0(\Omega)^d$. Additionally, we assume 
that the finite element
spaces are nested, i.e.,  
if $\grid_\star$ is a refinement of $\grid$, then
\begin{align}\label{eq:nested}
  \VG\subset\VG[\grid_\star]\qquad\text{and}\qquad\QG\subset\QG[\grid_\star].
\end{align}
Each of the above spaces is supposed to have a finite and locally
supported basis; e.\,g. for the discrete velocity space this means that
for $\grid\in\grids$
there exists an $N_\grid\in\N$ such that
\begin{align*}
  \VG=\operatorname{span}\{\vecV_1^\grid,\ldots,\vecV^\grid_{N_\grid}\}
\end{align*}
and for each basis function $\vecV_i^\grid$, $i=1,\ldots,N_\grid$, we have that if
there exists an $\elm\in\grid$ with $\vecV_i^\grid\not\equiv 0$
 on $\elm$, then
$
  \supp{\vecV_j^\grid}\subset
  \mathcal{U}^\grid(\elm).
$
We introduce the subspace 
$\VoG$ of  discretely divergence-free functions by
\begin{align*}
  \VoG\definedas \Big\{\vecV\in\VG\colon
  \int_\Omega Q\divo\vecV\dx=0~\text{for all}~Q\in\QG\Big\}
\end{align*}
and we define
\begin{align*}
  \QoG\definedas \Big\{Q\in\QG:\int_\Omega Q\dx=0\Big\}.
\end{align*}

It will be assumed throughout the paper that all pairs of velocity-pressure
finite element spaces considered possess the following properties.

\begin{ass}[Projector $\PGdiv$]\label{ass:Pndiv}
  We assume that for each $\grid\in\grids$
  there exists a linear projection operator
  $\PGdiv\colon W^{1,1}_0(\Omega)^{d}\to\VG$  such that, for all $s\in(1,\infty)$,
  \begin{itemize}[leftmargin=0.5cm]
  \item  $\PGdiv$ preserves divergence in ${\QG}^*$; \hbox{i.\,e.},\xspace for
    $\vecv\in W^{1,s}_0(\Omega)^{d}$ we have
    \begin{align*}
      \int_\Omega Q\divo \vecv\dx&=\int_\Omega Q\divo\PGdiv \vecv\dx\qquad\text{for
        all}~Q\in\QG.
    \end{align*}
  \item $\PGdiv$ is locally defined; i.e.,  for any other partition
    $\grid_\star\in\grids$ we have
    \begin{align}\label{eq:locality}
      \PGdiv[\grid_\star]\vecv|_{\mathcal{U}^{\grid_\star}(\elm)}=\PGdiv\vecv|_{\mathcal{U}^\grid(\elm)}
    \end{align}
     for all $~\vecv\in W^{1,s}_0(\Omega)^{d}$ and all $\elm\in\grid$
     with $\mathcal{N}^\grid(\elm)\subset \grid_\star$.
  \item  $\PGdiv$ is locally $W^{1,1}$-stable; \hbox{i.\,e.},\xspace there exists a $c_1>0$,
    independent of $\grid$, such that
    \begin{align}\label{eq:stability}
    \int_\elm \abs{\PGdiv \vecv}+\hG\abs{\nabla\PGdiv \vecv}\dx\le c_1
    \int_{\mathcal{U}^\grid(\elm)}\abs{ \vecv}+\hG\abs{\nabla
      \vecv}\dx
  \end{align}
  for all $~\vecv\in W^{1,s}_0(\Omega)^{d}$ and all $\elm\in\grid$.
 \end{itemize}
\end{ass}
As in \cite{DieningKreuzerSueli:2013,BelBerDieRu:12,DieningRuzicka:07b}, the local
$W^{1,1}$-stability property \eqref{eq:stability}  implies global
$W^{1,s}$-stability, i.e., for each $s\in[1,\infty]$, there exists
a $c_s>0$, such that
\begin{align}\label{eq:W1s-stab}
  \norm[1,s]{\PGdiv \vecv}\le c_s \norm[1,s]{ \vecv}\qquad \text{for
    all}~\vecv\in W^{1,s}_0(\Omega)^{d}.
\end{align}
Moreover, since $\VG$ contains piecewise affine functions, we have the
following interpolation error bound.
  For each $s\in[1,\infty]$ there exists a $c_s>0$ such that
  \begin{align}\label{eq:interpolation}
  \int_\elm \abs{\vecv-\PGdiv \vecv}^s+\hG^s\abs{\nabla\vecv-\nabla\PGdiv \vecv}^s\dx\le c_s
  h_\elm^{s(1+\delta)}\abs{
      \vecv}_{W^{1+\delta,s}(\mathcal{U}^\grid(\elm))}^s
  \end{align}
  for all $\elm\in \grid$ and $\vecv\in
  W^{1+\delta,s}(\Omega)^{d}\cap W^{1,s}_0(\Omega)^d$, $\delta\in\{0,1\}$.

 As a consequence, we deduce the following result for weak limits in
 nested spaces. Before stating the result, we adopt the following notational convention:
 we shall write $A \Cleq B$ to denote $A \leq C\cdot B$ with a constant $C>0$ that is independent
 of the discretization parameter $h$.

\begin{prpstn}\label{prop:Pnweak}
  Let $\{\vecv_k\}_{k\in\N}\subset W^{1,s}_0(\Omega)^{d}$,
  $s\in(1,\infty)$, be such that $\vecv_k\weak \vec{0}$ weakly in $W^{1,s}_0(\Omega)^{d}$ as
  $k\to\infty$ and let $\{\grid_k\}_{k\in\N}\subset \grids$ be a sequence of
  nested partitions of $\Omega$, i.e., $\gridk\le\grid_{k+1}$ for all $k\in\N$. Then,
\begin{align*}
  {\PGdivk\vecv_k\weak \vec{0}\quad\text{weakly in}~W^{1,s}_0(\Omega)^{d}~\text{as
  $k\to\infty$}.}
\end{align*}
\end{prpstn}

\begin{proof}
  Thanks to the uniform boundedness~\eqref{eq:W1s-stab} of the sequence of linear
    operators $\{\PGdivk\,:\, W^{1,s}_0(\Omega)^{d} \rightarrow
    \VG[\gridk] \subset W^{1,s}_0(\Omega)^{d}\}_{k\in\N}$, we have
  that there exists a not relabelled weakly converging subsequence of
  $\{\PGdivk\vecv_k\}_{k\in\N}$ in
  $W^{1,s}_0(\Omega)^{d}$.
  By the compact embedding $W^{1,s}_0(\Omega)^{d}\hookrightarrow
  \hookrightarrow L^s(\Omega)^{d}$ the sequence $\{\PGdivk\vecv_k\}_{k\in\N}$ converges
  strongly in $L^s(\Omega)^d$. Thanks to the uniqueness of the strong limit, it
  suffices to identify the limit of $\{\PGdivk\vecv_k\}_{k\in\N}$  in $L^s(\Omega)^{d}$.
  To this end, we introduce the sets
  \begin{align*}
    \grid^+_k:=\bigcap_{j\ge
      k}\grid_j\qquad\text{and}\qquad
    \mathring\grid_k^+:=\{\elm\in\gridk^+\colon\mathcal{N}^{\gridk}(\elm)\subset\gridk^+\},
  \end{align*}
  i.e., $\mathcal{N}^{\grid_j}(\elm)=\mathcal{N}^{\grid_k}(\elm)$ for all
  $j\ge k$ and $\elm\in\mathring\grid_k^+$.
  For $j\ge k$, we consider the decomposition
  \begin{align*}
    \PGdivk[j]\vecv_j = (\PGdivk[j]\vecv_j)\chi_{\Omega(\mathring\grid_k^+)}+(\PGdivk[j]\vecv_j)\chi_{\Omega(\gridk\setminus\mathring\grid_k^+)}.
  \end{align*}
  For the latter term, we have according to~\eqref{eq:interpolation} that
  \begin{align*}
    \norm[s]{(\vecv_j-\PGdivk[j]\vecv_j)\chi_{\Omega(\gridk\setminus\mathring\grid_k^+)}}&\Cleq
    \norm[\infty]{\hG[{\gridk[j]}]\chi_{\mathcal{U}(\gridk\setminus\mathring\grid_k^+)}}
    \norm[s]{\nabla
      \vecv_j}
   \le
    \norm[\infty]{\hG[{\gridk}]\chi_{\mathcal{U}(\gridk\setminus\mathring\grid_k^+)}}
    \norm[s]{\nabla
      \vecv_j}.
  \end{align*}
  Here we have used the monotonicity of the mesh-size under refinement in the last
  step.
  It follows from \cite[Corollary 4.1 and (4.15)]{MoSiVe:08} that
  $\norm[L^\infty(\Omega)]{\hG[\gridk]\chi_{\Omega(\gridk\setminus\mathring\grid_k^+)}}\to0$
  as $k\to\infty$.
  Thanks to the shape-regularity of $\grids$, this readily
  implies that
  \begin{align}\label{eq:h->0}
    \lim_{k\to\infty}\norm[L^\infty(\Omega)]{\hG[\gridk]\chi_{\mathcal{U}(\gridk\setminus\mathring\grid_k^+)}}=0.
  \end{align}
  By the compact embedding $W^{1,s}_0(\Omega)^{d}\hookrightarrow
  \hookrightarrow L^s(\Omega)^{d}$ we have that $\vecv_j\to \vec{0}$ strongly in
  $L^s(\Omega)^d$ as $j\to\infty$. Combining these observations, we deduce that for any $\epsilon>0$ there exists a $K_\epsilon>0$ such that
  \begin{align}\label{eq:G-G+}
    \norm[s]{(\PGdivk[j]\vecv_j)\chi_{\Omega(\gridk\setminus\mathring\grid_k^+)}}\le
    \epsilon\qquad\text{for all $j\ge k\ge K_\epsilon$.}
  \end{align}

  We next investigate the term $(\PGdivk[j]\vecv_j)\chi_{\Omega(\mathring\grid_k^+)}$.
  Thanks to the definition of $\mathring\grid_k^+$ and
  \eqref{eq:locality} we have
  \begin{align*}
    (\PGdivk[j]\vecv_j)|_{\Omega(\mathring\grid_k^+)}=(\PGdivk[k]\vecv_j)|_{\Omega(\mathring\grid_k^+)}\quad\text{for
      all}~j\ge k.
  \end{align*}
  Since a linear operator between two normed linear spaces is norm-continuous if and only if it is weakly continuous
(cf. Theorem 6.17 in \cite{AB}, for example,) we deduce that, for fixed $k\in\N$, we have
  \begin{align*}
    (\PGdivk[j]\vecv_j)|_{\Omega(\mathring\grid_k^+)}\weak
    \vec{0} \qquad\text{weakly
      in}~W^{1,s}(\Omega(\mathring\grid_k^+))^d~\text{as}~j\to\infty.
  \end{align*}
  By the compact embedding $W^{1,s}_0(\Omega)^{d}\hookrightarrow
  \hookrightarrow L^s(\Omega)^{d}$ this implies that
  \begin{align*}
    (\PGdivk[j]\vecv_j)\chi_{\Omega(\mathring\grid_k^+)} \to\vec{0}
    \qquad\text{strongly in}~L^s(\Omega)^d ~\text{as}~j\to\infty.
  \end{align*}
  Together with~\eqref{eq:G-G+} we have, for all $j\ge k \ge
  K_\epsilon$, that
  \begin{align*}
    {\norm[s]{\PGdivk[j]\vecv_j}}&\leq \norm[s]{(\PGdiv[j]\vecv_j)\chi_{\Omega(\gridk\setminus\mathring\grid_k^+)}}
    +\norm[s]{(\PGdivk[j]\vecv_j)\chi_{\Omega(\mathring\grid_k^+)}}
    \\
    &\le \epsilon +
    \norm[s]{(\PGdivk[j]\vecv_j)\chi_{\Omega(\mathring\grid_k^+)}}
    \to \epsilon\qquad\text{as $j\to\infty$}.
  \end{align*}
  Since $\epsilon>0$ was arbitrary, this proves the assertion.
 \end{proof}

  Next, we shall introduce a quasi-interpolation operator, which will be
  important for the treatment of the, generally non-polynomial, stress
  approximation.
  \begin{ass}\label{ass:PGS}
    We assume that for each $\grid\in\grids$ there exists a linear
    projection operator
    $\PGS:L^{1}(\Omega;\Rdds)
    \to\P_{\ell-1}(\grid;\Rdds)
    $, such that
    $\PGS$ is locally $L^1$ stable, \hbox{i.\,e.},\xspace there exists a $c>0$,
    depending on $\grid_0$, such that
    \begin{align*}
       \int_\elm\abs{\PGS \bsS}\,{\rm d}x \leq c \int_{\mathcal{U}^\grid(\elm)}\abs{\bsS}\,{\rm d}x \qquad \text{for
        all}~\bsS\in L^{1}(\Omega;\Rdds)
      .
    \end{align*}
    This implies that
    \begin{align}
      \label{eq:PDS-stab}
      \norm[s]{\PGS \bsS}\le c_s \norm[s]{\bsS}\qquad \text{for
        all}~\bsS\in L^{s}(\Omega;\Rdds)
      ,
    \end{align}
    with a constant $c_s$ depending on $\grid_0$ and $s$; compare also with \eqref{eq:W1s-stab}.
  \end{ass}

\begin{ass}[Projector $\PGQ$]\label{ass:PQ}
  We assume that for each $\grid\in\grids$ there exists a linear projection operator
  $\PGQ:L^{1}(\Omega)\to\QG$ such that, for all $s'\in(1,\infty)$,
  $\PGQ$ is locally $L^1$ stable, \hbox{i.\,e.},\xspace there exists a $c>0$,
    independent of $\grid$, such that
    \begin{align*}
      \int_\elm\abs{\PGQ q}\,{\rm d}x \leq c \int_{\mathcal{U}^\grid(\elm)}\abs{q}\,{\rm d}x \qquad \text{for
        all}~q\in L^{1}(\Omega)~\text{and all}~\elm\in\grid.
    \end{align*}

\end{ass}
  We may  argue similarly as for $\PGdiv$ to
  deduce that
  \begin{align}\label{eq:interpolationQ}
    \norm[s]{\PGQ q}\le c_s \norm[s]{q}\quad\text{and}\quad\int_\elm\abs{q-\PGQ q}^s\,{\rm d}x \leq c_s
      h_\elm^{\delta s}\abs{q}^{s}_{W^{\delta,s}(\mathcal{U}^\grid(\elm))}
    \end{align}
    for all $\elm\in\grid$ and $q\in W^{\delta,s}(\Omega)$, $\delta\in\{0,1\}$.

As a consequence of~\eqref{eq:inf-sup} and Assumption~\ref{ass:Pndiv}
(compare also with \eqref{eq:W1s-stab})
the following discrete counterpart of~\eqref{eq:inf-sup} holds; see \cite{BelBerDieRu:12}.
\begin{prpstn}[Inf-sup stability]\label{p:Dinf-sup} For all
  $s,s'\in(1,\infty)$ with $\frac1s+\frac1{s'}=1$, there exists
  a $\beta_{s}>0$,  independent of $\grid\in\grids$, such that
  \begin{align*}
    \sup_{0\neq \vecV\in\VG}\frac{\int_\Omega Q\divo \vecV\dx}{\norm[1,s]{\vecV}}\geq
    \beta_s \, \norm[s']{Q}\qquad\text{for all}\quad Q\in \QoG.
  \end{align*}
\end{prpstn}

Thanks to the above considerations, there exists a discrete Bogovski{\u\i}
operator, which has the following properties; compare also with
\cite[Corollary 9]{DieningKreuzerSueli:2013}.
\begin{crllr}[Discrete Bogovski{\u\i} operator]\label{c:Dbogovskii}
  The linear operator $\BogG:=\PGdiv\circ\mathfrak{B}:\divo \VG\to\VG$ satisfies
  \begin{align*}
    \divo (\BogG H)=H\qquad\text{and}\qquad
    \beta_s\,\norm[1,s]{\BogG H}\leq \sup_{Q\in\QG}\frac{\int_\Omega HQ\dx}{\norm[s']{Q}}
  \end{align*}
  for all $H\in \divo\VG$ and $s\in(1,\infty)$, with a positive constant
  $\beta_s,$ independent of $\grid\in\grids$.

  Moreover, let $\{\grid_k\}_{k\in\N}\subset \grids$ be a sequence of
  nested partitions of $\Omega$, i.e., $\grid_{k+1}\ge\gridk$ for all
  $k\in\N$, and let $\vecV_k\in\VG[\gridk]$ be
  such that $\vecV_k\weak \vec{0}$ weakly in $W^{1,s}_0(\Omega)^d$ as $k\to\infty$. We then have that
  \begin{align*}
    \BogG[\gridk]\divo\vecV_k \weak \vec{0}\quad\text{weakly in}~W^{1,s}_0(\Omega)^{d}~\text{as $k\to\infty$.}
  \end{align*}

\end{crllr}
\begin{proof}
  The claim follows as in \cite[Corollary
  10]{DieningKreuzerSueli:2012b} after replacing \cite[Proposition
  7]{DieningKreuzerSueli:2012b} in the proof by Proposition~\ref{prop:Pnweak} here.
\end{proof}

Upon integration by parts, it follows that
\begin{align}\label{eq:skew_sym}
  -\int_\Omega (\vecv\otimes\vecw):\nabla \vec{h} \dx = \int_\Omega
  (\vecv\otimes\vec{h}):\nabla\vecw+(\divo\vecv)(\vecw\cdot\vec{h})
  \dx
\end{align}
for all $\vecv,\vecw,\vec{h}\in
\mcD(\Omega)^d$. The last term vanishes provided  that $\divo
\vecv\equiv0$, i.e., the convection term is skew-symmetric with
respect to the second and the third argument, which implies that
\begin{align*}
  \int_\Omega (\vecv\otimes\vecv):\nabla\vecv \dx=0.
\end{align*}
It can be easily seen that this is not generally true for
finite element functions
$\vecV\in\VG$, even if
\begin{align}\label{eq:disc_divfree}
    \int_\Omega Q\divo\vecV\dx &=0\qquad\text{for all}~Q\in\QG,
\end{align}
\hbox{i.\,e.},\xspace if $\vec{V}$ is discretely divergence-free.  As in \cite{Temam:84}, we wish to ensure that the discrete counterpart
of the convection term inherits this skew-symmetry of the convection
term. To this end, we observe
from~\eqref{eq:skew_sym} that
\begin{align}\label{eq:trilin}
  \begin{split}
    -\int_\Omega (\vecv\otimes\vecw):\nabla\vec{h}
    \dx=\frac12\int_\Omega
    (\vecv\otimes\vec{h}):\nabla \vecw
    -(\vecv\otimes\vecw):\nabla \vec{h}
   \dx
    \asdefined
      \Trilin{\vecv}{\vecw}{\vec{h}}
  \end{split}
\end{align}
for all $\vecv,\vecw,\vec{h}\in W^{1,\infty}_{0,\divo}(\Omega)^d$.
 We extend this definition to $W^{1,\infty}(\Omega)^d$
in the obvious way and deduce that
\begin{align}
  \label{eq:SKEW_SYM}
  \Trilin{\vecv}{\vecv}{\vecv}=0\qquad\text{for all }\vecv\in
  W^{1,\infty}(\Omega)^d.
\end{align}

We further investigate this modified convection term for fixed
$r,r'\in(1,\infty)$ with $\frac1r+\frac1{r'}=1$; recall the definition
of $\tr$ from \eqref{eq:tr}.  We note that $\tr>1$ is equivalent to
the condition $r>\frac{2d}{d+2}$. In this case we can define its dual
$\tr'\in(1,\infty)$ by $\frac1\tr+\frac1{\tr'}=1$ and we note that the
Sobolev embedding
\begin{align}\label{eq:tr-emb}
  W^{1,r}(\Omega)^d\hookrightarrow L^{2\tr}(\Omega)^d
\end{align}
holds.  This is a crucial property in the continuous problem, which
guarantees that
\begin{align}\label{eq:cont_div0}
  \int_\Omega (\vecv\otimes\vecw):\nabla \vec{h}\dx\le c\,
  \norm[1,r]{ \vecv}\norm[1,r]{\vecw}\norm[1,\tr']{\vec{h}}
\end{align}
for all $\vecv,\vecw,\vec{h}\in
W^{1,\infty}(\Omega)^d$;
see~\cite{BulGwiMalSwi:09}. 
Because of the extension \eqref{eq:trilin} of the convection term to
functions that are not necessarily pointwise divergence-free, we have to adopt
the following stronger condition in order to ensure that the trilinear
form $\Tri[\cdot,\cdot,\cdot]$ is bounded on $W^{1,r}(\Omega)^d\times W^{1,r}(\Omega)^d\times
W^{1,\tr'}(\Omega)^d$. In particular,
let $r>\frac{2d}{d+1}$, in order to ensure that there exists an
$s\in(1,\infty)$ such that $\frac{1}{r}+\frac{1}{2\tr}+\frac1{s}=1$.
In other words, we have for $\vecv,\vecw,\vec{h}\in
W^{1,\infty}(\Omega)^d$ that
\begin{align*}
  \int_\Omega(\divo\vecv)\,(\vecw\cdot\vec{h})\dx\le \norm[r]{\divo
    \vecv}\norm[2\tr]{\vecw}\norm[s]{\vec{h}} \le c\,
  \norm[1,r]{ \vecv}\norm[1,r]{\vecw}\norm[1,\tr']{\vec{h}},
\end{align*}
with a constant $c$ depending on $r$, $\Omega$ and $d$.  Here we have used
the embeddings \eqref{eq:tr-emb} and
$W^{1,\tr'}_0(\Omega)^d\hookrightarrow L^s(\Omega)^d$.
Consequently, together
with \eqref{eq:cont_div0} we thus obtain
\begin{align}
  \label{eq:cont_trilin}
  \Trilin{\vecv}{\vecw}{\vec{h}}\le c\,
  \norm[1,r]{\vecv}\norm[1,r]{\vecw}\norm[1,\tr']{\vec{h}}.
\end{align}

In view of \eqref{eq:trilin}, for $\vecv  = (v_1, \ldots,v_d)^{\rm T}\in W^{1,r}_0(\Omega)^d$, the
convective term can be reformulated as
  \begin{align}\label{eq:frakB}
    \int_\Omega\frakB[\vecv,\vecv]\cdot \vecw\dx =
    \Trilin{\vecv}{\vecv}{\vecw},\qquad \vecw\in W^{1,\tilde r}_0(\Omega)^d,
  \end{align}
where $\frakB[\vecv,\vecv]\in L^{\tr}(\Omega)^d$ is defined by $(\frakB[\vecv,\vecv])_j = \frac{1}{2}\sum_{i=1}^d v_i \frac{\partial v_i}{\partial x_j}
+ \frac{\partial}{\partial x_i}(v_i v_j)$ for $j=1,\ldots,d$. In particular, for $\vecv=\vecV\in \VG$, we have that
$\frakB[\vecV,\vecV]\in\P_{2\ell-1}(\grid)^d$.

\begin{xmpl}\label{ex:dfree-fe}
  The following velocity-pressure pairs of finite elements satisfy Assumptions
  \ref{ass:Pndiv} and \ref{ass:PQ} for $d=2,3$ (see, e.g.,
  \cite{BelBerDieRu:12,GiraultLions:01,GiraultScott:03}):
  \begin{itemize}[leftmargin=0.5cm]
  \item The lowest order Taylor--Hood element;
  \item Spaces of continuous piecewise quadratic elements for the
    velocity and piecewise constants for the pressure; see e.g. \cite[\S VI
    Example 3.6]{BrezziFortin:91}.
  \end{itemize}
  We note that the MINI element
  and the conforming Crouzeix--Raviart Stokes element do not
  satisfy the nestedness hypothesis stated in \eqref{eq:nested}.
\end{xmpl}

\begin{rmrk} The boundedness of the trilinear form $\Trilin{\cdot}{\cdot}{\cdot}$
stated in \eqref{eq:cont_trilin} requires that $r>\frac{2d}{d+1}$. In \cite{DieningKreuzerSueli:2012b} and
\cite{DieningKreuzerSueli:2013} the set of admissible values of $r$
was the same range,  $r \in (\frac{2d}{d+2},\infty)$, as in the
existence theorem for the continuous problem in~\cite{BulGwiMalSwi:09};
however, for $r \in (\frac{2d}{d+2},\frac{2d}{d+1}]$
the finite element space for the velocity was
assumed in \cite{DieningKreuzerSueli:2012b} and
\cite{DieningKreuzerSueli:2013} to consist of pointwise divergence-free
functions, whose construction is more complicated.
For simplicity, we shall therefore confine ourselves here to the
limited range of $r>\frac{2d}{d+1}$ so as to be able to admit standard
discretely divergence-free (cf. \eqref{eq:disc_divfree}) finite
element velocity spaces.
\end{rmrk}

\subsection{The Galerkin approximation\label{ss:Galerkin}}
We are now ready to state the discrete problem. Let
$\{\VG,\QG\}_{\grid\in\grids}$ be the finite element spaces of Section
\ref{ss:fem_spaces}. 

For $n\in\N$ and $\grid\in\grids$ we call a
triple of
functions $\big(\UnG,\,P^n_\grid\big)\in\VG\times \QoG$ a Galerkin approximation of
\eqref{eq:explicit} if it satisfies
\begin{align}\label{eq:discrete}
   \begin{aligned}
      \int_\Omega\bSn(\cdot,\nablas\UnG):\nablas \vecV
      +\frakB[\UnG,\UnG]\cdot\vecV-P^n_\grid\divo \vecV\dx
      &=\int_\Omega\vecf\cdot\vecV\dx,
      \\
      \int_\Omega Q\divo\UnG\dx
      &=0,
    \end{aligned}
\end{align}
for all $\vecV\in\VG$ and $Q\in\QG$.

Restricting the test-functions to $\VoG$ the discrete problem
\eqref{eq:discrete} reduces to finding $\UnG\in\VoG$ such that
\begin{align}\label{eq:Vn0}
  \int_\Omega\bSn(\cdot,\nablas\UnG):\nablas \vecV\dx
      +\Trilin{\UnG}{\UnG}{\vecV}= \int_\Omega\vecf\cdot
        \vecV\dx
\end{align}
for all $\vecV\in\VoG$.
Thanks to~\eqref{eq:SKEW_SYM}, 
it follows from Assumption~\ref{ass:Sn}
and Korn's inequality~\eqref{eq:korn} that the nonlinear operator
defined on $\VoG$ by the left-hand side of \eqref{eq:Vn0} is coercive
and continuous
on $\VoG$.
Since the dimension of $\VG$ is finite, Brouwer's fixed point
theorem ensures the existence of a solution to~\eqref{eq:Vn0}. The
existence of a solution triple to~\eqref{eq:discrete} then follows by
the discrete inf-sup stability, Proposition~\ref{p:Dinf-sup}.
Of course, because of the
weak assumptions in the definition of the maximal monotone $r$-graph,
\eqref{eq:discrete} does not define the Galerkin approximation $(\UnG,P^n_\grid)$
uniquely.
However, supposing the axiom of choice,
for each $n\in\N$, $\grid\in\grids$, we may choose an arbitrary one among possibly
infinitely many solution triples and thus obtain
\begin{align}\label{eq:SEQ}
  \big\{\big(\UnG,P^n_\grid,\bSn(\cdot,\nablas\UnG) \big) \big\}_{n\in\N,\grid\in\grids}.
\end{align}
 From~\eqref{eq:discrete} we see that $\UnG$ is discretely
  divergence-free
  and thus, thanks to~\eqref{eq:Vn0} and~\eqref{eq:SKEW_SYM},
  we have that
  \begin{align*}
    \int_\Omega\bsS^n(\cdot,\nablas\UnG):\nablas \UnG\dx= \dual{\vecf}{\UnG}\leq
    \norm[-1,r']{\vecf}\norm[1,r]{\UnG}.
  \end{align*}
  The coercivity of $\bsS^n$ (Assumption~\ref{ass:Sn}) and Korn's inequality
  \eqref{eq:korn}
  imply
  that the sequence $\{\UnG\}_{n \in \N}$ is bounded
  in the norm of $W^{1,r}_0(\Omega)^d$, independently of
  $\grid\in\grids$ and $n\in\N$. This in turn implies,
  again by Assumption~\ref{ass:Sn}, the uniform boundedness of
  $\bsS^n(\cdot,\nablas\UnG)$ in $L^{r'}(\Omega;\Rdds)
  $. In other words, there exists a constant $c_{\vecf}>0$ depending on
  the data $\vecf$, such that
  \begin{align}\label{eq:bound}
    \norm[1,r]{\UnG}+ \norm[r']{\bsS^n(\cdot,\nablas\UnG)}\leq c_{\vecf},\qquad
    \text{for all}~ \grid\in\grids~\text{and}~n\in\N.
  \end{align}

For the sake of simplicity of the
presentation, if there is no risk of confusion, we will denote in what follows
$\bSn(\nablas\UnG)=\bSn(\cdot,\nablas\UnG)$.

\begin{rmrk}
  An alternative formulation of~\eqref{eq:Vn0}   
  is as follows: find a triple
  $\big(\UnG,\,P^n_\grid,\,\bSnG\big)\in\VG\times \QoG\times
\P_{\ell-1}(\grid;\Rdds)
$ such that
\begin{align*}
   \begin{aligned}
      \int_\Omega\bSnG:\nablas \vecV
      +\frakB[\UnG,\UnG]\cdot\vecV-P^n_\grid\,\divo \vecV\dx
      &=\int_\Omega \vecf\cdot\vecV\dx,
      \\
      \int_\Omega Q\,\divo\UnG\dx
      &=0,
      \\
      \int_\Omega \bSnG:\bsD\dx&=\int_\Omega
      \bSn(\nablas\UnG):\bsD\dx,
    \end{aligned}
\end{align*}
for all $\vecV\in\VG$, $Q\in\QG$, and
$\bsD\in\P_{\ell-1}(\grid;\Rdds)
$. Here  $\P_{\ell-1}(\grid;\Rdds)$ denotes the space of all
  piecewise polynomials of degree $\le\ell-1$ on $\grid$ with values in $\Rdds$.
In particular, if we define $\PGS:L^1(\Omega;\Rdds)
\to\P_{\ell-1}(\grid;\Rdds)
$ by
\begin{align*}
   \int_\Omega \PGS\bsS:\bsD\dx&=\int_\Omega
      \bsS:\bsD\dx,\quad\text{for
                                 all}~\bsD\in\P_{\ell-1}(\grid;\Rdds)
                                 ,
\end{align*}
then Assumption~\ref{ass:PGS} can be easily verified and
$\bSnG$ may take the role of $\PGS\bSn(\nablas\UnG)$ in the subsequent analysis.
\end{rmrk}

\subsection{Discrete Lipschitz truncation}
\label{ss:Lipschitz}
In this section we shall recall a discrete counterpart of Lipschitz truncation,
which acts on finite element spaces.
This {\em discrete Lipschitz
  truncation} is a composition of a continuous Lipschitz
truncation with a projection onto the finite element space.
The continuous Lipschitz truncation used here is based on results from
\cite{DieningMalekSteinhauer:08,BreDieFuc:12,BreDieSch:12},
which provides finer estimates than
the original Lipschitz truncation technique proposed by Acerbi and
Fusco in \cite{AcerbiFusco:88};  for details
consider~\cite{DieningKreuzerSueli:2013}.

We summarize the properties of the discrete Lipschitz truncation in
the following result. Similar results for Sobolev functions can be found
in~\cite{DieningMalekSteinhauer:08} and~\cite{BreDieFuc:12}.

\begin{prpstn}
  \label{prop:dremlip}
  Let $1< s< \infty$ and let $\{\vec{E}_k\}_{k\in\N}$ be
  a sequence such that for all $k\in\N$ we have $\vec{E}_k\in\VG[\gridk]$ for some $\gridk\in\grids$.
  In addition, assume that $\{\vec{E}_k\}_{k\in\N}\subset W^{1,s}_0(\Omega)^d$
  converges to zero weakly in $W^{1,s}_0(\Omega)^d$, as $k\to\infty$.

  Then, there exists a
  sequence $\{\lambda_{k,j}\}_{k,j\in\N}\subset\R$ with $2^{2^j} \leq
  \lambda_{k,j}\leq 2^{2^{j+1}-1}$ and Lipschitz truncated functions
  $\vec{E}_{k,j} = \vec{E}_{k,\lambda_{k,j}}$, $k,j\in\N$, with the following
  properties:
  %
  \begin{enumerate}[leftmargin=1cm,itemsep=1ex,label={\rm (\alph{*})}]
  \item \label{itm:dremlip1} $\vec{E}_{k,j}\in \VG[\gridk]$;
      \item \label{itm:dTlestW1h} $\norm[1,s]{\vec{E}_{k,j}}
    \leq c\, \norm[1,s]{\vec{E}_k}$ for $1< s \leq \infty$;
  \item \label{itm:dremlip2} $\|\nabla\vec{E}_{k,j}\|_\infty\leq
    c\,\lambda_{k,j}$;
  \item  \label{itm:dremlip3} $\vec{E}_{k,j} \to 0$ in $L^\infty(\Omega)^d$
     as $k \to \infty$;
  \item \label{itm:dremlip4} $\nabla\vec{E}_{k,j} \weak^\ast 0$ in $L^\infty(\Omega)^{d\times d}$
     as $k \to \infty$;
  \item \label{itm:dremlip5} For all $k,j \in \N$ we have
    $\norm[s]{\lambda_{k,j}\, \chi_{\{\vec{E}_k\neq\vec{E}_{k,j}\}}} \leq c\,
    2^{-\frac{j}{s}}\norm[s]{\nabla \vecE_k}$.
  \end{enumerate}
  The constants $c$ appearing in the inequalities \ref{itm:dTlestW1h}, \ref{itm:dremlip2} and \ref{itm:dremlip5} depend on
  $d$, $\Omega$, $\hat\P_\V$ and the shape-regularity of
  $\{\gridk\}_{k\in\N}$. The constants in \ref{itm:dTlestW1h} and \ref{itm:dremlip5} also depend on $s$.
\end{prpstn}

\begin{proof}
  The proof is exactly the same as that of \cite[Theorem 17 and Corollary
  18]{DieningKreuzerSueli:2012b}, replacing \cite[Proposition
  7]{DieningKreuzerSueli:2012b} by Proposition~\ref{prop:Pnweak}.
\end{proof}



\section{Error Analysis}
\label{sec:aposteriori}

\subsection{Graph approximation error}
\label{sec:graph-error}

In order to quantify the error committed in the approximation of the graph
$\mathcal{A}(x)$, $x\in\Omega$, we introduce the
following indicator. For $\bsD\in L^{r}(\Omega;\Rdds)
$, $\bsS\in
L^{r'}(\Omega;\Rdds)
$, we define
\begin{align}\label{eq:EA}
  \est_{\Ao} (\bsD,\bsS)\definedas\int_\Omega
  \inf_{(\vec{\delta},\vec{\sigma})\in\Ao(x)}\abs{\bsD-\vec{\delta}}^{r}
+\abs{\bsS-\vec{\sigma}}^{r'}\dx.
\end{align}
The following result shows that this indicator is well-defined.
\begin{prpstn}\label{prop:EstA-representation}
  Let $\bsD\in L^{r}(\Omega;\Rdds)
  $ and $\bsS\in
  L^{r'}(\Omega;\Rdds)
  $; then, the mapping
  \begin{align*}
    x\mapsto  \inf_{(\vec{\delta},\vec{\sigma})\in\Ao(x)}\abs{\bsD(x)-\vec{\delta}}^{r}
    +\abs{\bsS(x)-\vec{\sigma}}^{r'}
  \end{align*}
  is integrable. Moreover, there exist $\tilde\bsD\in L^{r}(\Omega;\Rdds)
  $ and $\tilde\bsS\in
  L^{r'}(\Omega;\Rdds)
  $ such that
  $\big(\tilde\bsD(x),\tilde\bsS(x)\big)\in \Ax(x)$ for
  a.e. $x\in\Omega$ and
  \begin{align*}
    \est_{\Ao} (\bsD,\bsS)=\int_\Omega
    \abs{\bsD-\tilde\bsD}^{r}
    +\abs{\bsS-\tilde\bsS}^{r'}\dx.
  \end{align*}
\end{prpstn}

\begin{proof}
  The first claim is an immediate consequence of the second one.
  The second assertion follows from \cite[Theorem 8.2.11]{AubinFrankowska:2009}
  by observing that the mapping
  \begin{align*}
    \Omega\times\Rdds\times\Rdds\ni\big(x;(\vec{\delta},\vec{\sigma})\big)\mapsto |\bsD(x)-\vec{\delta}|^r+|\bsS(x)-\vec{\sigma}|^{r'}
  \end{align*}
  is Carath\'eodory, i.e., $x\mapsto
  |\bsD(x)-\vec{\delta}|^r+|\bsS(x)-\vec{\sigma}|^{r'}$ is measurable
  for all $\vec{\delta},\vec{\sigma}\in\Rdds$ and $$(\vec{\delta},\vec{\sigma})\mapsto
  |\bsD(x)-\vec{\delta}|^r+|\bsS(x)-\vec{\sigma}|^{r'}$$ is continuous
  for a.e. $x\in\Omega$.
\end{proof}

\subsection{A posteriori finite element error estimates}
\label{sec:indicators}
In this section we shall prove bounds on the residual
$$\Res(\UnG,P_\grid^n,\bSn(\nablas\UnG))=(\Res^{\textsf{pde}}(\UnG,P_\grid^n,\bSn(\nablas\UnG)),\Res^{\textsf{ic}}(\UnG)) \in W^{-1,\tr}(\Omega)^{d}\times L^r_0(\Omega)$$
of \eqref{eq:explicit}. In particular, for $(\vec{v},q,\bsT)\in W^{1,r}_0(\Omega)^d\times
L^{\tr}_0(\Omega)\times L^{r'}(\Omega;\Rdds)
$ we have
\begin{align}
  \label{eq:residual}
  \begin{split}
    \dual{\Res(\vecv,q,\bsT)}{(\vecw,o)}&\definedas
    \dual{\Res^{\textsf{pde}}(\vecv,q,\bsT)}{\vecw}+
    \dual{\Res^{\textsf{ic}}(\vecv)}{o}
    \\
    &\definedas
    \int_\Omega\bsT:\nablas \vecw
    +\frakB[\vecv,\vecv]\cdot\vecw- q\divo
    \vecw-\vecf\cdot\vecw\dx
    \\
    &\quad-\int_\Omega o\divo\vecv \dx,
  \end{split}
\end{align}
where
$(\vecw,o)\in W^{1,\tr'}_0(\Omega)^d\times L^{r'}(\Omega)/\mathbb{R}$.
Although for the sake
of simplicity we restrict ourselves here
to residual-based estimates, we note that in principle other a posteriori
techniques, such as hierarchical estimates,
flux-equilibration or estimates based on local problems, can be
used as well; compare with \cite{MoSiVe:08,Siebert:11}. For $n\in\N$ and
  $\grid\in\grids$ let
$\big(\UnG,\,P^n_\grid\big)\in\VG\times \QoG$ be the Galerkin approximation defined in
\eqref{eq:discrete}. We begin with some preliminary observations.

The first part of the residual in \eqref{eq:residual}, $\Res^{\textsf{pde}}(\vecv,q,\bsT)\in W^{-1,\tr}(\Omega)^d$, provides
information about how well the functions $\vecv,q,\bsT$ satisfy the
first equation in \eqref{eq:explicit}. For the second part, we have
$\Res^{\textsf{ic}}(\vecv)\in (L_0^r(\Omega))^*$. We note that the space $(L_0^r(\Omega))^*$ is
isometrically isomorphic to $L^{r'}(\Omega)/\R$, which is, in turn, isomorphic to
$L_0^{r'}(\Omega)$ since $r\in(1,\infty)$. The term $\dual{\Res^{\textsf{ic}}(\vecv)}{o}$
provides information about the compressibility of~$\vecv$.

We emphasize that $\Res(\vecv,q,\bsT)=0$ if and only if
$\Res^{\textsf{pde}}(\vecv,q,\bsT)=0$ and $\Res^{\textsf{ic}}(\vecv)=0$,
but that  a vanishing residual itself
does not guarantee that
$\big(\nablas\vecv(x),\bsT(x)\big)\in\Ax(x)$
for almost every $x\in\Omega$. For this, additionally
$\est_{\Ao}(\nablas\vecv,\bsT)=0$ is needed.


For the rest of the paper let $t$ and $\tldt$ be such that
\begin{subequations}\label{eq:tsr}
  \begin{align}\label{eq:tsra}
    \frac{2d}{d+1}&<t< r &\quad&\text{and}&\quad \tldt &\definedas
    \frac12\frac{dt}{d-t},&\qquad&\text{if}~r\le \frac{3d}{d+2},
    \\\label{eq:tsrb}
    t&=r &\quad&\text{and}&\quad\tldt&=t'=\tr=r',&\qquad&\text{otherwise}.
  \end{align}
\end{subequations}
Note that \eqref{eq:tsra} implies that if $r \leq \frac{3d}{d+2}$, then $t<r$ and $\tldt<\tr$.

\begin{lmm}\label{l:Res=0}
The triple
  $(\vecu,p,\bsS)\in W^{1,r}_0(\Omega)^d\times
  L^{\tr}_0(\Omega)\times L^{r'}(\Omega;\Rdds)
  $ is a solution of
  \eqref{eq:implicit} if and only if
  \begin{align*}
    \Res(\vecu,p,\bsS)&=0\quad\text{in}~W^{-1,\tldt}(\Omega)^d\times
    L^{t}_0(\Omega)
    \qquad\text{and}\qquad
    \est_{\Ao}(\vecu,\bsS)=0.
  \end{align*}
\end{lmm}

\begin{proof}
  Thanks to the fact, that $W_0^{1,\tldt'}(\Omega)^d\times
  L^{t'}(\Omega)/\R$ is dense in  $W_0^{1,\tr'}(\Omega)^d\times
  L^{r'}(\Omega)/\R$ we have that
  $\Res(\vecu,p,\bsS)=0$ in $W^{-1,\tldt}(\Omega)^d\times
  L^{t}_0(\Omega)$ is equivalent to $\Res(\vecu,p,\bsS)=0$ in $W^{-1,\tr}(\Omega)^d\times
  L^{r}_0(\Omega)$. This is, in turn,
  equivalent to the fact that the triple $(\vecu,p,\bsS)$ satisfies
  the system of partial differential equations \eqref{eq:implicit}.

  On the other hand we have that $(\nablas\vecu(x),\bsS(x))\in\Ax(x)$ for
  almost every $x\in\Omega$ if and only if
  $\est_{\Ao}(\nablas\vecu,\bsS)=0$, and that completes the proof.
\end{proof}

Note that Lemma~\ref{l:Res=0} does not provide a quantitative relation
between the error and the residual. Even for simple $r$-Laplacian type
problems, such a relation requires complicated techniques and problem-adapted
error notions (e.g. a suitable quasi-norm);
cf. \cite{DieningKreuzer:08,BelenkiDieningKreuzer:12,LiuYan:01}.
However, because of the possible nonuniqueness of solutions to~\eqref{eq:implicit}, such a relation cannot
be guaranteed in our situation. We shall therefore restrict the a
posteriori analysis to bounding the residual of the problem
instead of bounding the error.

Recalling the quasi-interpolation $\PGS$ from
Assumption~\ref{ass:PGS} as well as the
representation of the discrete convective term in~\eqref{eq:frakB}, we define the
local indicators on $\elm \in\grid$ as follows:
\begin{subequations}
  \begin{alignat}{2}
    \label{eq:locEa}
    \begin{split}
    {   \est_{\grid}^{\textsf{pde}}\big(\UnG,P_\grid^{n},\bSn(\nablas\UnG);\elm\big)}
      & \definedas { \norm[\tldt,\elm]{\hG
        \big(-\divo\PGS\bSn(\nablas\UnG) +\frakB[\UnG,\UnG]
        +\nabla P_\grid^n-\vecf\big)}^{\tldt}}
       \\
       &\quad+\norm[\tldt,\partial\elm]{\hG^{1/\tldt}\jump{\PGS\bSn(\nablas\UnG)
           -\PnG\id}}^{\tldt}
       \\
       &\quad+\norm[\tldt,\elm]{\bSn(\nablas\UnG)-\PGS\bSn(\nablas\UnG)}^{\tldt},
    \end{split}
    \\
    \label{eq:locEb}
    \est_\grid^{\textsf{ic}}\big(\UnG;\elm\big)&\definedas\norm[t,E]{\divo \UnG}^t,
    \intertext{and}
    \est_{\grid}\big(\UnG,P_\grid^{n},\bSn(\nablas\UnG);\elm\big)&\definedas
    \est_{\grid}^{\textsf{pde}}\big(\UnG,P_\grid^{n},\bSn(\nablas\UnG);\elm\big)+
    \est_{\grid}^{\textsf{ic}}\big(\UnG;\elm\big).
 \end{alignat}
\end{subequations}
Here, for $S\in\sides(\grid)$, $\jump{\cdot}|_S$
denotes the normal jump across $S$ and
$\jump{\cdot}|_{\partial\Omega}:= 0$.  Moreover, we define the
error bounds to be the sums of the
local indicators, i.e., for $\mathcal{M}\subset\grid$, we have
\begin{gather*}
  \est_{\grid}^{\textsf{pde}}\big(\UnG,P_\grid^n,\bSn(\nablas\UnG);\mathcal{M}\big):=\sum_{\elm\in\mathcal{M}}
  \est_{\grid}^{\textsf{pde}}\big(\UnG,P_\grid^{n},\bSn(\nablas\UnG);\elm\big),
  \\
  \est_\grid^{\textsf{ic}}\big(\UnG;\mathcal{M}\big)\definedas\norm[t;\Omega(\mathcal{M})]{\divo \UnG}^t
  \intertext{and}
  \begin{aligned}
    \est_{\grid}\big(\UnG,P_\grid^{n},\bSn(\nablas\UnG)\big)&\definedas
    \est_{\grid}^{\textsf{pde}}\big(\UnG,P_\grid^{n},\bSn(\nablas\UnG)\big)+
    \est_\grid^{\textsf{ic}}\big(\UnG\big)
    \\
    &\definedas
    \est_{\grid}^{\textsf{pde}}\big(\UnG,P_\grid^n,\bSn(\nablas\UnG);\grid\big)+\est_\grid^{\textsf{ic}}\big(\UnG;\grid\big).
  \end{aligned}
\end{gather*}

\begin{thrm}[Upper bound on the residual] \label{Thm:bounds} Let $n\in\N$ and
  $\grid\in\grids$, and denote by
$\big(\Un,\,P^n_\grid\big)\in\VG\times \QoG$ a Galerkin approximation
of
\eqref{eq:discrete}.
We then have the following bounds:
\begin{subequations}
  \begin{align}
     \label{eq:upper}
      \norm[W^{-1,\tldt}(\Omega)]{\Res^{\textsf{pde}}\big(\UnG,P_\grid^n,\bSn(\nablas\UnG)\big)}
      &\leq C_1\,\est_{\grid}^{\textsf{pde}}\big(\UnG,P_\grid^n,\bSn(\nablas\UnG)\big)^{1/\tldt},
      \\
       \sup_{o\in L^{t'}(\Omega)/\R} \dual{\Res^{\textsf{ic}}\big(\UnG\big)}{\frac{o}{\inf_{c\in\R}\norm[t']{o-c}}}&=
      \est^{\textsf{ic}}_\grid\big(\UnG\big)^{1/t}.
     \label{eq:div_bnd}
    \end{align}
    \end{subequations}
    The constant $C_1>0$ depends only on the shape-regularity of
  $\grid$, $\tldt$, and on the dimension $d$.
\end{thrm}

\begin{proof}
  The assertions are proved using standard techniques; compare
  e.g. with  \cite{Verfuerth:96,AinsworthOden:00}.
  For the
  reader's convenience we  sketch the arguments.
  For arbitrary $(\vecv,q)\in W^{1,\tldt'}_0(\Omega)^d\times L^{t'}(\Omega)/\R$ with
  $\norm[1,\tldt']{\vecv}=\norm[t']{p}=1$ we deduce from
  \eqref{eq:discrete} that
  \begin{multline*}
    \!\!\!\!\dual{\Res^{\textsf{pde}}\big(\UnG,P_\grid^n,\bSn(\nablas\UnG)\big)}{\vecv}
    \\
    \begin{aligned}
      &= \int_\Omega\PGS\bSn(\nablas\UnG):\nablas (\vecv-\PGdiv\vecv)
      +\frakB[\UnG,\UnG]\cdot(\vecv-\PGdiv\vecv)-\vecf\cdot(\vecv-\PGdiv\vecv)\dx
      \\
      &\quad-\int_\Omega P_\grid^n\divo
        (\vecv-\PGdiv\vecv)\dx+\int_\Omega(\bSn(\nablas\UnG)-\PGS\bSn(\nablas\UnG)):\nablas (\vecv -\PGdiv\vecv) \dx.
    \end{aligned}
  \end{multline*}
  Thanks to~\eqref{eq:frakB}, local integration by parts and using H\"older's inequality,  we obtain
  \begin{multline*}
    \dual{\Res^{\textsf{pde}}\big(\UnG,P_\grid^n,\bSn(\nablas\UnG)\big)}{\vecv}
    \\
    \begin{aligned}
      &\le
      \sum_{\elm\in\grid}\Big\{\norm[\tldt,\elm]{-\divo\PGS\bSn(\nablas\UnG)
        +\frakB [\UnG,\UnG]+\nabla
        P_\grid^n-\vecf}\norm[\tldt',\elm]{\vecv-\PGdiv\vecv}
      \\
      &\quad\qquad+\frac12\norm[\tldt,\partial\elm]{\jump{\PGS\bSn(\nablas\UnG)
          -P_\grid^n\id}}\norm[\tldt',\partial
      \elm]{\vecv-\PGdiv\vecv}
      \\
      &\quad\qquad
      +\norm[\tldt,\elm]{\bSn(\nablas\UnG)-\PGS\bSn(\nablas\UnG)}\norm[\tldt',\elm]{
        \nablas\vecv} \Big\}
      \\
      &\le C\Big( \sum_{\elm\in\grid}\Big\{\norm[\tldt,\elm]{\hG
        \big(-\divo\PGS\bSn(\nablas\UnG) +\frakB[\UnG,\UnG]+\nabla
        P_\grid^n-\vecf\big) }^{\tldt}
      \\
      &\quad\qquad\qquad+\norm[\tldt,\partial\elm]{\hG^{1/\tldt}\jump{\PGS\bSn(\nablas\UnG)
          -P_\grid^n\id}}^{\tldt}
      \\
      &\quad\qquad\qquad+\norm[\tldt,\elm]{\bSn(\nablas\UnG)-\PGS\bSn(\nablas\UnG)}^{\tldt}
      \Big\}\Big)^{1/\tldt}\norm[1,\tldt']{\vecv}.
    \end{aligned}
    \end{multline*}
  Here, in the last inequality, we used the stability of
  $\PGdiv$ (see \eqref{eq:stability}), a scaled trace theorem, and the
  interpolation estimate for $\PGdiv$ in~\eqref{eq:interpolation}, as well as the finite overlapping
  of patches and a scaled trace theorem.

  To prove the bound~\eqref{eq:div_bnd}, we first deduce from
  $\int_\Omega 1\divo \UnG\dx=0$ and
  H\"older's inequality that for all $c\in\R$, we have
  \begin{align*}
    \int_\Omega o\divo \UnG\dx= \int_\Omega(o-c)\divo \UnG\dx  \le \norm[t]{\divo
      \UnG}\norm[t']{o-c}.
  \end{align*}
  Taking the infimum over all $c\in\R$ and then the supremum over all
  $o\in L^{r'}(\Omega)$ proves `$\le$' in~\eqref{eq:div_bnd}.
  In order to prove `$\ge$', we observe that
    \begin{align*}
     \est^{\textsf{ic}}_\grid\big(\UnG\big)&=\int_\Omega\divo \UnG\, |\divo
       \UnG|^{t-2}\divo\UnG\dx
      \\
     &\le\sup_{o\in L^{t'}(\Omega)/\R}
     \dual{\Res^{\textsf{ic}}\big(\UnG\big)}{\frac{o}{\inf_{c\in\R}\norm[t']{o-c}}}
    \norm[t']{|\divo
       \UnG|^{t-2}\divo\UnG}.
  \end{align*}
  Together with the definition of $\est^{\textsf{ic}}_\grid\big(\UnG\big)$ and noting that
  \begin{align*}
    \norm[t']{|\divo\UnG|^{t-2}\divo\UnG}=\norm[t]{\divo\UnG}^{t-1}
    =\est^{\textsf{ic}}_\grid\big(\UnG\big)^{1-\frac1t}=\est^{\textsf{ic}}_\grid\big(\UnG\big)^{\frac{1}{t'}},
  \end{align*}
  this yields~\eqref{eq:div_bnd}.
\end{proof}

\begin{crllr}\label{c:bounds}
  Under the conditions of Theorem~\ref{Thm:bounds}, we have
    \begin{multline*}
      \dual{\Res^{\textsf{pde}}\big(\UnG,P_\grid^n,\bSn(\nablas\UnG)\big)}{\vecv}
      \\\le C_1 \sum_{\elm\in\grid}\est_{\grid}^{\textsf{pde}}\big
      (\UnG,P_\grid^n,\bSn(\nablas\UnG);\elm\big)^{1/\tldt}
      \norm[\tldt',\mathcal{U}^\grid(\elm)]{\nabla\vecv}
    \end{multline*}
    and
    \begin{align*}
      \dual{\Res^{\textsf{ic}}\big(\UnG\big)}{q}&\le
      \sum_{\elm\in\grid}\est^{\textsf{ic}}_\grid\big(\UnG;\elm\big)^{1/t}
      \norm[t',\mathcal{U}^\grid(\elm)]{q}
    \end{align*}
  for all $\vecv\in W^{1,\tldt'}_0(\Omega)^d$ and $q\in L^{t'}(\Omega)$.
\end{crllr}

\begin{thrm}[Lower bound on the residual]\label{Thm:lbounds}
  Under the conditions of Theorem~\ref{Thm:bounds}, we have
    \begin{align} \label{eq:lower}
      \begin{split}
        c_2\,\est_{\grid}^{\textsf{pde}}&\big(\UnG,P_\grid^n,\bSn(\nablas\UnG)\big)^{1/\tldt}
        \\
        &\le
        \norm[W^{-1,\tldt}(\Omega)]{\Res^{\textsf{pde}}\big(\UnG,P_\grid^n,\bSn(\nablas\UnG)\big)}
        +\osc_\grid\big(\UnG,\bSn(\nablas\UnG)\big)^{1/\tldt}.
      \end{split}
    \end{align}
 The constant $c_2>0$ depends solely on the shape-regularity of
  $\grid$, $\tldt$, and on the dimension $d$. The oscillation term is defined by
  \begin{align*}
    \osc_\grid\big(\UnG,\bSn(\nablas\UnG)\big)&
    \definedas
    \sum_{\elm\in\grid}\osc\big(\UnG,\bSn(\nablas\UnG),\elm\big)
    \\
    &\definedas
    \sum_{\elm\in\grid}\Big\{\min_{\vecf_\elm\in\P^d_{2\ell-1}}\norm[\tldt,\elm]{\hG
      \big(       \vecf-\vecf_\elm\big) }^{\tldt}
      \\
    &\qquad+\norm[\tldt,\elm]{\bSn(\nablas\UnG)-\PGS\bSn(\nablas\UnG)}^{\tldt}\Big\}.
  \end{align*}
\end{thrm}

\begin{proof}
  Let $\elm\in\grid$ and let $S\in\sides(\grid)$, i.e.,  there exist $\elm_1,\elm_2\in\grid$, $\elm_1\neq \elm_2$, such that
  $S=\elm_1\cap \elm_2$. Let $\vecf_\elm\in\P_{2\ell-1}^d$ be arbitrary; for convenience we
  use the notation
  \begin{align*}
    R_\elm&\definedas -\divo\PGS\bSn(\nablas\UnG) +\frakB[\UnG,\UnG]+\nabla
  P_\grid^n-\vecf_\elm\in\P_{2\ell-1}^d,
  \intertext{and}
  J_S&\definedas
  \jump{\PGS \bSn(\nablas\UnG)-\PnG\id}|_S\in\P_{m}^{d\times
    d},\quad\text{where}~m=\max\{\ell-1, \jmath\}.
  \end{align*}

  It is well known that there exist local bubble functions
  $\pzb_\elm,\pzb_S\in
  W^{1,\infty}_0(\Omega)$, such that
  \begin{subequations}\label{eq:bubble}
    \begin{align}\label{eq:bubblea}
      0\le \pzb_\elm,\pzb_S\le 1,\qquad\supp \pzb_\elm =
      \elm\qquad\text{and}\qquad \supp \pzb_S = \omega_S:=\elm_1\cup
      \elm_2.
    \end{align}
    Moreover, we have that there exist $\rho_\elm\in \P_{2\ell-1}^d$
    and $\rho_S\in\P_m^{d\times d}$, with
    $\norm[\tldt',\elm]{\rho_\elm}=1=\norm[\tldt',S]{\rho_S}$, such that
    \begin{align}\label{eq:bubbleb}
      \begin{gathered}
        \norm[\tldt,\elm]{R_\elm}\le C \int_\elm R_\elm
        \pzb_\elm\rho_\elm\dx,\qquad \norm[\tldt',\elm]{\nabla(\pzb_\elm
          \rho_\elm)}\le C \norm[\tldt',\elm]{\hG^{-1}\rho_\elm},
        \\
        \norm[\tldt,S]{J_S}\le C \int_S J_S \pzb_S\rho_S\dx,\qquad
        \norm[\tldt',\omega_S]{\nabla(\pzb_S \rho_S)}\le C
        \norm[\tldt',S]{\hG^{-1/\tldt'}\rho_S},
        \\
        \text{and}\qquad \norm[\tldt',\omega_S]{\pzb_S \rho_S}\le C
        \norm[\tldt',S]{\hG^{1/\tldt'}\rho_S};
      \end{gathered}
    \end{align}
  \end{subequations}
  compare, for example, with \cite[Chapter 3.6]{Verfuerth:13}.
  Here the constants only depend on $r$, the polynomial degree
  of $R_\elm$, respectively $J_S$, on the shape-regularity of $\grid$, and on
  the dimension $d$. Hence,  for the element residual, we deduce
  that
  \begin{align*}
     \norm[\tldt,\elm]{R_\elm}
      &\le
      C\,\Big\{\dual{\Res^{\textsf{pde}}(\UnG,\PnG,\bSn(\nablas\UnG))}{\pzb_\elm
          \rho_\elm}
        \\
    &\qquad+\dual{\PGS\bSn(\nablas\UnG)-\bSn(\nablas\UnG)}{\bsD(\pzb_\elm
          \rho_\elm)}+ \norm[\tldt,\elm]{\vecf-\vecf_\elm} \Big\},
  \end{align*}
  where we have used H\"older's inequality and that $0 \leq \pzb_\elm\le1$. Together
  with a triangle inequality and~\eqref{eq:bubbleb}, this implies that
  \begin{align}\label{eq:Elest}
    \begin{split}
      \norm[\tldt,\elm]{\hG \big(-\divo\PGS\bSn(\nablas\UnG) +\frakB[\UnG,\UnG]+\nabla
        P_\grid^n-f\big)} &
      \\
      &\hspace{-6.5cm}
      \le
      C\,\Big\{\dual{\Res^{\textsf{pde}}(\UnG,\PnG,\bSn(\nablas\UnG))}{\hG\pzb_\elm
          \rho_\elm}
        \\
        &\hspace{-6.5cm}\qquad
        + \norm[\tldt,\elm]{\PGS\bSn(\nablas\UnG)-\bSn(\nablas\UnG)}
        + \norm[\tldt,\elm]{\hG(\vecf-\vecf_\elm)}\Big\}.
    \end{split}
\end{align}

For the jump residual, we deduce from~\eqref{eq:bubblea} and
integration by parts, that
 \begin{align*}
   \norm[\tldt,S]{J_S}&\le C\int_S\jump{\PGS\bSn(\nablas\UnG)-\PnG\id}\pzb_S\rho_S\ds
   \\
   &=C\Big\{\dual{\Res^{\textsf{pde}}(\UnG,\PnG,\bSn(\nablas\UnG))}{\pzb_S\rho_S}
   \\
   &\qquad+
   \sum_{i=1,2}\int_{\elm_i}\big(\divo\PGS\bSn(\nablas\UnG) -\frakB[\UnG,\UnG]-\nabla
   P_\grid^n+\vecf\big) \pzb_S\rho_S\dx\Big\}.
 \end{align*}
 Therefore, we obtain, with~\eqref{eq:bubbleb}, H\"older's inequality and~\eqref{eq:Elest}, that
\begin{align}
  \label{eq:Jest}
  \begin{split}
    \norm[\tldt,S]{\hG^{1/\tldt}\jump{\bSn(\nablas\UnG)-\PnG\id}}&
    \\
    &\hspace{-3.5cm}\le
    C\Big\{\dual{\Res^{\textsf{pde}}(\UnG,\PnG,\bSn(\nablas\UnG))}{h_S^{1/\tldt}\pzb_S\rho_S}
    \\
    &\hspace{-3.5cm}\qquad+\sum_{i=1,2}\Big[
    \dual{\Res^{\textsf{pde}}(\UnG,\PnG,\bSn(\nablas\UnG))}{\hG\pzb_{\elm_i}
      \rho_{\elm_i}}
    \\
    &\hspace{-3.5cm}\qquad\qquad
    + \norm[\tldt,\elm_i]{\PGS\bSn(\nablas\UnG)-\bSn(\nablas\UnG)}
    + \norm[\tldt,\elm_i]{\hG(\vecf-\vecf_{\elm_i})}\Big]\Big\}.
  \end{split}
\end{align}
We define the constants
$\alpha_\elm:= \norm[\tldt,\elm]{\hG \big(-\divo\PGS\bSn(\nablas\UnG) +\frakB[\UnG,\UnG]+\nabla
   P_\grid^n-\vecf\big)}^{\tldt-1}$, $\elm\in\grid$, and $\beta_S:=
 \norm[\tldt,S]{\hG^{1/\tldt}\jump{\PGS\bSn(\nablas\UnG)-\PnG\id}}^{\tldt-1}$, $S\in\sides(\grid)$. Then,
combining~\eqref{eq:Elest} and~\eqref{eq:Jest} and  summing over all
$\elm\in\grid$, $S\in\sides(\grid)$,  yields
 \begin{multline*}
   \est_{\grid}^{\textsf{pde}}\big(\UnG,P_\grid^n,\bSn(\nablas\UnG)\big)
   \\
   \begin{aligned}
     &=\sum_{\elm\in\grid}\alpha_\elm \norm[\tldt,\elm]{\hG
       \big(-\divo\PGS\bSn(\nablas\UnG) +\frakB[\UnG,\UnG]+\nabla
       P_\grid^n-\vecf\big)}
     \\
     &\quad+\sum_{s\in\sides(\grid)}\beta_S
     \norm[\tldt,S]{\hG^{1/\tldt}\jump{\bSn(\nablas\UnG)-\PnG\id}}
     \\
     &\le
     C\,\Bigg\{\dual{\Res^{\textsf{pde}}(\UnG,\PnG,\bSn(\nablas\UnG))}{\sum_{\elm\in\grid}\Big(\alpha_\elm+\sum_{S\subset\partial\elm\cap\Omega}\beta_S\Big)
       \hG\pzb_\elm \rho_\elm}
     \\
     &\qquad+\dual{\Res^{\textsf{pde}}(\UnG,\PnG,\bSn(\nablas\UnG))}{\sum_{S\in\sides(\grid)}\beta_S
       h_S^{1/\tldt}\pzb_S \rho_S}
     \\
     &\qquad+\sum_{\elm\in\grid}
     \Big(\alpha_\elm+\sum_{S\subset\partial\elm\cap\Omega}\beta_S\Big)\osc\big(\UnG,\bSn(\nablas\UnG),\elm\big)^{1/\tldt}\Bigg\}.
   \end{aligned}
   \end{multline*}
 Here we have used in the last step that $\vecf_\elm\in \P^d_{2\ell-1}$, with $\elm\in\grid$,
 are arbitrary.

 Thanks to the fact that the $\supp\pzb_\elm$,
 $\elm\in\grid$, are mutually disjoint up to a null-set, together
 with~\eqref{eq:bubbleb},  we have
 that
 \begin{multline*}
   \Big\|\sum_{\elm\in\grid}\Big(\alpha_\elm+\sum_{S\subset\partial\elm\cap\Omega}\beta_S\Big) \nabla(\hG\pzb_\elm
     \rho_\elm)\Big\|^{\tldt'}_{\tldt'}
     \\
   \begin{aligned}
     &=\sum_{\elm}\Big(\alpha_\elm+\sum_{S\subset\partial\elm\cap\Omega}\beta_S\Big)^{\tldt'}\int_\elm
     \big|\hG\nabla(\pzb_\elm \rho_\elm)\big|^{\tldt'}\dx
     \\
     &\le C
     \Big(\sum_{\elm}\alpha_\elm^{\tldt'}+\sum_{S\in\sides(\grid)}\beta_S^{\tldt'}\Big)
     \le
     C\,\est_{\grid}^{\textsf{pde}}\big(\UnG,P_\grid^n,\bSn(\nablas\UnG)\big)^{1/\tldt'},
   \end{aligned}
 \end{multline*}
 where we have used that each element $\elm$ has at most $(d+1)$ sides
 $S\in\sides(\grid)$ with $S\subset \partial\elm$. The constants $C$
 depend only on the shape-regularity of $\grid$.
Analogously, we deduce from the fact that only finitely many of the $\supp\pzb_S$,
 $S\in\sides(\grid)$, overlap, that
 \begin{align*}
   \Big\|\nabla\Big(\sum_{S\in\sides(\grid)}\beta_S \hG^{1/\tldt}\pzb_S
   \rho_S\Big)\Big\|_{\tldt'}^{\tldt'}&\le C
   \sum_{S\in\sides(\grid)}\beta^{\tldt'}\int_{\omega_S}\big|\hG^{1/\tldt}\nabla\pzb_S\rho_S\big|^{\tldt'}
   \le C \sum_{S\in\sides(\grid)}\beta^{\tldt'}
   \\
   &\le C \est_{\grid}^{\textsf{pde}}\big(\UnG,P_\grid^n,\bSn(\nablas\UnG)\big)^{1/\tldt'}.
 \end{align*}
Combining H\"older's inequality with similar arguments yields for the
last term
\begin{multline*}
  \sum_{\elm\in\grid}
   \Big(\alpha_\elm+\sum_{S\subset\partial\elm\cap\Omega}\beta_S\Big)\osc\big(\UnG,\bSn(\nablas\UnG);\elm\big)^{1/\tldt}
   \\
  \begin{aligned}
    &\le C\Big(\sum_{\elm\in\grid}\alpha_\elm^{\tldt'}+
    \sum_{S\in\sides(\grid)}\beta_S^{\tldt'}\Big)^{1/\tldt'}\osc(\UnG,\bSn(\nablas\UnG))^{1/\tldt}
    \\
    &\le C
    \,\est_{\grid}^{\textsf{pde}}\big(\UnG,P_\grid^n,\bSn(\nablas\UnG)\big)^{1/\tldt'}
    \osc(\UnG,\bSn(\nablas\UnG))^{1/\tldt}.
  \end{aligned}
 \end{multline*}
 Altogether, we have thus proved that
 \begin{multline*}
   \est_{\grid}^{\textsf{pde}}\big(\UnG,P_\grid^n,\bSn(\nablas\UnG)\big)
   \\
   \le
   C\,\Big\{\norm[W^{-1,\tldt'}(\Omega)]{\Res^{\textsf{pde}}(\UnG,\PnG,\bSn(\nablas\UnG))}\est_{\grid}^{\textsf{pde}}\big(\UnG,P_\grid^n,\bSn(\nablas\UnG)\big)^{1/\tldt'}
   \\
   +\osc(\UnG,\bSn(\nablas\UnG))
   \,\est_{\grid}^{\textsf{pde}}\big(\UnG,P_\grid^n,\bSn(\nablas\UnG)\big)^{1/\tldt'}\Big\}.
 \end{multline*}
 This is the desired bound.
\end{proof}

The following result states the local stability of the error bound and is
referred to as \textit{local lower bound} in the context of linear elliptic problems.

\begin{crllr}[Local stability]\label{c:stabEstG}
  Suppose the conditions of Theorem~\ref{Thm:bounds} and let
  $\mathcal{M}\subset\grid$; then,
  there exists a constant $C$, depending solely on the
  shape-regularity of $\grid$, $\tldt$, $d$ and $\Omega$,
  such that
  \begin{multline*}
    \est_{\grid}^{\textsf{pde}}\big(\UnG,P_\grid^{n},\bSn(\nablas\UnG);\mathcal{M}\big)^{1/\tldt}
    \\
    \begin{aligned}
      &\le C\Big(
      \norm[W^{-1,\tldt}(\mathcal{U}^\grid({\mathcal{M}}))]{\Res^{\textsf{pde}}
        (\UnG,\PnG,\bSn(\nablas\UnG))}
      +\osc(\UnG,\bSn(\nablas\UnG);\mathcal{M})^{1/\tldt}\Big)
      \\
      &\le C\Big(\norm[1,t;\mathcal{U}^\grid(\mathcal{M})]{\UnG}+
      \norm[1,t;\mathcal{U}^\grid(\mathcal{M})]{\UnG}^2+
      \norm[\tldt;\mathcal{U}^\grid(\mathcal{M})]{\PnG}+
      \norm[\tldt,\mathcal{U}^\grid(\mathcal{M})]{\vecf}+
      \norm[r,\mathcal{U}^\grid(\mathcal{M})]{\tilde  k}\Big).
    \end{aligned}
  \end{multline*}
\end{crllr}

\begin{proof}
  The first bound follows as in the proof of
  Theorem~\ref{Thm:lbounds}. In order to prove the second bound,
  let $\vecv\in
  W^{1,\tldt'}_0(\mathcal{U}^\grid(\mathcal{M}))^d$. We then have with
  H\"older's inequality and~\eqref{eq:cont_trilin}, with $t$ and $\tldt$ instead of $r$ and
  $\tr$, 
  that
  \begin{multline*}
    \dual{\Res^{\textsf{pde}}(\UnG,\PnG,\bSn(\nablas\UnG))}{\vecv}
    \\
    \begin{aligned}
      &= \int_\Omega\PGS\bSn(\nablas\UnG):\nablas \vecv
      +\frakB[\UnG,\UnG]\cdot\vecv-\PnG\divo \vecv-\vecf\cdot\vecv\dx
      \\
      &
      \begin{aligned}
        \le
        C\Big(\norm[\tldt;\mathcal{U}^\grid(\mathcal{M})]{\PGS\bSn(\nablas\UnG)}&+
        \norm[1,t;\,\mathcal{U}^\grid(\mathcal{M})]{\UnG}^2
        \\
        &+ \norm[\tldt;\,\mathcal{U}^\grid(\mathcal{M})]{\PnG}+
        \norm[\tldt;\,\mathcal{U}^\grid(\mathcal{M})]{\vecf}\Big)
        \norm[1,\tldt';\,\mathcal{U}^\grid(\mathcal{M})]{\vecv}.
      \end{aligned}
    \end{aligned}
  \end{multline*}
  Note that in the case of \eqref{eq:tsra}, we have
  \begin{align*}
    \tldt=\frac12\frac{dt}{d-t}<\frac12\frac{dt}{d-r}=\tr\frac{t}{r}\le r'\frac{t}{r}=\frac{t}{r-1}=:s'.
  \end{align*}
  Hence, H\"older's inequality,  the stability of
  $\PGS$ (Assumption~\ref{ass:PGS})  and
  Assumption~\ref{ass:Sn} yield
  \begin{align*}
    \norm[\tldt;\,\mathcal{U}^\grid(\mathcal{M})]{\PGS\bSn(\nablas\UnG)}&\le
      \abs{\mathcal{U}^\grid(\mathcal{M})}^{\frac{s'-\tldt}{s'\tldt}}
      \norm[s';\,\mathcal{U}^\grid(\mathcal{M})]{\PGS\bSn(\nablas\UnG)}
      \\
    &\le
      \abs{\Omega}^{\frac{s'-\tldt}{s'\tldt}}
      \norm[s';\,\mathcal{U}^\grid(\mathcal{M})]{\PGS\bSn(\nablas\UnG)}
      \\
      &\le C \Big(
      \norm[t;\,\mathcal{U}^\grid(\mathcal{M})]{\nablas\UnG}+
      \norm[r';\,\mathcal{U}^\grid(\mathcal{M})]{\tilde  k}\Big).
  \end{align*}
  The oscillation term can be bounded above similarly, and the assertions follows.
\end{proof}

\begin{rmrk}\label{rem:BerroneSueli}
Corollary~\ref{c:stabEstG} states the stability properties of the
estimator, which are required in order to apply the convergence theory
in \cite{Siebert:11,MoSiVe:08}; compare with
\cite[(2.10b)]{Siebert:11}, for example. The stability of the estimator is also of importance
for the efficiency of the estimator. If Corollary~\ref{c:stabEstG} fails to hold,
it may happen that the a posteriori error estimator is unbounded even though the sequence
of discrete solutions is convergent; in particular, $\divo \bSn(\nablas\UnG)$ need not
belong to $L^{r'}(\elm)$ when $1<r<2$.
This problem already appears
in the a posteriori analysis of quadratic finite element
approximations of the $r$-Laplacian, or the $r$-Stokes problem (cf. \cite{BerroneSuli:08})
for $1<r<2$.
In order to avoid this, we use $\PGS\bSn(\nablas\UnG)$ instead of
$\bSn(\nablas\UnG)$ in the element residual~\eqref{eq:locEa}.
This is compensated by the term
$\norm[\tldt]{\bSn(\nablas\UnG)-\PGS\bSn(\nablas\UnG)}^{\tldt}$ in the
a posteriori bounds~\eqref{eq:upper} and~\eqref{eq:lower};
cf. Appendix~\ref{sec:appendix} for further details. 
\end{rmrk}

\section{Convergent Adaptive Finite Elements}
\label{sec:afem+conv}
This section is concerned with the proof of convergence of an adaptive finite element algorithm for the implicit
constitutive model under consideration.

\subsection{The adaptive finite element method (AFEM)}\label{sec:AFEM}
In this section, we shall introduce an adaptive finite element
method for \eqref{eq:implicit}.
\begin{algorithm}[h]\caption{AFEM}
   \baselineskip=15pt
  \flushleft{Let $k=0$, $n_0=1$, and let $\grid_0$ be a given partition
    of $\Omega$.}

  \begin{algorithmic}[1]
    \LOOP
    \STATE let $\bSk=\bSn[n_k]$.
    \STATE $(\Uk,P_k,\bSk(\cdot,\nablas\Uk))=\SOLVE (n_k,\grid_k)$
    \STATE compute 
    $\big\{\est_{\grid_k}\big(\Uk,P_k,\bSk(\cdot,\nablas\Uk);\elm\big)\big\}_{\elm\in\grid_k}$,
    and $\est_{\Ao}(\nablas\Uk,\bSk(\cdot,\nablas\Uk))$
    \IF{$\est_{\grid_k}\big(\Uk,P_k,\bSk(\cdot,\nablas\Uk)\big)\ge\est_{\Ao}(\nablas\Uk,\bSk(\cdot,\nablas\Uk))$\label{afem:if}}
    \STATE  \label{afem:G}
    $\Mk=\MARK\big(\big\{\est_{\grid_k}
      \big(\Uk,P_k,\bSk(\cdot,\nablas\Uk);\elm\big)\big\}_{\elm\in\grid_k},\grid_k\big)$
    \STATE $\gridk[k+1]=\REFINE\big(\gridk,\Mk\big)$\hfill\% mesh-refinement
    \STATE $n_{k+1}=n_k$
    \ELSE
    \STATE $n_{k+1}=n_k+1$\label{afem:nk}\hfill\% graph-refinement
    \ENDIF
    \STATE $k=k+1$
    \ENDLOOP
  \end{algorithmic}
\end{algorithm}
The details of the subroutines used in the process are listed below:

\paragraph{The routine \SOLVE}
We assume that for arbitrary $n\in\N$, $\grid\in\grids$, the
  routine $\SOLVE(n,\grid)\!=\!(\UnG,P_n,\bSn(\cdot,\nablas\UnG))$
  computes an exact solution $(\UnG,P_n)\in\VG\times \QG$ of
  \eqref{eq:discrete}.

\paragraph{The routine \MARK}
For a fixed function $g:\R^+_0\to\R^+_0$, which is continuous at $0$
with $g(0)=0$, we assume that the set $\marked=\MARK\big(\big\{\est_{\grid}
\big(\UnG,\PnG,\bsS^n(\nablas\UnG);\elm\big)\big\}_{\elm\in\grid},\grid\big)$
satisfies
\begin{align}\label{eq:marking}
  \begin{split}
    &\max\big\{\est_{\grid}
    \big(\UnG,\PnG,\bsS^n(\cdot,\nablas\UnG);\elm\big):\elm\in\grid\setminus\marked\big\}
    \\
    &\hspace{3cm}\le g\big(\max\big\{\est_{\grid}
    \big(\UnG,\PnG,\bsS^n(\cdot,\nablas\UnG);\elm\big):\elm\in\marked\big\}\big).
  \end{split}
\end{align}
Hence the marking criterion guarantees that all indicators in $\grid$
are controlled by the maximal indicator in $\marked$. Note that this
criterion covers most commonly used marking strategies with $g(s)=s$;
cf. \cite{Siebert:11,MoSiVe:08}.

For the definition of \textbf{the routine \REFINE} see Section~\ref{sec:grid}.

For the sake of simplicity of the
presentation, in the following, we will suppress the dependence on $x$
in our notation and
write $\bSk(\nablas\Uk)=\bSk(\cdot,\nablas\Uk)$ if there is no risk of confusion.

\subsection{Convergence of the AFEM}
\label{sec:convergence}
Let $\{\grid_k\}_{k\in\N}\subset\grids$ be the sequence of meshes
produced by AFEM.
For $s\in(1,\infty]$, we define
\begin{align}\label{df:Vinfty}
  \V_\infty^s\definedas
  \overline{\bigcup_{k\ge0}\VG[\grid_k]}^{\norm[1,s]{\cdot}}\subset W^{1,s}_0(\Omega)^d\quad\text{and}\quad
  \Q^s_\infty\definedas \overline{\bigcup_{k\ge0}\QG[\grid_k]}^{\norm[s]{\cdot}}\subset
  L^{s}_0(\Omega).
\end{align}

\begin{lmm}\label{l:conv}
   Let $\big\{\big(\Uk,P_k,\bSk(\nablas\Uk)\big)\big\}_{k\in\N}\subset
  W_0^r(\Omega)^d\times L^{\tr}_0(\Omega)\times L^{r'}(\Omega;\Rdds)
  $ be
  the sequence produced by AFEM; then, at least for a not relabelled subsequence, we have
  \begin{alignat*}{2}
    \Uk&\weak \vecu_\infty &\qquad&\text{weakly in}~W^{1,r}_0(\Omega)^d,
    \\
    P_k&\weak p_\infty&\qquad&\text{weakly in}~L^{\tr}_0(\Omega),
    \\
    \bSk(\nablas\Uk)&\weak \bsS_\infty &\qquad&\text{weakly in}~L^{r'}(\Omega;\Rdds)
    ,
  \end{alignat*}
  for some $(\vecu_\infty,p_\infty,\bsS_\infty)\in \V_\infty^r\times
  \Q_\infty^{\tr}\times L_0^{r'}(\Omega)$. Moreover, we have that
  \begin{align*}
    \Res(\Uk,P_k,\bsS_k(\nablas\Uk))\weak^*\Res(\vecu_\infty,p_\infty,\bsS_\infty)\qquad\text{weakly*
    in}~W^{-1,\tr}(\Omega)^d
  \end{align*}
  and
  \begin{align*}
    \dual{\Res(\vecu_\infty,p_\infty,\bsS_\infty)}{(\vecv,q)}=0\qquad\text{for all}~q\in\Q_\infty^{r'},\vecv\in\V_\infty^{\tr'}.
  \end{align*}
\end{lmm}

\begin{proof}
The proof is postponed to Section~\ref{ssec:Proof-conv}.
\end{proof}

\begin{crllr}\label{c:estG}
  Let $\big\{\big(\Uk,P_k,\bSk(\nablas\Uk)\big)\big\}_{k\in\N}\subset
  W_0^r(\Omega)^d\times L^{\tr}_0(\Omega)\times L^{r'}(\Omega;\Rdds)
  $ be a
  not
  relabelled subsequence with weak limit  $(\vecu_\infty,p_\infty,\bsS_\infty)\in \V_\infty^r\times
  \Q_\infty^{\tr}\times L^{r'}(\Omega;\Rdds)
  $ as in Lemma~\ref{l:conv}.
  Then,
  \begin{align*}
    \est_{\gridk}(\Uk,P_k,\bSk(\nablas\Uk))\to 0,\quad
    \text{as}~k\to\infty;
  \end{align*}
  implies that
  \begin{align*}
    \Res(\vecu_\infty,p_\infty,\bsS_\infty)=0\in W^{-1,\tr}(\Omega)^d.
  \end{align*}
\end{crllr}

\begin{proof}
  The upper bound, Theorem~\ref{Thm:bounds}, together with
  $\est_{\gridk}(\Uk,P_k,\bSk(\nablas\Uk))\to 0$ as $k\to\infty$,
  implies that
  \begin{align*}
    \Res(\Uk,P_k,\bSk(\nablas\Uk))\to 0\quad\text{strongly in}~ W^{-1,\tldt}(\Omega)^d.
  \end{align*}
  Thus the assertion follows from
  Lemma~\ref{l:conv} and the uniqueness of the limit.
\end{proof}



\begin{lmm}\label{l:estA}  Let $\big\{\big(\Uk,P_k,\bSk(\nablas\Uk)\big)\big\}_{k\in\N}\subset
  W_0^r(\Omega)^d\times L^{\tr}_0(\Omega)\times L^{r'}(\Omega;\Rdds)
  $ be a
  not relabelled subsequence with weak limit  $(\vecu_\infty,p_\infty,\bsS_\infty)\in \V_\infty^r\times
  \Q_\infty^{\tr}\times L^{r}(\Omega;\Rdds)
  $ as in Lemma~\ref{l:conv}.
  Assume that
  \begin{align*}
    \est_{\Ao}(\nablas\Uk,\bSk(\nablas\Uk))\to 0\quad
    \text{as}~k\to\infty;
  \end{align*}
  then,
  \begin{align*}
    (\nablas\vecu_\infty(x),\bsS_\infty(x))\in\Ao(x)\quad\text{  for almost every $x\in\Omega$.}
  \end{align*}

\end{lmm}

\begin{proof}
  The proof of this lemma is postponed to Section~\ref{ssec:Proof-estA} below.
\end{proof}

\begin{lmm}\label{l:estGconv}
  Assume that the sequence $\{n_k\}_{k\in\N}$ satisfies $n_k\to
  N<\infty$ as $k\to\infty$. We then have that
  \begin{align*}
    \est_{\grid_k}\big(\Uk,P_k,\bSk(\nablas\Uk))\to0\qquad\text{as}~k\to\infty.
  \end{align*}
\end{lmm}
\begin{proof}
  The proof of this lemma is postponed to Section \ref{ssec:Proof-estGconv} below.
\end{proof}

We further assume that the graph approximation is uniform with respect
to the graph approximation indicator.
\begin{ass}\label{ass:estAconv}
  For every $\epsilon>0$, 
  there exists an $N=N(\epsilon)\in\N$, such that
  \begin{align*}
    \est_{\Ao}(\nablas\vecv,\bSn(\cdot,\nablas\vecv))<\epsilon \qquad\text{for
      all}~\vecv\in W^{1,r}_0(\Omega)^d~\text{and}~n>N.
  \end{align*}
\end{ass}
We note that this and~Assumption~\ref{ass:Sn} are the only strong assumptions
among the ones we have made; Assumption~\ref{ass:estAconv} is, however, only used in the proof of
the next theorem, and is not required for any of the preceding results.

\begin{thrm}\label{Thm:AFEM}
  Suppose that Assumption~\ref{ass:estAconv} holds and let
  $\{(\Uk,P_k,\bSk(\nablas\Uk))\}$ be the sequence of function triples
  produced by the AFEM. We then have that
  \begin{align*}
    \est_{\Ao}(\nablas\Uk,\bSk(\nablas\Uk))&\to0 \qquad \text{as}~k\to\infty
    \intertext{and, for a not relabelled subsequence, we have
      that}
    \est_{\grid_k}(\Uk,P_k,\bSk(\nablas\Uk))&\to0\qquad \text{as}~k\to\infty.
  \end{align*}

\end{thrm}

\begin{proof}
  We argue by contradiction. First assume that there exists an $\epsilon>0$ such
  that, for some subsequence, we have that
  \begin{align*}
    \est_{\Ao}(\nablas\Uk[k_\ell],\bSk[k_\ell](\nablas\Uk[k_\ell]))>\epsilon\qquad\text{for all}~\ell\in\N.
  \end{align*}
  Consequently, by Assumption \ref{ass:estAconv}, we have that
  $n_{k_\ell}= N$, for some $\ell_0,N\in\N$, and
  all $\ell\ge \ell_0$. Moreover, thanks to Lemma~\ref{l:conv}, there
  exists  a not relabelled subsequence
 $\big\{\big(\Uk[k_\ell],P_{k_\ell},\bSk[k_\ell](\nablas\Uk[k_\ell])\big)\big\}_{\ell\in\N}$ that
 converges weakly  in
 $W_0^r(\Omega)^d\times L^{\tr}_0(\Omega)\times
  L^{r'}(\Omega;\Rdds)
  $. Combining these facts, we deduce with Lemma~\ref{l:estGconv}
   that
  \begin{align*}
    \est_{\grid_{k_\ell}}(\Uk[k_\ell],P_{k_\ell},\bSk[k_\ell](\nablas\Uk[k_\ell]))\to0.
  \end{align*}
  In particular, there exists an $\ell>\ell_0$, such that
  $\est_{\grid_{k_\ell}}(\Uk[k_\ell],P_{k_\ell},\bSk[k_\ell](\nablas\Uk[k_\ell]))<\epsilon$. Therefore, by
  line~\ref{afem:nk} of AFEM we have that $n_{k_\ell+1}=N+1$, a
  contradiction.
  Consequently, we have (for the full sequence) that
  \begin{align*}
    \est_{\Ao}(\nablas\Uk,\bSk(\nablas\Uk))\to0\qquad\text{as}~k\to\infty.
  \end{align*}
  This proves the first claim.

  Assume now that there exists an $\epsilon>0$ such
  that we have that
  \begin{align}\label{eq:17}
    \est_{\grid_k}(\Uk,P_k,\bSk(\nablas\Uk))>\epsilon\qquad\text{for all}~k\in\N.
  \end{align}
  By the above considerations, there exists a $k_0\in\N$ such that
  $\est_{\Ao}(\nablas\Uk,\bSk(\nablas\Uk))<\epsilon$ for all $k\ge k_0$. Therefore, according to
  line~\ref{afem:if} of AFEM,
  we have that $n_k=n_{k_0}$ for all $k\ge k_0$.
  Consequently, Lemma~\ref{l:estGconv} contradicts~\eqref{eq:17}.

  Combining the two cases proves the assertion.
\end{proof}

\begin{crllr}\label{C:convergence}
  Let
  $\{(\Uk,P_k,\bSk(\nablas\Uk))\}$ be the sequence of function triples
  produced by the AFEM. Then, there exists a not relabelled
  subsequence with weak limit  $(\vecu_\infty,p_\infty,\bsS_\infty)\in W^{1,r}_0(\Omega)^d\times
  L_0^{\tr}(\Omega) \times L^{r'}(\Omega;\Rdds)
  $ such that
  \begin{align*}
    \est_{\Ao}(\nablas\Uk,\bSk(\nablas\Uk))\to0 \qquad
    \text{and}\qquad
    \est_{\grid_k}(\Uk,P_k,\bSk(\nablas\Uk))&\to0,
  \end{align*}
  $\text{as}~k\to\infty$ and $(\vecu_\infty,p_\infty,\bsS_\infty)$
  solves~\eqref{eq:implicit}.
\end{crllr}

\begin{proof}
  The claim follows from Theorem~\ref{Thm:AFEM}, Lemma~\ref{l:estA},
  Corollary~\ref{c:estG}, and Lemma~\ref{l:Res=0}.
\end{proof}

\begin{rmrk} \label{R:filter-conv}
  We emphasize that even in the case when the exact solution of
  \eqref{eq:implicit} is unique, we do not have that the statement of
  Corollary~\ref{C:convergence} is true for the full sequence. This is
  due to the fact that the finite element error estimator is
  not necessarily decreasing with respect to the refinement of the graph
  approximation. However,
  when the exact solution is unique, it is easy to select a
  converging subsequence with the help of the estimators; one can choose, for example, a subsequence,
  such that $
    \est_{\grid_k}(\Uk[k_\ell],P_{k_\ell},\bSk[k_\ell](\nablas\Uk[k_\ell]))
    $ is monotonic decreasing in $\ell$.
  \end{rmrk}

\section{The Proofs of the Auxiliary Results }
\label{sec:aux}

\subsection{Proof of Lemma~\ref{l:conv}}
\label{ssec:Proof-conv}
  We recall~\eqref{eq:bound} and observe that
  the spaces $W^{1,r}_0(\Omega)^d$ and
  $L^{r'}(\Omega;\Rdds)
  $, $r\in(1,\infty)$, are
  reflexive. Therefore,
  there exist $\vecu_\infty\in\V_\infty^r$ and  $\bsS_\infty\in
  L^{r'}(\Omega;\Rdds)
  $ such that for a not
  relabelled subsequence we have
  \begin{align}
    \Uk &\weak \vecu_\infty\qquad\text{weakly in}\quad
    W^{1,r}_0(\Omega)^d\label{eq:weakUn} \intertext{and}
    \bsS_k(\nablas\Uk)&\weak\bsS_\infty\qquad\text{weakly in}\quad
    L^{r'}(\Omega;\Rdds)
                        ,\label{eq:weakSn}
  \end{align}
  as      $k\to\infty$.
  The function $\vecu_\infty$ is discretely divergence-free with respect to
  $\Q_\infty^{r'}$, i.e.,
  \begin{align*}
    \int_\Omega q\divo \vecu_\infty\dx=\lim_{k\to\infty}\int_\Omega(\PGQ[\gridk] q)\divo\Uk\dx=0\qquad\text{for all}~q\in\Q_\infty^{r'}.
  \end{align*}
  This follows from~\eqref{eq:interpolationQ} as in the proof of \cite[Lemma
  19]{DieningKreuzerSueli:2013},
  replacing \cite[(3.5)]{DieningKreuzerSueli:2013} with
  the density of the union of the discrete pressure spaces in $\Q_\infty^{r'}$.

  Moreover, using compact embeddings of Sobolev
  spaces, we have that
  \begin{align}\label{eq:strongUn}
    \Uk\rightarrow \vecu_\infty\qquad\text{strongly in}\quad
    L^s(\Omega)^d\quad\text{for all}\quad
    \begin{cases}
      s\in\left(1,\tfrac{rd}{d-r}\right), \quad&\text{if}~ r<d,
      \\
      s\in(1,\infty),\quad&\text{otherwise}.
    \end{cases}
  \end{align}
  This implies, for arbitrary $\vecv\in W^{1,\infty}_0(\Omega)^d$, that
  \begin{align*}
    \Trilin{\Uk}{\Uk}{\vecv}&\to
    \Trilin{\vecu_\infty}{\vecu_\infty}{\vecv},\intertext{or
      equivalently,}
    \frakB[\Uk,\Uk]&\to\frakB[\vecu_\infty,\vecu_\infty]\quad\text{weakly in}~W^{-1,1}(\Omega)^d
  \end{align*}
  as $k\to\infty$; compare also with \cite[Lemma 19]{DieningKreuzerSueli:2013}.

  We now prove convergence of the pressure. Thanks to
  \eqref{eq:cont_trilin},  
  we have
  \begin{align*}
    \int_\Omega P_k\divo
      \vecV\dx&=\int_\Omega\bsS_k(\nablas\Uk):\nablas\vecV+\frakB[\Uk,\Uk]\cdot\vecV -\vecf\cdot\vecV\dx
    \\
    &\le \norm[r']{\bsS_k(\nablas\Uk)}\norm[r]{\nablas\vecV}+
    c\,\norm[1,r]{\Uk}^2\norm[1,\tr']{\vecV} +
    \norm[-1,r']{\vecf}\norm[1,r]{\vecV}
  \end{align*}
  for all $\vecV\in\VGk$.  By~\eqref{eq:bound} and the discrete inf-sup
  condition stated in Proposition~\ref{p:Dinf-sup}, it follows that the sequence
  $\{P_k\}_{k\in\N}$ is bounded in the reflexive Banach space
  $L^{\tr}_0(\Omega)$. Hence, there exists a $p_\infty\in
  \Q_\infty^{\tr}\subset L^{\tr}_0(\Omega)$ such that, for a (not relabelled) subsequence,
  \begin{align*}
    P_k\weak p_\infty\qquad\text{weakly in}~L^{\tr}_0(\Omega).
  \end{align*}
On the other
  hand we deduce for an arbitrary $\vecv\in \V^\infty_\infty\subset W^{1,\infty}_0(\Omega)^d$
  that
  \begin{align*}
    \int_\Omega p_\infty\divo\vecv\dx&\leftarrow\int_\Omega P_k\divo \vecv\dx
    =\int_\Omega P_k\divo\PGdivk\vecv\dx
    \\
    &\quad= \int_\Omega\bsS_k(\nablas\Uk):\nablas \PGdivk\vecv-\vecf\cdot\PGdivk\vecv +\frakB[\Uk,\Uk]\cdot\PGdivk \vecv\dx
    \\
    &\to \int_\Omega \bsS_\infty:\nablas\vecv\dx+\frakB[\vecu_\infty,\vecu_\infty]\cdot{\vecv}-\vecf\cdot\vecv\dx
  \end{align*}
  as $k\to\infty$, where we have
  used~\eqref{eq:weakSn}, 
  the properties of $\PGdivk$ together with the density of the union
  of the discrete velocity spaces in $\V_\infty^\infty$ and the boundedness of the sequence
  $\{P_k\}_{k\in\N}$ in~$L^{\tr}_0(\Omega)$. The assertion for all
  $\vecv\in \V_\infty^{\tr'}$ then follows from the density of
  $\V_\infty^\infty$ in $\V_\infty^{\tr'}$. \qed

\subsection{Proof of Lemma~\ref{l:estA}}
\label{ssec:Proof-estA}
  According to Proposition~\ref{prop:EstA-representation}, for $k\in\N$ there
  exist $\tilde\bsD_k\in L^{r}(\Omega;\Rdds)
  $ and
  $\tilde\bsS_k\subset  L^{r'}(\Omega;\Rdds)
  $, such that
  $(\tilde \bsD_k(x), \tilde\bSk(x))\in\Ao(x)$ for a.e.
  $x\in\Omega$ and
  \begin{align}\label{eq:4}
    \norm[r]{\nablas\Uk-\tilde\bsD_k}^r+
    \norm[r']{\bSk(\nablas\Uk)-\tilde\bsS_k}^{r'}= \est_{\Ax}(\nablas\Uk,\bSk(\nablas\Uk)) \to0\quad\text{as}~k\to\infty.
  \end{align}
  Thanks to \eqref{eq:bound}, the sequences $\{\tilde
  \bsD_k\}_{k\in\N}$ and $\{\tilde
  \bsS_k\}_{k\in\N}$ are bounded in
  $L^r(\Omega;\Rdds)
  $ and $L^{r'}(\Omega;\Rdds)
  $
  respectively. Since both spaces are reflexive, together with the
  uniqueness of the limit, we obtain that
  \begin{subequations}\label{eq:weaktilde}
    \begin{alignat}{2}\label{eq:weaktildea}
      \tilde \bsD_k&\weak \bsD\vecu_\infty&\qquad&\text{weakly
        in}~L^{r}(\Omega;\Rdds)
      ,
      \\\label{eq:weaktildeb}
      \tilde\bsS_k&\weak \bsS_\infty&\qquad&\text{weakly
        in}~L^{r'}(\Omega;\Rdds)
    \end{alignat}
  \end{subequations}
  as $k\to\infty$.

  Let $\bsS^\star:\Omega\times\Rdds\to\Rdds$ be a measurable selection
  with $(\vec{\delta},\bsS^\star(x,\vec{\delta}))\in \Ax(x)$ for
  a.e. $x\in\Omega$ and thus
  $(\nablas\vecu_\infty(x),\bsS^\star(x,\nablas\vecu_\infty(x)))\in\Ao(x)$
   for a.e. $x\in\Omega$; compare with
   Remark~\ref{rem:Aprop}.
   Consequently, for every bounded sequence
   $\{\phi_k\}_{k\in\N}\in
   L^\infty(\Omega)$ of nonnegative functions, we have (recall \ref{A3}) that
 \begin{align}\label{eq:5}
   \begin{split}
     0&\le \limsup_{k\to\infty}\int_\Omega \abs{\tilde\bsS_k-
     \bsS^\star(\nablas\vecu_\infty)\big):(\tilde\bsD_k-\bsD\vecu_\infty)}\phi_k\dx
     \\
     &=\limsup_{k\to\infty}\int_\Omega \big(\tilde\bsS_k-
     \bsS^\star(\nablas\vecu_\infty)\big):(\tilde\bsD_k-\bsD\vecu_\infty)\,\phi_k\dx
     \\
       &=\limsup_{k\to\infty}\int_\Omega\underbrace{
         (\bSk(\nablas\Uk)-\bsS^\star(\nablas\vecu_\infty)):(\nablas\Uk-\bsD\vecu_\infty)}_{=:a_k(x)}\phi_k\dx,
     \end{split}
 \end{align}
 where we have used  \eqref{eq:4} in the last step.
 We assume for the moment that $a_k\to 0$ in
  measure and therefore
  \begin{align}\label{eq:ak->0}
    a_k\to 0 \quad\text{a.e. in}~\Omega
 \end{align}
 for at least a subsequence of $a_k$. Since $a_k$ is bounded in
 $L^1(\Omega)$, we obtain with the biting Lemma (Lemma~\ref{l:biting})
 and Vitali's theorem, that there exists a nonincreasing sequence
 of measurable subsets $E_j\subset\Omega$ with
  $\abs{E_j}\to 0$ as $j\to\infty$, such that for all $j\in \N$, we
  have
  \begin{align*}
    a_k\to 0 \qquad\text{strongly in}~L^1(\Omega\setminus E_j)\quad\text{as}~k\to\infty.
  \end{align*}
  This, together with \eqref{eq:5} and \eqref{eq:weaktilde},
  implies for all nonnegative $\phi\in
  L^\infty(\Omega\setminus E_j)\subset L^\infty(\Omega)$ (extend
  $\phi$ by zero on $E_j$) and each
  fixed $j\in\N$,
  that
  \begin{align*}
    \lim_{k\to\infty}\int_{\Omega\setminus E_j}\tilde\bsS_k:\tilde\bsD_k\phi\dx=
    \int_{\Omega\setminus E_j}\bsS_\infty:\bsD\vecu_\infty\phi\dx.
  \end{align*}
  Consequently, since the graph is monotone and $(\tilde \bsD_k(x), \tilde\bSk(x))\in\Ao(x)$ for a.e.
  $x\in\Omega$, we observe for
  arbitrary $\vec{\delta}\in\Rdds$ and all nonnegative $\phi\in
  L^\infty(\Omega\setminus E_j)$, that
  \begin{align*}
    0&\le \lim_{k\to\infty}\int_{\Omega\setminus
      E_j}\big(\tilde\bsS_k-\bsS^\star(\cdot,\vec{\delta})\big):(\tilde\bsD_k-\vec{\delta})\phi\dx 
    \\
    &=
    \int_{\Omega\setminus
      E_j}\big(\bsS_\infty -
    \bsS^\star(\cdot,\vec{\delta})\big):(\bsD\vecu_\infty-\vec{\delta})\phi\dx.
  \end{align*}
  Since $\phi$ was arbitrary, we have that
  \begin{align*}
    \big(\bsS_\infty -
    \bsS^\star(\cdot,\vec{\delta})\big):(\bsD\vecu_\infty-\vec{\delta})\ge0\qquad\text{for
      all $\vec{\delta}\in\Rdds$ and a.e. $x\in\Omega\setminus E_j$.}
  \end{align*}
  According to Lemma~\ref{lem:Sstar},
    this implies that
  \begin{align*}
    (\nablas\vecu_\infty(x),\bsS_\infty(x))\in\Ao(x)\quad\text{  for
      almost every $x\in\Omega\setminus E_j$.}
  \end{align*}
  The assertion then follows from $\abs{E_j}\to0$ as $j\to\infty$.

  It remains to verify that $a_k\to 0$ in measure as $k\to\infty$. We divide the proof into four
  steps.

  \textbf{Step 1:} First, we introduce some preliminary facts
  concerning discrete Lipschitz truncations. For
  convenience we use the notation
  \begin{align*}
    \vecE_k \definedas \PGdivk(\Uk-\vecu_\infty)=\Uk-\PGdivk\vecu_\infty\in\VG[\gridk]
  \end{align*}
  and let $\{\vecE_{k,j}\}_{k,j\in\N}$ be the sequence of
  Lipschitz-truncated  finite element functions
  according to
  Proposition~\ref{prop:dremlip}.  Recall from Lemma~\ref{l:conv}
  that $\vecE_k \weak 0$ weakly in $W^{1,r}_0(\Omega)^d$, \hbox{i.\,e.},\xspace we are
  exactly in the situation of Proposition~\ref{prop:dremlip}.
  Although $\vecE_k\in\VoG[\gridk]$, \hbox{i.\,e.},\xspace $\vecE_k$ is discretely divergence-free,
  this does not necessarily imply that $\vecE_{k,j}\in\VoG[\gridk]$
  and thus we need to modify $\vecE_{k,j}$ in order to be able
  to use it as a test function in \eqref{eq:Vn0}. With the
  discrete Bogovski{\u\i} operator $\Bogk:=\BogG[\gridk]$ from
  Corollary~\ref{c:Dbogovskii}, we define
  \begin{subequations}\label{df:PhiPsi}
    \begin{align}
      \vec{\Psi}_{k,j}\definedas \Bogk(\divo \vecE_{k,j})\in \VG[\gridk].
    \end{align}
    The `corrected' function
    \begin{align}
      \vec{\Phi}_{k,j}\definedas \vecE_{k,j}-\vec{\Psi}_{k,j}\in\VoG[\gridk]
    \end{align}
  \end{subequations}
  is then discretely divergence-free. We need to control the correction in
  a norm.  To this end we recall from Section \ref{ss:fem_spaces} that
  $\QG[\gridk]=\operatorname{span}\{Q^k_1,\ldots,Q^k_{N_k}\}$ for a certain locally
  supported basis. Then, thanks to properties of the discrete Bogovski{\u\i}
  operator and Corollary~\ref{c:Dbogovskii}, we have that
    \begin{align*}
      \beta_r\norm[1,r]{\vec{\Psi}_{k,j}}&\leq
      \sup_{Q\in\QG[\gridk]}\frac{\int_\Omega Q\divo \vecE_{k,j}\dx}{\norm[r']{Q}}
      =\sup_{Q\in\QG[\gridk]}\frac{\int_\Omega Q\divo \vecE_{k,j}-\divo
          \vecE_k\dx}{\norm[r']{Q}}\\
      &=\sup_{Q=\sum_{i=1}^{N_k}\rho_iQ_i^k}\Bigg( \sum_{\supp
        Q_i^k\subset \{\vecE_{k,j}=\vecE_{k}\}} \frac{\int_\Omega \rho_iQ_i^k\divo(
          \vecE^{k,j}- \vecE_k)\dx}{\norm[r']{Q}}\\
      &\hspace{2.5cm}+ \sum_{\supp
        Q_i^k\cap\{\vecE_{k,j}\neq\vecE_{k}\}\neq\emptyset}
      \frac{\int_\Omega \rho_iQ_i^k\divo (\vecE_{k,j}-
          \vecE_k)\dx}{\norm[r']{Q}}\Bigg)\\
      &=\sup_{Q=\sum_{i=1}^{N_k}\rho_iQ_i^n}\Bigg( \sum_{\supp
        Q_i^k\cap\{\vecE_{k,j}\neq\vecE_{k}\}\neq\emptyset}
      \frac{\int_\Omega \rho_iQ_i^k\divo (\vecE_{k,j}-
          \vecE_k)\dx}{\norm[r']{Q}}\Bigg)
      \\
      &=\sup_{Q=\sum_{i=1}^{N_k}\rho_iQ_i^k}\Bigg( \sum_{\supp
        Q_i^k\cap\{\vecE_{k,j}\neq\vecE_{k}\}\neq\emptyset}
      \frac{\int_\Omega \rho_iQ_i^k\divo \vecE_{k,j}\dx}{\norm[r']{Q}}\Bigg)
      \\
      &\le \norm[r]{\divo
        \vecE_{k,j}\chi_{\Omega^k_{\{\vecE_{k,j}\neq\vecE_k\}}}}\sup_{Q=\sum_{i=1}^{N_k}\rho_iQ_i^k}\frac{\norm[r']{\sum_{\supp
        Q_i^k\cap\{\vecE_{k,j}\neq\vecE_{k}\}\neq\emptyset}\rho_iQ_i^k}}{\norm[r']{Q}}
      \\
      &
        \leq  c\,\norm[r]{\divo
        \vecE_{k,j}\,\chi_{\Omega^k_{\{\vecE_{k,j}\neq\vecE_k\}}}}\le
      c\,\norm[r]{\nabla
        \vecE_{k,j}\,\chi_{\Omega^k_{\{\vecE_{k,j}\neq\vecE_k\}}}},
  \end{align*}
  where $\chi_{\Omega^k_{\{\vecE_{k,j}\neq\vecE_k\}}}$ is the characteristic function of the
  set
  \begin{align*}
    \Omega^k_{\{\vecE_{k,j}\neq\vecE_k\}}\definedas\bigcup\left\{\Omega^k_\elm\mid
      \elm\in\gridk~\text{such that}~\elm\subset
      \overline{\{\vecE_{k,j}\neq\vecE_k\}} \right\}.
  \end{align*}
  Note that in the penultimate step of the above estimate we have used
  norm equivalence
  on the reference space $\hat\P_\Q$ from \eqref{df:Qn}. In
  particular, we see by means of standard scaling arguments that for
  $Q=\sum_{i=1}^{N_k}\rho_iQ_i^k$  the norms
  \begin{align*}
    Q \mapsto \Big(\sum_{i=1}^{N_k}\abs{\rho_i}^{r'}\norm[r']{Q_i^k}^{r'}\Big)^{1/r'}
    \qquad\text{and}
    \qquad Q \mapsto \norm[r']{Q}
  \end{align*}
  are equivalent with constants depending on the shape-regularity of
  $\gridn$ and $\hat\P_\Q$ only. This directly implies the desired estimate.

  Observe that
  $\abs{\Omega_\elm^k}\le c\abs{\elm}$ for all  $\elm\in\gridk$,
  $k\in\N$, with a shape-dependent constant $c>0$; hence,
  $\abs{\Omega^k_{\{\vecE_{k,j}\neq\vecE_k\}}}\le
  c\abs{{\{\vecE_{k,j}\neq\vecE_k\}}}$, and
  it follows from Proposition~\ref{prop:dremlip}\ref{itm:dremlip5} that
   \begin{align}
     \label{est:Psi}
     \beta_r\norm[1,r]{\vec{\Psi}_{k,j}}
     &\leq c\,\norm[r]{\lambda_{k,j}\chi_{\Omega^k_{\{\vecE_{k,j}\neq\vecE_k\}}}
     } \leq c\, 2^{-j/r} \norm[r]{\nabla \vec{E}_k}.
   \end{align}
  Moreover, we have from Proposition~\ref{prop:dremlip} and the continuity properties of $\Bogn$ (see
  Corollary~\ref{c:Dbogovskii}) that
  \begin{subequations}\label{eq:PhiPsi}
    \begin{alignat}{3}
      \vec{\Phi}_{k,j},\vec{\Psi}_{k,j}&\weak 0&\qquad& \text{weakly
        in } W^{1,s}_0(\Omega)^d&\quad&\text{for
        all}~s\in[1,\infty),\label{eq:PhiPsia}
      \\
      \vec{\Phi}_{k,j},\vec{\Psi}_{k,j}&\to 0&\qquad &\text{strongly
        in } L^s(\Omega)^d&\quad&\text{for
        all}~s\in[1,\infty),\label{eq:PhiPsib}
    \end{alignat}
  \end{subequations}
  as $k\to\infty$.

\textbf{Step 2}:
We shall prove (recall the last line of~\eqref{eq:5} for the
definition of $a_k$) that
 \begin{align*}
    \underset{n\to\infty}{\lim\sup}
    \int_{\{\vecE_k=\vecE_{k,j}\}}\abs{a_k}\,\d x\leq c\,2^{-j/r},
  \end{align*}
  with a constant $c>0$ independent of $j$.
  To see this we first observe that $\abs{a_k}=a_k+2a_k^-$
  with the usual notation $a_k^-(x)=\max\{-a_k(x),0\}$, $x\in\Omega$.
  Therefore, we have that
   \begin{align}\label{eq:an=}
    \begin{split}
      \underset{k\to\infty}{\lim\sup}
      \int_{\{\vecE^n=\vecE^{n,j}\}}\abs{a_k}\,\d x&\le\underset{k\to\infty}{\lim\sup}
      \int_{\{\vecE^n=\vecE^{n,j}\}}a_k\,\d x
      \\
      &\qquad+2\,\underset{k\to\infty}{\lim\sup}
      \int_{\{\vecE^n=\vecE^{n,j}\}}a_k^-\,\d x.
    \end{split}
  \end{align}
  By choosing $\phi_k:=\chi_{\supp(a_k^-)}\in L^{\infty}(\Omega)$
  in~\eqref{eq:5}, we observe that the latter term is zero.
  In order to bound the first term, we recall
  \eqref{df:PhiPsi} and observe that
  \begin{align*}
    \int_{\{\vecE_k=\vecE_{k,j}\}}a_k\,\d x
    &=\int_{\{\vecE_k=\vecE_{k,j}\}}
    \left(\bsS_k-\bsS^\star(\cdot,\nablas
      \vecu_\infty)\right): (\nablas \PGdivk\vecu_\infty-\nablas \vecu_\infty)\,\dx
    \\
    &\quad+\int_\Omega
    \bSk: \nablas \vec{\Phi}_{k,j}\,\d x +
    \int_\Omega \bSk: \nablas
    \vec{\Psi}_{k,j}\,\d x
    \\
    &\quad-\int_\Omega \bsS^\star(\cdot,\nablas \vecu_\infty): \nablas
    \vecE_{k,j}\,\d x
    \\
    &\quad+ \int_{\{\vecE_k\neq\vecE_{k,j}\}}
    \big(\bsS^\star(\cdot,\nablas\vecu_\infty)-\bSk\big):
    \nablas \vecE_{k,j}\,\d x
    \\
    &= \text{I}_{k,j}+
    \text{II}_{k,j}+\text{III}_{k,j}+\text{IV}_{k,j}+\text{V}_{k,j}.
  \end{align*}
  Thanks to \eqref{df:Vinfty} and \eqref{eq:bound} we have that
  \begin{align*}
    \abs{\text{I}_{k,j}}&\le \int_{\{\vecE_k=\vecE_{k,j}\}}
    \abs{\bSk(\cdot,\nablas\Uk)-\bsS^\star(\cdot,\nablas
      \vecu_\infty)}\abs{\nablas \PGdivk\vecu_\infty-\nablas \vecu_\infty}\,\dx
    \\
    &\le \norm[r']{\bSk(\cdot,\nablas\Uk)-\bsS^\star(\cdot,\nablas
      \vecu_\infty)}\norm[r]{\nablas \PGdivk\vecu_\infty-\nablas \vecu_\infty}\to 0
  \end{align*}
  as $k \rightarrow \infty$.
  In order to estimate $\text{II}_{k,j}$ we recall that
  $\vec{\Phi}_{k,j}\in\VoG[\gridk]$ is discretely divergence-free, and
  we can therefore use it as a test function in \eqref{eq:Vn0} to deduce that
  \begin{align*}
    \text{II}_{k,j}=-\Trilin{\Uk}{\Uk}{\vec{\Phi}_{k,j}}
    +\int_\Omega\vecf\cdot\vec{\Phi}_{k,j}\dx\to 0\qquad\text{as}~ k\to\infty.
  \end{align*}
  Indeed, the second term vanishes thanks to \eqref{eq:PhiPsia}. The
  first term vanishes thanks to \eqref{eq:strongUn} and the weak convergence \eqref{eq:PhiPsia} of
  $\vec{\Phi}_{k,j}$.
  The term $\text{III}_{k,j}$ can be
  bounded by means of \eqref{est:Psi}; in particular,
   \begin{align*}
     \underset {k\to\infty} {\lim\,\sup}~\abs{\text{III}_{k,j}} \leq
     \underset {k\to\infty} {\lim\,\sup}\,
     \norm[r']{\bsS(\cdot,\nablas
       \Uk)}\norm[r]{\nablas\vec{\Psi}_{k,j}}\le c\, 2^{-j/r},
   \end{align*}
   where we have used \eqref{eq:bound}.
   Proposition~\ref{prop:dremlip} implies that
   \begin{align*}
     \lim_{k\to\infty}\text{IV}_{k,j}=0.
   \end{align*}
   Finally, by \eqref{eq:bound} and Proposition \ref{prop:dremlip}, we have that
   \begin{align*}
     \underset {k\to\infty} {\lim\,\sup}~\abs{\text{V}_{k,j}}&\leq
     \underset {k\to\infty}
     {\lim\,\sup}\,\big(\norm[r']{\bsS^\star(\cdot,\nablas
       \vecu_\infty)}+\norm[r']{\bSk(\cdot,\nablas\Uk)}
     \big)\norm[r]{\nablas
       \vecE_{k,j}\chi_{\{\vecE_k\neq\vecE_{k,j}\}}}
     \\
     &\leq c\,2^{-j/r}.
   \end{align*}
   In view of~\eqref{eq:an=}, this completes Step 2.

   \textbf{Step 3}:
  We prove, for any $\vartheta\in(0,1)$, that
  \begin{align}\label{eq:atheta->0}
    \lim_{k\to\infty}\int_\Omega\abs{a_k}^\vartheta\,\d x=0,
  \end{align}
  which then implies the assertion   $a_k\to 0$ in measure as
  $k\to\infty$.

  Using H{\"o}lder's inequality, we easily obtain that
  \begin{align*}
    \int_\Omega\abs{a_k}^\vartheta\,\d
    x&=\int_{\{\vecE_k=\vecE_{k,j}\}}\abs{a_k}^\vartheta\,\d x
    +\int_{\{\vecE_k\neq\vecE_{k,j}\}}\abs{a_k}^\vartheta\,\d x
    \\
    &\leq
    \abs{\Omega}^{1-\vartheta}
    \left(\int_{\{\vecE_k=\vecE_{k,j}\}}\abs{a_k}\,\d x\right)^{\vartheta}
    + \left(\int_\Omega\abs{a_k}\,\d x\right)^{\vartheta}
    \abs{\{\vecE_k\neq\vecE_{k,j}\}}^{1-\vartheta}.
  \end{align*}
  Thanks to \eqref{eq:bound}, we have that
  $(\int_\Omega\abs{a_k}\,\d x)^{\vartheta}$ is bounded uniformly in $k$
  and by Proposition \ref{prop:dremlip} we have that
  \begin{align*}
    \abs{\{\vecE_k\neq\vecE_{k,j}\}}\leq c\,
    \frac{\norm[1,r]{\vecE_k}^r}{\lambda_{k,j}^r}
    \leq \frac{c}{2^{2^jr}},
  \end{align*}
  where we have used that $\{\vecE_k \}_{k\in\N}$ is
  bounded in $W^{1,r}_0(\Omega)^d$ according to \eqref{eq:bound} and
  Assumption \ref{ass:Pndiv}. Consequently, from Step 2 we deduce that
  \begin{align*}
    \underset {k\to\infty}
    {\lim\,\sup}\int_\Omega\abs{a_k}^\vartheta\,\d x&\le
    c\abs{\Omega}^{1-\vartheta} 2^{-j\vartheta/r} +
    \frac{c}{2^{2^jr(1-\vartheta)}}.
  \end{align*}
  The left-hand side is independent of $j$ and we can thus pass to the
  limit $j\to\infty$. This proves~\eqref{eq:atheta->0}.\qed

\subsection{Proof of Lemma~\ref{l:estGconv}}
\label{ssec:Proof-estGconv}

Since $n_k\to N$ as $k\to\infty$, we may, w.l.o.g., assume that
$n_k=N$ for all $k\in\N$.

\textbf{Step 1:} We shall first prove that, in this case, we have
that the (sub)sequences in Lemma~\ref{l:conv} do actually converge
strongly, i.e.,  \begin{align}
  \begin{alignedat}{2}
    \Uk&\to \vecu_\infty &\qquad&\text{in}~W^{1,t}_0(\Omega)^d,
    \\
    P_k&\to p_\infty&\qquad&\text{in}~L^{\tldt}_0(\Omega),
    \\
    \bSk(\nablas\Uk)=\bsS^N(\nablas\Uk)&\to \bsS_\infty=\bsS^N(\nablas
    \vecu_\infty) &\qquad&\text{in}~L^{\tldt}(\Omega;\Rdds)
    .
  \end{alignedat}\label{eq:strongconv}
\end{align}
To this end, we investigate
  \begin{align*}
    a_k:=\big(\bsS^N(\nablas\Uk)-\bsS^N(\nablas\vecu_\infty)\big):
    \nablas\,(\Uk-\vecu_\infty)\ge0
  \end{align*}
  (compare with Assumption~\ref{ass:Sn}) distinguishing two cases: $r\le\frac{3d}{d+2}$ and $r>\frac{3d}{d+2}$.

  If $r\le\frac{3d}{d+2}$, then  we can deduce,
  as in the proof of
  Lemma~\ref{l:estA} in Section ~\ref{ssec:Proof-estA}, that
  \begin{align}\label{eq:7}
    0\le\int_\Omega|a_k|^\vartheta\dx= \int_\Omega
    \Big(\big(\bsS^N(\nablas\Uk)-\bsS^N(\nablas\vecu_\infty)\big):
    \nablas\,(\Uk-\vecu_\infty)\Big)^\vartheta\dx\to 0,
  \end{align}
  where we have used that  $a_k^-=0$
  almost everywhere in $\Omega$. Thus, recalling that
 $\bsS^N$ is strictly monotone, we obtain that
 \begin{align}\label{eq:8}
   \nablas\Uk&\to \nablas\vecu_\infty\qquad\text{and}\qquad
   \bsS^N(\nablas\Uk)\to \bsS^N(\nablas\vecu_\infty)
   \qquad\text{a.e. in}~ \Omega,
 \end{align}
 at least for a subsequence of $k\to\infty$.
 Since $1\le t<r$ and  $1\le \tldt<r'$ (compare with~\eqref{eq:tsra}), we obtain with Lemma \ref{l:conv} and  Vitali's
 theorem that
 \begin{align*}
   \Uk\to
   \vecu_\infty\quad\text{in}~W^{1,t}_0(\Omega),
   \quad\text{and}\quad
   \bsS^N(\nablas\Uk)\to \bsS^N(\nablas\vecu_\infty)=
   \bsS_\infty\quad\text{in}~L^{\tldt}(\Omega;\Rdds)
   .
 \end{align*}
 Indeed, using H\"older's inequality, we obtain with \eqref{eq:bound}
 that
 \begin{align*}
   \norm[t,\omega]{\nablas \Uk-\nablas\vecu_\infty} \le |\omega|^{\frac{r-t}{rt}} \norm[r,\omega]{\nablas \Uk-\nablas\vecu_\infty}
 \end{align*}
 for all measurable $\omega\subset \Omega$; i.e.,  $\{|\nablas \Uk-\nablas
 \vecu_\infty|^t\}_{k\in\N}$ is uniformly integrable and the claim
 follows from Vitali's theorem and~\eqref{eq:8}. The convergence of
 the stress sequence follows analogously and the claim
 $\bsS^N(\nablas\vecu_\infty)=\bsS_\infty$ is a consequence of the
 uniqueness of the limit.

  If $r>\frac{3d}{d+2}$, then $t=r$ and $\tldt=\tr=r'$; compare
  with~\eqref{eq:tsrb}.  We deduce from Lemma~\ref{l:conv} and
  \eqref{eq:discrete} that
  \begin{align*}
    0&\le \limsup_{k\to\infty}\int_\Omega
    \big(\bsS^N(\nablas\Uk)-\bsS^N(\nablas\vecu_\infty)\big):
    \nablas\,(\Uk-\vecu_\infty)\dx
    \\
    &= \limsup_{k\to\infty}\int_\Omega
    \bsS^N(\nablas\Uk):\nablas\Uk-\bsS_\infty:
    \nablas\vecu_\infty\dx
    \\
    &= \limsup_{k\to\infty}\int_\Omega
    \vecf\cdot\Uk-\bsS_\infty:
    \nablas\vecu_\infty\dx=\int_\Omega\vecf\cdot\vecu_\infty-\bsS_\infty:
    \nablas\vecu_\infty\dx=0.
  \end{align*}
  As before, thanks to the strict monotonicity of
 $\bsS^N$, we have that
 \begin{align}\label{eq:convergences}
   \Uk&\to \vecu_\infty\qquad\text{and}\qquad
   \bsS^N(\nablas\Uk)\to\bsS^N(\nablas\vecu_\infty) \qquad\text{a.e. in}~ \Omega
 \end{align}
 at least for a subsequence of $k\to\infty$.
 Moreover,
   \begin{align*}
     \lim_{k\to\infty}\int_\Omega\bsS^N(\nablas\Uk):\nablas\Uk\dx=
     \int_\Omega\bsS^N(\nablas\vecu_\infty):\nablas\vecu_\infty\dx,
   \end{align*}
   and thanks to \eqref{eq:bound}, we have that
   $\bsS^N(\nablas\Uk):\nablas\Uk$ is bounded in
   $L^1(\Omega)$, hence $\bsS^N(\nablas\Uk):\nablas\Uk\bweak
   \bsS^N(\nablas\vecu_\infty):\nablas\vecu_\infty$; compare with
   Lemma~\ref{l:biting}.
   Recalling Assumption~\ref{ass:Sn} we have $0\le\tilde
   m+\bsS^N(\nablas\Uk):\nablas\Uk$ almost everywhere in $\Omega$.
   Combining these properties, it follows from Lemma~\ref{l:bit=>L1}
   that $\tilde m+\bsS^N(\nablas\Uk):\nablas\Uk\weak
   \tilde m+\bsS^N(\nablas\vecu_\infty):\nablas\vecu_\infty$ in
   $L^1(\Omega)$ and thus  $\bsS^N(\nablas\Uk):\nablas\Uk\weak
   \bsS^N(\nablas\vecu_\infty):\nablas\vecu_\infty$ in
   $L^1(\Omega)$.
   Consequently, by the Dunford--Pettis theorem, $\{\bsS^N(\nablas\Uk):\nablas\Uk\}_{k\in\N}$ is
   uniformly integrable.
   Thanks to the coercivity of $\bsS^N$, we have that $\{|\nablas\Uk|^r\}_{k\in\N}$
   and $\{|\bsS^N(\nablas\Uk)|^{r'}\}_{k\in\N}$ are uniformly integrable
   and hence we deduce from~\eqref{eq:convergences}, 
   with Vitali's theorem, that
   \begin{align*}
     \Uk\to
     \vecu_\infty\quad\text{in}~W^{1,r}_0(\Omega)^d
     \quad\text{and}\quad
     \bsS^N(\nablas\Uk)\to
     \bsS_\infty\quad\text{in}~L^{r'}(\Omega;\Rdds)
     .
   \end{align*}

   It remains to prove the strong convergence of the pressure
   (sub)sequence $\{P_k\}_{k\in\N}$ in $L^{\tldt}_0(\Omega)$. Thanks to
   \eqref{eq:inf-sup}, for $k\in\N$ there exists a $\vecv_k\in W^{1,\tldt'}_0(\Omega)^d$ with
   $\norm[{\tldt'}]{\vecv_k}=1$, such that
   \begin{align*}
     \alpha_{{\tldt'}}\norm[\tldt]{p_\infty-P_k}\le \int_\Omega(p_\infty-P_k)\divo\vecv_k\dx.
   \end{align*}
   Since $\{\vecv_k\}_{k\in\N}$ is bounded in $W^{1,\tldt'}_0(\Omega)^d$,
   there exists a not relabelled weakly converging
   subsequence with weak limit $\vecv\in W^{1,\tldt'}_0(\Omega)^d$.
   Taking $\vecw_k\definedas \vecv_k-\vecv\weak0$ in
   $W^{1,\tldt'}_0(\Omega)^d$ as $k\to\infty$, 
   we deduce using the properties of $\PGdivk$ (see Assumption~\ref{ass:Pndiv}) that
   \begin{align*}
      \alpha_{\tr}\norm[\tldt]{p_\infty-P_k}&\le
                                              \int_\Omega(p_\infty-P_k)\divo\vecw_k\dx - \int_\Omega(p_\infty-P_k)\divo\vecv\dx
      \\
     &=\int_\Omega(p_\infty-P_k)\divo \PGdivk\vecw_k\dx+\int_\Omega
       p_\infty\divo(\vecw_k-\PGdivk\vecw_k)\dx
       \\
     &\quad-\int_\Omega(p_\infty-P_k)\divo\vecv\dx .
   \end{align*}
   The second term vanishes in the limit of $k \rightarrow \infty$ as, thanks to Proposition~\ref{prop:Pnweak},
  $\vecw_k \rightharpoonup 0$ in $W^{1,\tldt'}_0(\Omega)^d$ implies that $\PGdivk \vecw_k
  \rightharpoonup 0$ in $W^{1,\tldt'}_0(\Omega)^d$;
  while the third term vanishes
  because $P_k\weak p_\infty$ weakly in $L^{\tr}_0(\Omega)$ thanks to Lemma~\ref{l:conv}.
  For the first term, we have, thanks to Lemma~\ref{l:conv} and
  \eqref{eq:discrete}, that
  \begin{align*}
    &\int_\Omega(p_\infty-P_k)\divo \PGdivk\vecw_k\dx\\
    &\qquad=\int_\Omega \big(\bsS_\infty-\bsS^N(\nablas
\Uk)\big):\nablas\PGdivk\vecw_k+\big(\frakB[\vecu_\infty,\vecu_\infty]-\frakB[\Uk,\Uk]\big)\cdot
      \PGdivk\vecw_k\dx
      \\
    &\qquad\to 0,
 \end{align*}
 as $k\to\infty$.
   Here we have used in the last step the strong convergence of
   $\{\Uk\}_{k\in\N}$ and
    $\{\bsS^N(\nablas\Uk)\}_{k\in\N}$ in $W^{1,t}_0(\Omega)^d$
    respectively $L^{\tldt}(\Omega;\Rdds)
    $, as well as that the former
   result implies that $ \frakB[\Uk,\Uk]\to
   \frakB[\vecu_\infty,\vecu_\infty]$ strongly in
   $L^{\tldt'}(\Omega)^d$. 
   This completes the proof of~\eqref{eq:strongconv}.


\textbf{Step 2:} As a consequence of~\eqref{eq:strongconv} we shall
prove that
\begin{align}\label{eq:stepResconv}
  \Res\big(\Uk,P_k,\bsS^N(\nablas\Uk)\big)\to
  \Res\big(\vecu_\infty,p_\infty,\bsS^N(\nablas\vecu_\infty)\big)\quad\text{strongly
    in }~W^{-1,\tldt}(\Omega)^d.
\end{align}
To this end we observe, for $\vecv\in W^{1,\tldt'}_0(\Omega)^d$, that
\begin{multline*}
  \dual{\Res^{\textsf{pde}}\big(\Uk,P_k,\bsS^N(\nablas\Uk)\big)-
  \Res^{\textsf{pde}}\big(\vecu_\infty,p_\infty,\bsS_\infty\big)}{\vecv}
  \\
  \begin{aligned}
    &= \int_\Omega\big(\bsS^N(\nablas\Uk)-\bsS_\infty\big):\nablas\vecv+\big(\frakB[\Uk,\Uk]-
    \frakB[\vecu_\infty,\vecu_\infty]\big)\cdot{\vecv}\dx
    \\
    &\quad + \int_\Omega (P_k-p_\infty)\divo \vecv\dx
    \\
    &\le\Big\{\norm[\tldt]{\bsS^N(\nablas\Uk)-\bsS_\infty}+\norm[\tldt]{\frakB[\Uk,\Uk]-
      \frakB[\vecu_\infty,\vecu_\infty]}+\norm[\tldt]{P_k-p_\infty}\Big\}\norm[1,\tldt']{\vecv}.
  \end{aligned}
\end{multline*}
Hence, thanks to~\eqref{eq:tr}, \eqref{eq:strongconv},
and the fact, that $ \frakB[\Uk,\Uk]\to
   \frakB[\vecu_\infty,\vecu_\infty]$ strongly in
   $L^{\tldt'}(\Omega)$ (see Step 1),
   this proves the
   assertion for $\Res^{\textsf{pde}}$.
   The assertion for $\Res^{\textsf{ic}}$ is an
   immediate consequence of \eqref{eq:strongconv}.

\textbf{Step 3:} We use the techniques of \cite{Siebert:11} to prove that
\begin{align}\label{eq:stepRes=0}
  \Res\big(\vecu_\infty,p_\infty,\bsS^N(\nablas\vecu_\infty)\big)=0.
\end{align}
To this end, we first need to recall some results from
\cite{Siebert:11}.

For each $x\in\Omega$, the mesh-size $\hG[\gridk](x)$
is monotonically decreasing and bounded from below by zero; hence,
there exists an $h_\infty\in L^\infty(\Omega)$, such that
\begin{align}
  \label{eq:h_k->h_infty}
  \lim_{k\to\infty}\hG[\gridk]=h_\infty\quad\text{in}~L^\infty(\Omega);
\end{align}
compare e.g. with \cite[Lemma 3.2]{Siebert:11}. We next split the
domain $\Omega$ according to
\begin{align*}
  \grid_k^+\definedas \bigcap_{i\ge
    k}\grid_i=\{\elm\in\gridk:\elm \in\grid_i~\text{for all}~i\ge
  k\}\quad\text{and}\quad
  \grid_k^0:=\gridk\setminus\gridk^+,
\end{align*}
i.e., setting $\Omega_k^+:=\Omega(\gridk^+)$
and $\Omega_k^0:=\Omega(\gridk^0)$, we have
$  \bar\Omega=\Omega^+_k\cup\Omega^0_k$.
It is proved in \cite[Corollary 3.3]{Siebert:11} (compare also~\eqref{eq:h->0})
that, in the limit, the mesh-size
function $\hG$ vanishes on $\Omega_k^0$, i.e.,
\begin{align}\label{eq:Omega^0}
  \lim_{k\to\infty}\norm[\infty]{\hG[\gridk]\chi_{\Omega_k^0}}=0=\lim_{k\to\infty}\norm[\infty]{\hG[\gridk]\chi_{\Omega_k^\star}}.
\end{align}
Here $\Omega^\star_k:=
  \mathcal{U}^\grid(\Omega_k^0)$ and we have
  used the local quasi-uniformity of meshes for the latter limit.
Since $\grid_i^+\subset\grid_k^+\subset\gridk$ for any $k\ge i$, we
have
\begin{align*}
  \Omega_i^0=\Omega(\grid_i^0)=\Omega(\gridk\setminus\grid_i^+).
\end{align*}

Now, fix $\vecv\in W^{2,\tldt'}(\Omega)^d\cap W^{1,\tldt'}_0(\Omega)^d$ and $q\in
W^{1,t'}(\Omega)$, with $\norm[2,\tldt']{\vecv}=1=\norm[1,t']{q}$.
We shall prove that
\begin{align*}
 \dual{\Res\big(\Uk,P_k,\bsS^N(\nablas\Uk)\big)}{(\vecv,q)} &=
 \dual{\Res^{\textsf{pde}}\big(\Uk,P_k,\bsS^N(\nablas\Uk)\big)}{\vecv-\PGdiv[k]\vecv}
 \\
  &\quad+\dual{\Res^{\textsf{ic}}(\Uk)}{q-\PGQ[k]q}
 \end{align*}
vanishes as $k\to\infty$.
Here we use the abbreviations $\PGdiv[k]:=\PGdiv[\grid_k]$ and
$\PGQ[k]:=\PGQ[\gridk]$.
Then,~\eqref{eq:stepRes=0} follows from
\eqref{eq:stepResconv} and the density of $W^{2,\tldt}(\Omega)^d\cap
W^{1,\tldt}_0(\Omega)^d$ in $W^{1,\tldt}_0(\Omega)^d$ and of
$W^{1,t'}(\Omega)$ in $L^{t'}(\Omega)$.
We shall estimate the two terms on the right-hand side separately. For
the first one, we have with
Corollary~\ref{c:bounds} that
 \begin{multline*}
   \dual{\Res^{\textsf{pde}}\big(\Uk,P_k,\bsS^N(\nablas\Uk)\big)}{\vecv-\PGdiv[k]\vecv}
   \\
   \begin{aligned}
     &\Cleq \sum_{\elm\in\gridk}\est_{\gridk}^{\textsf{pde}}\big
     (\Uk,P_k,\bsS^N(\nablas\Uk);\elm\big)^{1/\tldt}
     \norm[\tldt',\mathcal{U}^{\gridk}(\elm)]{\nabla\vecv-\nabla\PGdiv[k]\vecv}
     \\
     &\Cleq \sum_{\elm \in\gridk\setminus\grid_i^+}
     \est_{\gridk}^{\textsf{pde}}\big
     (\Uk,P_k,\bsS^N(\nablas\Uk);\elm\big)^{1/\tldt}
     \norm[\tldt',\mathcal{U}^{\gridk}(E)]{\nabla\vecv-\nabla\PGdiv[k]\vecv}
     \\
     &\quad + \sum_{\elm \in\grid_i^+}
     \est_{\gridk}^{\textsf{pde}}\big
     (\Uk,P_k,\bsS^N(\nablas\Uk);\elm\big)^{1/\tldt}
     \norm[\tldt',\mathcal{U}^{\gridk}(E)]{\nabla\vecv-\nabla\PGdiv[k]\vecv}
     \\
     &\Cleq \est_{\gridk}^{\textsf{pde}}\big
     (\Uk,P_k,\bsS^N(\nablas\Uk);\gridk\setminus\grid_i^+\big)^{1/\tldt}
     \norm[\tldt',\Omega_i^\star]{\nabla\vecv-\nabla\PGdiv[k]\vecv}
     \\
     &\quad + \est_{\gridk}^{\textsf{pde}}\big
     (\Uk,P_k,\bsS^N(\nablas\Uk);\grid_i^+\big)^{1/\tldt}
     \norm[\tldt',\mathcal{U}^\grid(\grid_i^+)]{\nabla\vecv-\nabla\PGdiv[k]\vecv},
   \end{aligned}
 \end{multline*}
where we have used H\"older's inequality and the finite overlapping of the
$\mathcal{U}^{\gridk}(\elm)$, $\elm\in\grid_k$. In view
of Lemma~\ref{l:conv} and Corollary~\ref{c:stabEstG} we obtain that
\begin{align*}
  \est_{\gridk}^{\textsf{pde}}\big(\Uk,P_k,\bsS^N(\nablas\Uk);\gridk\setminus\grid_i^+\big)
   \le\est_{\gridk}^{\textsf{pde}}\big(\Uk,P_k,\bsS^N(\nablas\Uk)\big)\Cleq 1.
\end{align*}
Recalling~\eqref{eq:interpolation}, we thus obtain from the
monotonicity of the mesh-size function that
\begin{multline*}
   \dual{\Res^{\textsf{pde}}\big(\Uk,P_k,\bsS^N(\nablas\Uk)\big)}{\vecv-\PGdiv[k]\vecv}
  \\
  \Cleq
   \norm[\infty]{\hG[\grid_i]\chi_{\Omega_i^\star}}
   + \est_{\gridk}^{\textsf{pde}}\big
   (\Uk,P_k,\bsS^N(\nablas\Uk);\grid_i^+\big)^{1/\tldt} .
\end{multline*}
A similar argument shows that
\begin{align*}
   \dual{\Res^{\textsf{ic}}(\Uk)}{q-\PGQ[k]q}
   \Cleq
   \norm[\infty,\Omega_i^\star]{\hG[\grid_i]}
   + \est_{\gridk}^{\textsf{ic}}\big
   (\Uk;\grid_i^+\big)^{1/t'} .
\end{align*}
Thanks to~\eqref{eq:Omega^0}, for $\epsilon>0$ there exists an $i\in\N$ such that
\begin{align*}
  \dual{\Res\big(\Uk,P_k,\bsS^N(\nablas\Uk)\big)}{(\vecv,q)} \Cleq
  \epsilon + \est_{\gridk}
    \big(\Uk,P_k,\bsS^N(\nablas\Uk);\grid_i^+\big),
\end{align*}
and it therefore remains to prove that
\begin{align}\label{eq:E+conv}
    \est_{\gridk}
    \big(\Uk,P_k,\bsS^N(\nablas\Uk);\grid_i^+\big)
    \to 0\quad\text{as $k\to\infty$}
 \end{align}
in order to deduce~\eqref{eq:stepRes=0}.  To this end, let
\begin{align*}
  \est_k:=\est_{\gridk}
    \big(\Uk,P_k,\bsS^N(\nablas\Uk);\elm_k):=\max\Big\{\est_{\gridk}
    \big(\Uk,P_k,\bsS^N(\nablas\Uk);\elm\big)\colon\elm\in \mathcal{M}_k\Big\}.
\end{align*}
Then, by the stability estimate, Corollary~\ref{c:stabEstG},
and~\eqref{eq:locEb} we have that
\begin{align*}
  \est_k&\Cleq \norm[1,t;\mathcal{U}^\grid(\elm_k)]{\Uk}^{\tldt}+
      \norm[1,t;\mathcal{U}^\grid(\elm_k)]{\Uk}^{2\tldt}+
      \norm[\tldt;\mathcal{U}^\grid(\elm_k)]{P_k}^{\tldt}+
      \norm[\tldt,\mathcal{U}^\grid(\elm_k)]{\vecf}^{\tldt}+
      \norm[r,\mathcal{U}^\grid(\elm_k)]{\tilde  k}^{\tldt}
      \\
      &\quad+\norm[t;\mathcal{U}^\grid(\elm_k)]{\divo \Uk}^{t}
      \\
      &\Cleq
      \norm[1,t;\mathcal{U}^\grid(\elm_k)]{\Uk-\vecu_\infty}^{\tldt}+
      \norm[1,t;\mathcal{U}^\grid(\elm_k)]{\Uk-\vecu_\infty}^{2\tldt}+
      \norm[\tldt;\mathcal{U}^\grid(\elm_k)]{P_k-p_\infty}^{\tldt}
      \\
      &\quad+\norm[t;\mathcal{U}^\grid(\elm_k)]{\divo\Uk}^{t}
      \\
      &\quad
      +\norm[1,t;\mathcal{U}^\grid(\elm_k)]{\vecu_\infty}^{\tldt}+
      \norm[1,t;\mathcal{U}^\grid(\elm_k)]{\vecu_\infty}^{2\tldt}+
      \norm[\tldt;\mathcal{U}^\grid(\elm_k)]{p_\infty}^{\tldt}+
      \norm[\tldt,\mathcal{U}^\grid(\elm_k)]{\vecf}^{\tldt}+
      \norm[t,\mathcal{U}^\grid(\elm_k)]{\tilde  k}^{\tldt}.
\end{align*}
The first line of this bound vanishes thanks
to~\eqref{eq:strongconv}. Since $\elm_k\in\Omega_k^0$, we have that $|\elm_k|^{1/d}\Cleq
\norm[\infty;\Omega_k^0]{\hG[\gridk]}$ and the remaining terms therefore
vanish thanks to ~\eqref{eq:Omega^0} and the
observation that $\elm_k\in\Omega_k^0$.
Therefore, we deduce with~\eqref{eq:marking} that
\begin{align*}
  \est_{\gridk}\big(\Uk,P_k,\bsS^N(\nablas\Uk);\grid_i^+\big)
    &\le \#\grid_i^+
    \max\Big\{\est_{\gridk}\big(\Uk,P_k,\bsS^N(\nablas\Uk);\elm\big)\colon\elm\in
    \grid_i^+\Big\}
    \\
  &\le \#\grid_i^+ g(\est_k)\to 0\quad\text{as}~k\to\infty,
  \end{align*}
where we have used the continuity of $g$ at zero and that $\grid_i^+\subset
\grid_k^+\subset\grid_k\setminus\mathcal{M}_k$.
Combining these observations proves~\eqref{eq:stepRes=0}.


\textbf{Step 4:}
In this step, we shall prove that
\begin{align*}
  \est_{\grid_k}\big(\Uk,P_k,\bsS^N(\nablas\Uk))\to0\qquad\text{as}~k\to\infty.
\end{align*}
To this end, we observe from~\eqref{eq:E+conv} that it suffices to
prove that
\begin{align*}
  \est_{\grid_k}\big(\Uk,P_k,\bsS^N(\nablas\Uk);\gridk\setminus\grid_i^+)\to0\qquad\text{as}~k\to\infty
\end{align*}
for some fixed $i\ge0$. In view of Corollary~\ref{c:stabEstG}, \eqref{eq:stepResconv}
and~\eqref{eq:stepRes=0}, it thus suffices to show that
\begin{align*}
  \osc(\Uk,\bsS^N(\nablas\Uk);\grid_k\setminus\grid_i^+)\to0\qquad\text{as}~k\to\infty.
\end{align*}
This is a consequence of the properties of the
oscillation,~\eqref{eq:Omega^0}, \eqref{eq:strongconv}, and
Assumption~\ref{ass:PGS}, noting that
\begin{align*}
   \norm[\tldt;\Omega_i^0]{\bsS^N(\nablas\Uk)-\PGS[\grid_k] \bsS^N(\nablas\Uk)}&\leq
   \norm[\tldt;\Omega_i^0]{\PGS[\gridk]\bsS_\infty-\PGS[\gridk]\bsS^N(\nablas\Uk)}
   \\
   &\quad+\norm[\tldt;\Omega_i^0]{\PGS[\gridk]\bsS_\infty-\bsS_\infty}
 +\norm[\tldt;
\Omega_i^0]{\bsS_\infty-\bsS^N(\nablas\Uk)}.
\end{align*}
Observing that this readily implies that the estimator vanishes
  on the whole sequence completes the proof.\qed

\section{Graph approximation}
\label{sec:graph-approximation}

In this section we shall discuss the approximation of certain typical maximal monotone
graphs satisfying Assumption~\ref{ass:estAconv}. Admittedly, for particular problems the
approximations suggested here might not always represent the best possible choices, and in
the context of discrete nonlinear solvers, such as Newton's method, properties of
the smoothness of the approximation may become important as well. We believe however
that the following examples provide a reasonable guideline for constructing graph
approximations with properties that are required in applications.


\subsection{Discontinuous stresses}
\label{sec:disc-stress}
\begin{figure}[htbp]
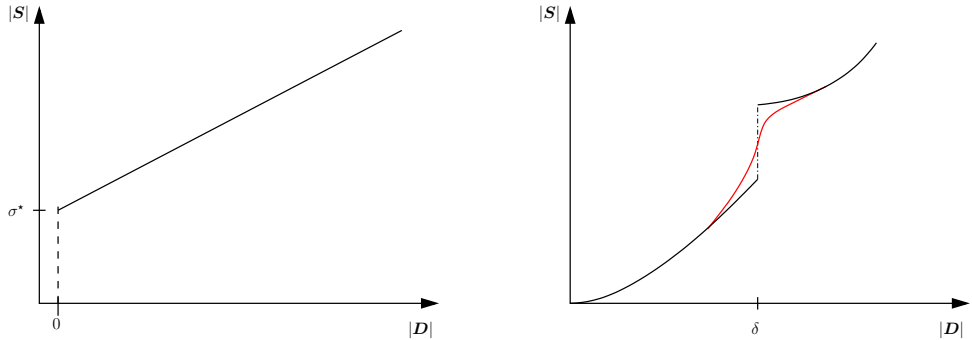

  \resizebox{0.45\textwidth}{!}{\input{approx_bingham.pspdftex}}\hfill
    \resizebox{0.45\textwidth}{!}{\input{approx_jump.pspdftex}}
       \caption{Bingham fluid (left) and schematic approximation of a more
         complex law (right).} \label{F:jump}
\end{figure}

Typical examples of discontinuous dependence of the stress on the
shear rate are Bingham or Herschel--Bulkley fluids. In this case, the fluid
behaves like a rigid body when the shear stress is below a certain critical
value 
and like a Navier--Stokes fluid, respectively power-law fluid,
otherwise; compare with Figure~\ref{F:jump}.
To be more precise, for some yield stress $\sigma\ge0$, we have
\begin{align}\label{eq:bingham}
  \begin{split}
    |\bsS|\le \sigma\quad&\Leftrightarrow \quad\bsD =\boldsymbol{0},\\
    |\bsS|> \sigma\quad&\Leftrightarrow \quad\bsS
    =\sigma\frac{\bsD}{|\bsD|}+2\nu(|\bsD|^2)\bsD;
  \end{split}
\end{align}
where $\nu>0$ denotes the viscosity $\nu>0$; see \cite{DuvautLions:76}.
A selection of the corresponding maximal monotone graph is given, for example, by
\begin{subequations}\label{eq:scalar}
  \begin{align}\label{eq:1dstress}
    \bsS^\star(\bsD)\definedas
    S^\star(\abs{\bsD})\frac{\bsD}{\abs{\bsD}},\qquad\text{with}~
    S^\star(D):=
    \begin{cases}
      0,\quad&\text{if}~D=0
      \\
      \sigma+2\nu(D) D,&\text{otherwise.}
    \end{cases}
  \end{align}
  For the sake of simplicity of presentation, we restrict ourselves in the
  following to $\nu>0$ being a constant. However, we emphasize, that the
  approximation techniques presented below can be generalized to more
  complex relations such as, for example,
  \begin{align}\label{eq:complex}
    S^\star(D)=
    \begin{cases}
      S^\star_1(D),\qquad&\text{if}~D<\delta
      \\
      S^\star_2(D),\qquad&\text{otherwise,}
    \end{cases}\qquad
    S^\star_i(D)=c_i(\kappa_i^2+D^2)^{\frac{q_i-2}2}D,
  \end{align}
\end{subequations}
for $D\ge0$.
Here $\delta\ge0$ and $c_1,c_2,\kappa_1,\kappa_2\ge0$, $q_1>1$, $q_2=r$, such that
$S_1(\delta)\le S_2(\delta)$.

We denote the maximal monotone graph
containing $\{(D,S^\star(D))\colon t\ge0\}$ by $\ao$ and observe that
$(\vec{\delta},\vec{\sigma})\in\Ao$ if and only if
$(|\vec{\delta}|,|\vec{\sigma}|)\in\ao$. Therefore, the approximation of the monotone graph reduces to
approximating the
univariate function $S^\star$ by some smooth $S^n:\R^+_0\to\R^+_0$.
The explicit smooth approximation of $\bsS^\star$ is then obtained
by setting
\begin{align}
  \label{eq:bsSn-Sn}
  \bsS^n(\bsD):=S^n(|\bsD|)\frac\bsD{|\bsD|}\qquad\text{for all}~\bsD\in\Rdds.
\end{align}


\paragraph{A simple approximation} 
A simple approach to approximating $S^\star$ in~\eqref{eq:1dstress} is to
use the following smooth explicit law (cf. \cite{HronMalekStebelTouska:13}):
\begin{align*}
  S^\tau(D):=\bigg(2\nu+\frac{\sigma}{D_\tau}\bigg)D,\qquad\text{where}~
  D_\tau:=\sqrt{D^2+\tau^2}.
\end{align*}
First assume that $\sigma\ge S^\tau(D)$; then, $(0,S^\tau(D))\in \ao$ according
to~\eqref{eq:bingham}. If $D\le\tau$, then
$$|S^\tau(D)-S^\tau(D)|^2+|D-0|^2=|D|^2\le\tau^2.$$ Otherwise, we
have
\begin{align*}
  \sigma\ge \bigg(2\nu+\frac{\sigma}{D_\tau}\bigg)D
  \quad\Leftrightarrow \quad \sigma\tau^2 \ge 2\nu
  D\,D_\tau(D+D_\tau)\ge 4\nu D^3,
\end{align*}
and hence
\begin{align*}
  |S^\tau(D)-S^\tau(D)|^2+|D-0|^2=|D|^2 \le \Big(\frac{\sigma}{4\nu}\Big)^{1/3}\tau^{2/3}.
\end{align*}
Assume now that $0 \leq \sigma<S^\tau(D)$; then, $D\le
\tau$ implies
\begin{align*}
  2\nu
  D+\sigma\frac{D}{D_\tau}=S^\tau(D)>\sigma\quad\Leftrightarrow\quad
  \frac{2\nu}{\sigma}D>\frac{\tau^2}{D_\tau(D+D_\tau)}\ge \frac1{4}.
\end{align*}
In other words this case can occur only for `large' $\tau\ge \frac{\sigma}{8\nu}$.
If $D\ge\tau$ then we have similarly
\begin{align*}
  \sigma \tau^2< 2\nu D\,D_\tau(D+D_\tau)\le 8\nu D^3
  \quad\Leftrightarrow\quad D>\tau^{2/3}
  \Big(\frac{\sigma}{8\nu}\Big)^{1/3}.
\end{align*}
Therefore, we obtain
\begin{align*}
  |S^\star(D)-S^\tau(D)|=\sigma
  \frac{D_\tau-D}{D_\tau}=\frac{\sigma\tau^2}{D_\tau(D+D_\tau)}
  \le \frac{\sigma\tau^2}{D^2}<4\nu^{2/3}\sigma^{1/3}\tau^{2/3}.
\end{align*}
Combining the above cases shows the validity of Assumption~\ref{ass:estAconv} with $\tau=\frac1n$, for
example. The verification of Assumption~\ref{ass:Sn} is left to
the reader.

\paragraph{Approximation by mollification}
We can extend $S^\star$ to an odd function
on the whole real axis by setting
\begin{align*}
  S^\star(D)\definedas -S^\star(-D)\quad\text{for}~ D<0.
\end{align*}
Then, for $n\in\N$, we define an approximation of $S^\star$ by
\begin{align*}
  S^n(t)\definedas \int_{-\infty}^\infty S^\star(s)\,\eta^n(s-t)\ds,
\end{align*}
with $\eta^n(t)=n\eta(n t)$; here $\eta\in C^0(\R)$ is a nonnegative even
function with support
$(-1,1)$ such that $\int_\R\eta(s)\ds=1$. Consequently, the function
$S^n\in C(\R)$ is odd and thus $S^n(0)=0$.

For $D\in\R^+_0$ we have, by the monotonicity of $S^\star$
and the definition of the function $S^n$, that there
exists a $D^\star$ with $0\le D^\star\in
(D-\frac1n,D+\frac1n)$, such that
$(D^\star,S^n(D))\in\ao$. Therefore, we have
\begin{align*}
  |S^n(D)-S^n(D)|^{r'}+|D-D^\star|^{r}\le 0+\frac1{n^r},
\end{align*}
and consequently
\begin{align*}
  \est_{\Ao}(\vec{\delta},\bsS^n(\vec{\delta}))\le
  \frac1{n^r}\to0\quad \text{as}~n\to\infty
\end{align*}
for all $\vec{\delta}\in\Rdds$.
This shows that Assumption~\ref{ass:estAconv} holds. Moreover, $\phi_n$
satisfies Assumption~\ref{ass:Sn}; compare e.g. with
\cite{BulGwiMalSwi:09,GwiaMalSwier:07,GwiazdaZatorska:07}.

\subsection{Monotone graph with plateaus}
\label{sec:const-graph}
\begin{figure}[htbp]
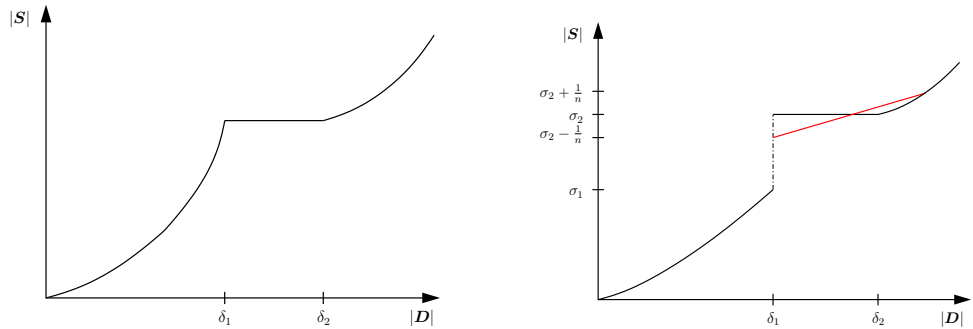
\centering
    \resizebox{0.45\textwidth}{!}{\input{approx_const.pspdftex}}\hfill
    \resizebox{0.45\textwidth}{!}
    {\input{approx_jc.pspdftex}}
     \caption{Graph with plateau (left) and plateau and jump (right).} \label{F:const}
\end{figure}

Similarly to \eqref{eq:scalar}, we consider a maximal monotone graph
with selection $\bsS^\star(\bsD)=S^\star(|\bsD|)\frac{\bsD}{|\bsD|}$,
but now assume that $S^\star:\R^+_0\to \R^+_0$ is continuous with
\begin{align*}
  S^\star(D)&=
  \begin{cases}
    S^\star_1(D),\qquad&\text{if}~D<\delta_1
    \\
    \sigma=const.,\qquad&\text{if}~\delta_1\le D<\delta_2
    \\
    S_2^\star(D),\qquad&\text{else,}
  \end{cases}
  \intertext{with}
  S^\star_i(D)&=c_i(\kappa_i^2+D^2)^{\frac{q_i-2}2}D,\qquad i=1,2.
\end{align*}
Here $c_1,c_2,\kappa_1,\kappa_2\ge0$, $q_1>1$, $q_2=r$, such that
$S^\star_1(\delta_1^\star)=S_2^\star=S_2^\star(\delta_2^\star)$;
compare with Figure~\ref{F:const} (left).
In this case, we are
basically in the same situation as in Section~\ref{sec:disc-stress}
with interchanged roles of $S$ and $D$. Therefore, using the
approximation techniques of Section~\ref{sec:disc-stress}, we can
construct an approximation of the monotone graph where the shear rate
depends explicitly on the shear stress. However, in a practical
numerical method this relation typically has to be inverted, which may cause
additional computational difficulties.

Another approach is to use an
approximation of the form
\begin{align*}
  \tilde S^n(D) :=
  \begin{cases}
    S_1(D)\qquad&\text{if}~S_1(D)<\sigma-\frac1n\\
    S_2(D)\qquad&\text{if}~S_1(D)>\sigma+\frac1n\\
    S_\sigma^n(D)\qquad&\text{otherwise,}
  \end{cases}
\end{align*}
where $S_\sigma^n$ is the linear interpolant between $\sigma-\frac1n$
and $\sigma+\frac1n$ with corresponding values for $D$.

Combined with an approximation strategy as in Section \ref{sec:disc-stress}, this
procedure can also be applied to cases where jumps and plateaus are both
present; compare with Figure~\ref{F:const} (right).

\begin{rmrk}
  The arguments of \S\ref{sec:disc-stress} and \S\ref{sec:const-graph}
  can be obviously extended to finitely many jumps/plateaus
  and even to cases with countably many jumps/plateaus.
\end{rmrk}

\providecommand{\bysame}{\leavevmode\hbox to3em{\hrulefill}\thinspace}
\providecommand{\MR}{\relax\ifhmode\unskip\space\fi MR }
\providecommand{\MRhref}[2]{%
  \href{http://www.ams.org/mathscinet-getitem?mr=#1}{#2}
}
\providecommand{\href}[2]{#2}


\appendix
\section{Stability of the a posteriori Estimator}
\label{sec:appendix}
Referring to Remark~\ref{rem:BerroneSueli} we shall illustrate that
$\divo \bsS^n(\nablas \UnG)$ does not need to belong to
$L^{r'}_{\text{loc}}(\Omega)$ when $1<r<2$.
In fact, this problem already appears
in the a posteriori analysis of quadratic finite element
approximations of the $r$-Laplacian. To see this, we
consider the following problem:
for $f\in L^{r'}(\Omega)$, find $u\in W^{1,r}_0(\Omega)$, such that
\begin{align*}
  -\divo \big(|\nabla u|^{r-2}\nabla u\big) = f \in W^{-1,r'}(\Omega).
\end{align*}
For $\grid\in\grids$, $\ell\ge2$,  let $U_\grid\in \VG=\{V\in C(\bar\Omega):
V|_{\elm}\in \P_\ell~\text{for all}~\elm\in\grid\}$.
Then, for the residual
\begin{align*}
  \Res(U)=f+ \divo(|\nabla U|^{r-2}\nabla U)\in
  W^{-1,r'}(\Omega)
\end{align*}
and $\tr\le r'$, we have the bound
\begin{align*}
  \norm[W^{-1,\tr}(\Omega)]{\divo |\nabla U_\grid|^{r-2}\nabla
    U_\grid+f}^{\tr}&\le C \sum_{\elm\in\grid}\norm[\tr,\elm]{\hG(-\divo |\nabla U_\grid|^{r-2}\nabla
    U_\grid+f)}^{\tr}
  \\
  &\quad+ \sum_{\elm\in\grid}\norm[\tr,\partial\elm]{\jump{\hG^{1/r'}|\nabla U_\grid|^{r-2}\nabla
      U_\grid}}^{\tr}
  \\
  &:=\sum_{\elm\in\grid}\est(U,\grid,\elm)^{\tr}=:\est(U,\grid)^{\tr};
\end{align*}
compare e.g. with \cite{BerroneSuli:08}.
In fact, it can be shown as in \cite[Chapter 3.8]{Verfuerth:13}, that
\begin{align*}
  \normbig[-1,\tr]{\Res(U)}\Cleq \est(U,\grid) \Cleq
  \normbig[-1,\tr]{\Res(U)} + \osc(U,\grid).
\end{align*}
However, it is not clear, whether the element residual
$$\sum_{\elm\in\grid}\norm[\tr,\elm]{\hG(-\divo |\nabla U_\grid|^{r-2}\nabla
  U_\grid+f)}^{\tr}$$
(and thus  $\est(U,\grid)$ and $\osc(U,\grid)$) is bounded.  In fact,
elementary calculations show that on an element $\elm\in\grid$,
we have
\newcommand{\eins}{\ensuremath{\mathbbm{1}}}
\begin{multline*}
  \divo(|\nabla U|^{r-2}\nabla U)
  \\
  =\sum_{i=1}^d\Big[\abs{\nabla
    U}^{r-2}\Big(\eins-\frac{\nabla U}{|\nabla U|^2}(\nabla U)_i\Big)+
  (r-2)|\nabla U|^{r-2}\frac{\nabla U}{|\nabla U|^2}(\nabla
  U)_i\Big] \frac{\partial^2 U}{\partial x_i^2}.
\end{multline*}
As a consequence of $U|_\elm\in\P_\ell$, we have that
$\norm[\infty,\elm]{\frac{\partial^2 U}{(\partial
    x_i)^2}}$ is bounded. Since
$|\nabla U|$ and $\frac{\partial^2 U}{\partial x_i^2}$ may
have different zeros, we need to check when
\begin{align*}
  \int_\elm |\nabla U|^{(r-2)\tr}\dx<\infty.
\end{align*}
Since $\nabla U\in\P_{\ell-1}$, for this to hold we need to have that
\begin{align}\label{eq:ell}
  \begin{split}
    (\ell -1) (r-2)\tr&>-1\qquad
    \Leftrightarrow \qquad\ell < 1+\frac{1}{\tr(2-r)}
  \end{split}
\end{align}
where we have used that $1<r<2$.
This puts a restriction on the polynomial degree of the finite
element space.

We consider some special cases:
\begin{itemize}
\item $d=1$ and $\ell =2$, $1<r<2$ and $\tr=r'$. For simplicity, we assume that
  $\elm=[0,1]$ and that $U=x^2$. Thus $\nabla U= U'= 2x$, $\divo \nabla
  U = U'' =2$, and we have that
  \begin{align*}
    \int_{0}^1 \big|\divo (|\nabla U|^{r-2}\nabla U )\big|^{r'}\dx&=
     \int_0^1 \big|\divo (|2x|^{r-2}2x )\big|^{r'}\dx=
     \\
     &=2^{(r-1)r'}\int_0^1 \big|\divo (x^{r-1} )\big|^{r'}\dx
     \\
     &=2^{(r-1)r'}(r-1)^{r'}\int_0^1 x^{(r-2)r'}\dx.
  \end{align*}
  This term is finite if and only if
  \begin{align*}
    &\qquad\quad(r-2)r'>-1
    \\
    &\Leftrightarrow\quad
    r(r-2)>1-r
    \\
    &\Leftrightarrow\quad
    r^2 -r-1>0
    \\
    &\Leftrightarrow\quad
    r\not\in[\frac{1-\sqrt{5}}2,\frac{1+\sqrt{5}}2].
  \end{align*}
  This restriction applies to any dimension, since one could
  always choose $U=[x_1^2,0,\ldots,0]^t$ as a quadratic function on the
  standard $d$-simplex.

\item $1<r<d=2$ and $\tr = \frac12
  \frac{dr}{d-r}$. This case is relevant for power-law type models
  when $r$ is close to $1$, or, more precisely,
  when $r\le \frac{3d}{d+2}$ respectively $r'\ge \frac12
  \frac{dr}{d-r}$; compare with \eqref{eq:tr}. Then,
  \eqref{eq:ell} becomes
  \begin{align*}
    \ell < 1 + 2 \frac{d-r}{d r(2-r)}= 1+\frac1r.
  \end{align*}
  Since $r>1$, in this case the residual is in general not in
  $L^{\tr}_{\text{loc}}(\Omega)$
  even for quadratic elements ($\ell=2$).
\end{itemize}

\end{document}